\documentclass[12pt]{amsart}

\usepackage{mathabx}
\usepackage{multicol}
\usepackage[inline]{asymptote}
\usepackage[hidelinks]{hyperref}
\usepackage[shortlabels]{enumitem}
\usepackage[letterpaper, width=165mm, top=25mm,bottom=25mm,bindingoffset=6mm]{geometry}
\usepackage{url}

\theoremstyle{plain}
\newtheorem{theorem}{Theorem}
\numberwithin{theorem}{section}
\newtheorem{lemma}[theorem]{Lemma}
\newtheorem{fgl}[theorem]{Factor Group Lemma}
\newtheorem{corollary}[theorem]{Corollary}
\newtheorem{proposition}[theorem]{Proposition}
\newtheorem{caseprop}{Proposition}
\numberwithin{caseprop}{subsection}

\theoremstyle{definition}
\newtheorem{definition}[theorem]{Definition}
\newtheorem{assumption}[theorem]{Assumption}
\newtheorem{notation}[theorem]{Notation}

\newtheorem{case}{Case}
\numberwithin{case}{section}
\newtheorem{subcase}{Subcase}
\numberwithin{subcase}{case}

\newcommand{\cartprod}{\mathbin{\square}}

\DeclareMathOperator{\Cay}{Cay}

\allowdisplaybreaks

\newcommand{\gq}{ a_{q} }
\newcommand{\gt}{ a_{3} }

\makeatletter

\renewcommand{\tocsection}[3]{%
  \indentlabel{\@ifnotempty{#2}{\bfseries\ignorespaces#1 #2\quad}}\bfseries#3}
\renewcommand{\tocsubsection}[3]{%
  \indentlabel{\@ifnotempty{#2}{\ignorespaces#1 #2\quad}}#3}

\newcommand\@dotsep{4.5}
\def\@tocline#1#2#3#4#5#6#7{\relax
  \ifnum #1>\c@tocdepth 
  \else
    \par \addpenalty\@secpenalty\addvspace{#2}%
    \begingroup \hyphenpenalty\@M
    \@ifempty{#4}{%
      \@tempdima\csname r@tocindent\number#1\endcsname\relax
    }{%
      \@tempdima#4\relax
    }%
    \parindent\z@ \leftskip#3\relax \advance\leftskip\@tempdima\relax
    \rightskip\@pnumwidth plus1em \parfillskip-\@pnumwidth
    #5\leavevmode\hskip-\@tempdima{#6}\nobreak
    \leaders\hbox{$\m@th\mkern \@dotsep mu\hbox{.}\mkern \@dotsep mu$}\hfill
    \nobreak
    \hbox to\@pnumwidth{\@tocpagenum{\ifnum#1=1\bfseries\fi#7}}\par
    \nobreak
    \endgroup
  \fi}
\AtBeginDocument{%
\expandafter\renewcommand\csname r@tocindent0\endcsname{0pt}
}
\def\l@subsection{\@tocline{2}{0pt}{2.5pc}{5pc}{}}
\makeatother

\setlistdepth{5}
\newlist{myEnumerate}{enumerate}{5}

\begin{document}

\title[Cayley Graphs of Order $6pq$  are Hamiltonian]{Cayley Graphs of Order $6pq$  are Hamiltonian}
\author[F. Maghsoudi]{Farzad Maghsoudi}
\address{Department of Mathematics and Computer Science, University of Lethbridge, Lethbridge, Alberta, T1K 3M4, Canada}
\email{farzad.maghsoudi@uleth.ca}
\subjclass[2010]{05C25, 05C45
}
\keywords{Cayley graphs, Hamiltonian cycles}
\begin{abstract}
    Assume $ G $ is a finite group, such that $ |G| = 6pq $ or $ 7pq $, where $ p $ and $ q $ are distinct prime numbers, and let $ S $ be a generating set of $ G $. We prove there is a Hamiltonian cycle in the corresponding Cayley graph $ \Cay(G;S) $.
\end{abstract}

{\mathversion{bold} \maketitle}
\tableofcontents

\section[Introduction]{Introduction}
\label{ch:Introduction}
All graphs in this paper are undirected.
\begin{definition}[{cf.~\cite[p.~34]{Godsil}}]
Let $ S $ be a subset of a finite group $ G $. The \textit{Cayley graph} $ \Cay(G;S) $ is the graph whose vertices are elements of $ G $, with an edge joining $ g $ and $ gs $, for every $ g\in G $ and $ s\in S $.
\end{definition}

There have been many papers on the topic of Hamiltonian cycles in Cayley graphs, but it is still an open question whether every connected Cayley graph has a Hamiltonian cycle. (See survey papers~\cite{Curran,Twelve,Pak} for more information. We ignore the trivial counterexamples on $1$ or $2$ vertices.) The following result combines the main result of this paper with the previous work of several authors~(C.~C.~Chen and N.~Quimpo \cite{Forth}, S.~J.~Curran, J.~Morris and D.~W.~Morris \cite{Fifth}, E.~Ghaderpour and D.~W.~Morris \cite{Seventh,Sixth}, D.~Jungreis and E.~Friedman \cite{twenty}, Kutnar et al. \cite{Tenth}, K.~Keating and D.~Witte \cite{Ninth}, D.~Li \cite{thirty}, D.~W.~Morris and K.~Wilk \cite{Fourteen}, and D.~Witte \cite{Eighteen}).

\begin{theorem}[\cite{Tenth,Fourteen,Eighteen}] \label{theorem 1.2}
Let $ G $ be a finite group. Every connected Cayley graph on $ G $ has a Hamiltonian cycle if $ |G| $ has any of the following forms (where $ p $, $ q $, and $ r $ are distinct primes):
  \begin{enumerate}
      \item $ kp $, where $ 1 \leq k \leq 47 $,\label{theorem 1.2.1}
      \item $ kpq $, where $ 1 \leq k \leq 7 $,\label{theorem 1.2.2}
      \item $ pqr $, \label{theorem 1.2.3}
      \item $ kp^{2} $, where $ 1 \leq k \leq 4 $,
      \item $ kp^{3} $, where $ 1\leq k \leq 2 $,
      \item $ p^k $, where $1\leq k<\infty$. \label{theorem 1.2.6}
  \end{enumerate}
\end{theorem}

The new part of Theorem~\ref{theorem 1.2} is the upper bound in~(\ref{theorem 1.2.2}). Previously, that part of the result was only known for $1 \leq k \leq 5$, but we improve this condition: we show that $5$ can be replaced with $7$. The hard part is when $k = 6$:

 \begin{theorem} \label{theorem1.1}
Assume $ G $ is a finite group of order $ 6pq $, where $ p $ and $ q $ are distinct prime numbers. Then every connected Cayley graph on $ G $ contains a Hamiltonian cycle.
\end{theorem}

This generalizes \cite{Sixth}, which considered only the case where $q = 5$. The proof takes up all of Section~\ref{ch:6}, after some preliminaries in Section~\ref{ch:Preliminaries}.

Unlike~Theorem~\ref{theorem1.1}, the following observation follows easily from known results, and may be known to experts. The proof is on page~\pageref{7pq}.

\begin{proposition} \label{pro 1.1.4}
Assume $ G $ is a finite group of order $ 7pq $, where $ p $ and $ q $ are distinct prime numbers. Then every connected Cayley graph on $ G $ contains a Hamiltonian cycle.
\end{proposition}

The Introduction of the author's masters thesis \cite{Farzad} provides additional background and a description of the methods that are used in the proof of the main theorem.

\numberwithin{theorem}{subsection}

\section{Preliminaries}\label{ch:Preliminaries}

This section establishes basic terminology and notation, and proves a number of technical results that will be used in the proof of Theorem~\ref{theorem1.1}.  In particular, it is shown we may assume that $ |G| $ is square-free, so the Sylow subgroups of $ G $ are $ \mathcal{C}_2 $, $ \mathcal{C}_3 $, $ \mathcal{C}_p $, and $ \mathcal{C}_q $, and that $ |G'| $ has precisely 2 prime factors, so $ G' $ is either $ \mathcal{C}_p\times\mathcal{C}_q $ or $ \mathcal{C}_3\times\mathcal{C}_p $. 

\subsection{Basic notation and definitions}

Throughout the paper, we have used standard terminology of graph theory and group theory that can be found in textbooks, such as \cite{Godsil,Hall}.

The following notation is used through the paper:
\begin{itemize}
   
    \item The commutator of $ g $ and $ h $ is denoted by $ [g,h] = ghg^{-1}h^{-1} $.
    \item We will always let $G' = [G,G]$ be the commutator subgroup of $G$.
    \item We define $ \overline{G} = G/G' $, $ \overline{g} = gG' $ for any $ g \in G $, and $ \overline{S} = \{\overline{g};g\in S\} $ for any $ S \subseteq G $.
    \item $ C_{G'}(S) $ denotes the centralizer of $ S $ in $ G' $. 
    \item $ G \ltimes H $ denotes a semidirect product of groups $ G $ and $ H $.
    \item $ D_{2n} $ denotes the dihedral group of order $ 2n $.
    \item $ e $ denotes the identity element of $ G $.
    \item For $ S\subseteq G $, a sequence $ (s_{1}, s_{2},\ldots,s_{n}) $ of elements of $ S\cup S^{-1} $ specifies the walk in the Cayley graph $ \Cay(G;S) $ that visits the vertices: $ e, s_{1}, s_{1}s_{2},\ldots, s_{1}s_{2}\cdots s_{n} $. Also, $ (s_{1}, s_{2},\ldots,s_{n})^{-1} = (s_{n}^{-1}, s_{n-1}^{-1},\ldots,s_{1}^{-1}) $.
    \item We use $ (\overline{s_{1}}, \overline{s_{2}},\ldots,\overline{s_{n}}) $ to denote the image of this walk in the quotient $ \Cay(G/G';\overline{S}) \\= \Cay(\overline{G};\overline{S}) $.
    \item For $ k\in\mathbb{Z^{+} } $, we use $ (s_{1}, s_{2},\ldots,s_{m})^{k} $ to denote the concatenation of $ k $ copies of the sequence $ (s_{1}, s_{2},\ldots,s_{m}) $.
    \item $ p $ and $ q $ are distinct prime numbers.
     \item $\mathcal{C}_n$ denotes the cyclic group of order $n$. When $|G| = 6pq$ and it is square free~(as is usually the case in Section~\ref{ch:6}), the Sylow subgroups are $\mathcal{C}_2$, $\mathcal{C}_3$, $\mathcal{C}_p$, and $\mathcal{C}_q$. Also, the commutator subgroup $G'$ will usually be either $\mathcal{C}_p\times\mathcal{C}_q$ or $\mathcal{C}_3\times\mathcal{C}_p$, so $\mathcal{C}_p$ is a normal subgroup and either $\mathcal{C}_q$ or $\mathcal{C}_3$ is also a normal subgroup.
     \item $\widehat{G} = G/\mathcal{C}_p$, we also let $\widecheck{G} = G/\mathcal{C}_q$ when $\mathcal{C}_q$ is a normal subgroup, and let $\overleftrightarrow{G} = G/\mathcal{C}_3$ when $\mathcal{C}_3$ is a normal subgroup.
    \item We let $ a_2 $, $ \gt $, $ \gamma_p $, and $ \gq $ be elements of $ G $ that generate $ \mathcal{C}_2 $, $ \mathcal{C}_3 $, $ \mathcal{C}_p $, and $ \mathcal{C}_q $, respectively.
\end{itemize}

\subsection{Basic methods}\label{sec 1.2}

In this subsection we explain some of the key ideas in the proof of our main result~(Theorem~\ref{theorem1.1}).

It is easy to see that $\Cay(G;S)$ is connected if and only if $S$ generates $G$~(\cite[Lemma 3.7.4]{Godsil}). Also, if $S$ is a subset of $S_0$, then $\Cay(G;S)$ is a subgraph of $\Cay(G;S_0)$ that contains all of the vertices. Therefore, in order to show that every connected Cayley graph on $G$ contains a Hamiltonian cycle, it suffices to consider $\Cay(G;S)$, where $S$ is a generating set that is \emph{minimal}, which means that no proper subset of $S$ generates $G$.

The following well known (and easy) result handles the case of Theorem~\ref{theorem1.1} where $G$ is abelian. 

\begin{lemma}[{\cite[Corollary on page 257]{Quimpo}}] \label{abelain group}
Assume $ G $ is an abelian group. Then every connected Cayley graph on $ G $ has a Hamiltonian cycle.
\end{lemma}

Note $ \Cay(\mathcal{C}_2;\{a\}) $ is a Cayley graph with two vertices, where $ \mathcal{C}_2 = \langle a \rangle $. We consider $ (a,a) $ as its Hamiltonian cycle which is:
\begin{align*}
    e \stackrel{a}{\rightarrow} a \stackrel{a}{\rightarrow} a^2 = e.
\end{align*}

Although graph theorists would not typically consider this a cycle, it satisfies the basic property of visiting each vertex exactly once. In some of our inductive proofs, we require a Hamiltonian cycle in a Cayley graph on a quotient group. When this quotient group is $ \mathcal{C}_2 $, this Hamiltonian cycle provide the structure we need for our inductive arguments to work.

\begin{theorem}[Maru\v{s}i\v{c}~\cite{Marusic}, Durnberger~\cite{Durnberger1,Durnberger2}, and Keating-Witte~{\cite{Ninth}}] \label{Keating-Witte}
If the commutator subgroup $ G' $ of $ G $ is a cyclic $p$-group, then every connected Cayley graph on $G$ has a Hamiltonian cycle.
\end{theorem}

\begin{theorem}[Chen-Quimpo~{\cite{Third}}] \label{Chen-Quimpo}
Let $v$ and $w$ be two distinct vertices of a connected Cayley graph $\Cay(G;S)$. Assume $G$ is abelian, $|G|$ is odd, and the valency of $\Cay(G;S)$ is at least $3$. Then $\Cay(G;S)$ has a Hamiltonian path that starts at $v$ and ends at $w$.
\end{theorem}

The following lemma (and its corollary) often provide a way to lift this Hamiltonian cycle to a Hamiltonian cycle in $\Cay(G;S)$. Before stating the results, we introduce a useful piece of notation.

\begin{notation}
Suppose $N$ is a normal subgroup of $G$, and $ C = (s_1,s_2,\ldots,s_n) $ is a walk in $\Cay(G;S)$. If the walk $ (s_1 N, s_2 N, \ldots,s_n N)$ in $\Cay(G/N;SN/N)$ is closed, then its \textit{voltage} is the product $\mathbb{V}(C) = s_1 s_2 \cdots s_n $. This is an element of $N$. In particular, if $C = (\overline{s}_1,\overline{s}_2,\ldots,\overline{s}_n)$ is a Hamiltonian cycle in $\Cay(\overline{G},\overline{S})$, then $ \mathbb{V}(C) = s_1 s_2\cdots s_n $.
\end{notation}

\begin{fgl}[{\cite[Section 2.2]{Twelve}}] \label{FGL}
Suppose:
\begin{itemize}
    \item $ S $ is a generating set of $ G $,
    \item $ N $ is a cyclic normal subgroup of $ G $,
    \item $ \overline{G} = G/N $,
    \item $ C =(\overline{s_{1}},\overline{s_{2}},\ldots,\overline{s_{n}})$ is a Hamiltonian cycle in $ \Cay(G/N;\overline{S}) $, and
    \item the voltage $ \mathbb{V}(C) $  generates $ N $.
\end{itemize}
Then there is a Hamiltonian cycle in $ \Cay(G;S) $.
\end{fgl}

\begin{corollary}[{\cite[Corollary 2.3]{Sixth}}] \label{corollary 5.2}
Suppose:
\begin{itemize}
    \item $ S $ is a generating set of $ G $,
    \item $ N $ is a normal subgroup of $ G $, such that $ |N| $ is prime,
    \item $ sN = tN $ for some $ s,t\in S $ with $ s \neq t $, and
    \item there is a Hamiltonian cycle in $ \Cay(G/N;\overline{S}) $ that uses at least one edge labeled~$ \overline{s} $.
\end{itemize}
Then there is a Hamiltonian cycle in $ \Cay(G;S) $.
\end{corollary}

\begin{lemma} \label{lemma 2.5.2}
Assume $G = H\ltimes(\mathcal{C}_p\times\mathcal{C}_q)$, where $ G' = \mathcal{C}_p\times\mathcal{C}_q $, and let $S$ be a generating set of $G$. As usual, let $\overline{G} = G/G' \cong H$. Assume there is a unique element $c$ of $S$ that is not in $H\ltimes\mathcal{C}_q $, and $ C $ is a Hamiltonian cycle in $ \Cay(\overline{G};\overline{S}) $ such that $c$ occurs precisely once in $C$. Then the subgroup generated by~$ \mathbb{V}(C) $~contains~$ \mathcal{C}_p $. 
\end{lemma}

\begin{proof}
Write $C = (\overline{s}_1,\overline{s}_2,\cdots,\overline{s}_n)$, and let $H^+ = H \ltimes\mathcal{C}_q$. By assumption, there is a unique $k$, such that $s_k = c$, and all other elements of $S$ are in $H^+$. Therefore,
\begin{align*}
 \mathbb{V}(C) = s_1s_2...s_n \in H^+\cdot H^+\cdots H^+\cdot c\cdot H^+\cdot H^+\cdots H^+ = H^+cH^+.    
\end{align*}
Since $c \notin H^+$, we conclude that $\mathbb{V}(C) \notin H^+$.

On the other hand, since $\mathbb{V}(C)$ is an element of $G' = \mathcal{C}_p\times\mathcal{C}_q$, we have $\mathbb{V}(C) = \gq^i\gamma_p^j \in H^+\gamma_p^j$.  Since $\mathbb{V}(C) \notin H^+$, this implies $j \not \equiv 0 \pmod{p}$, so $\langle \gq^i\gamma_p^j\rangle $ contains~$\mathcal{C}_p$.
\end{proof}

\begin{definition}
The \textit{Cartesian product} $ X_1 \cartprod X_2 $ of graphs $ X_1 $ and $ X_2 $ is a graph such that the vertex set of $ X_1 \cartprod X_2 $ is $ V(X_1)\times V(X_2) = \{(v,v'); v\in V(X_1), v'\in V(X_2) \} $, and two vertices $ (v_1,v_2) $ and $ (v_1',v_2') $ are adjacent in $ X_1 \cartprod X_2 $ if and only if either
\begin{itemize}
\item $ v_1 = v_1' $ and $v_2$ is adjacent to $v_2'$ in $ X_2 $ or
\item $ v_2 = v_2' $ and $v_1$ is adjacent to $v_1'$ in $ X_1 $.
\end{itemize}
\end{definition}

\begin{lemma}[{\cite[Lemma 5 on page 28]{Third}}] \label{lemma 2.3.13}
The Cartesian product of a path and a cycle is Hamiltonian. 
\end{lemma}

\begin{corollary}[cf.~{\cite[Corollary on page 29]{Third}}]\label{Cartesian}
The Cartesian product of two Hamiltonian graphs is Hamiltonian.
\end{corollary}

\begin{lemma}[{}{\cite[Lemma 2.27]{Tenth}}] \label{lemma 5.6} 
Let $ S $ generate the finite group $ G $, and let $ s \in S $, such that $ \langle s \rangle \triangleleft  G $. If $ \Cay(G/\langle s \rangle ;\overline{S}) $ has a Hamiltonian cycle, and either 
\begin{enumerate}
    \item $ s \in Z(G) $ \label{lemma 5.6 n1}, or
    \item $ Z(G)$  $ \cap$ $\langle s \rangle = \{e\}$, \label{lemma 5.6 n2}
\end{enumerate}
then $ \Cay(G;S) $ has a Hamiltonian cycle.
\end{lemma}

\subsection{Some facts from group theory}

In this subsection we state some facts in group theory, which are used to prove our main result. The following lemma often makes it possible to use Factor Group Lemma~\ref{FGL} for ﬁnding Hamiltonian cycles in connected Cayley graphs of $G$.

\begin{lemma}[{}{\cite[Corollary 4.4]{Fifth}}] \label{lemma 5.12}
Assume $ G = \langle a,b \rangle $ and $G'$ is cyclic. Then~$ G' = \langle [a,b] \rangle $.
\end{lemma}

\begin{corollary} \label{lemma 5.13.}
Assume $ G = \langle a,b \rangle $ and $ \gcd(k,|a|) = 1 $, where $ k \in \mathbb{Z} $, and $G'$ is cyclic. Then $ G' = \langle [a^k,b] \rangle $.
\end{corollary}

\begin{lemma} \label{lemma 5.13} Assume $ G = (\mathcal{C}_p\times\mathcal{C}_q)\ltimes (\mathcal{C}_r\times\mathcal{C}_t) $, where $ p,q,r $ and $ t $ are distinct primes. If $ |\overline{a}| = pq $, then $ |a| = pq $.
\end{lemma}
\begin{proof}
Suppose $ |a| \neq pq $. Without loss of generality, assume $ |a| $ is divisible by $ r $. Then (after replacing $a$ by a conjugate) the abelian group $ \langle a \rangle $ contains $\mathcal{C}_p\times\mathcal{C}_q$ and $\mathcal{C}_r$, so $\mathcal{C}_r$ centralizes $ \mathcal{C}_p\times\mathcal{C}_q $.  Since $\mathcal{C}_r$ also centralizes $\mathcal{C}_t$, this implies that $\mathcal{C}_r\subseteq Z(G)$. This contradicts the fact that $G' \cap Z(G) = \{e\}$~(see Proposition~\ref{Hall Theorem}(\ref{Hall Theorem2})). 
\end{proof}

\begin{lemma}[{\cite[Exercise 19 on page 4]{Robinson}}] \label{lemma 1.5.6}
Assume $ |G| = 2k $, where $ k $ is odd. Then $G$ has a subgroup of index~$ 2 $.
\end{lemma}

\begin{corollary} \label{odd order commutator}
Assume $ |G| = 2k $, where $ k $ is odd. Then $ |G'| $ is odd.
\end{corollary}

\begin{proof}
By Lemma~\ref{lemma 1.5.6}, there is a normal subgroup $ H $ of $ G $ such that $ [G:H] = 2 $. Now since $ G/H $ has order $ 2 $, then $ G/H $ is abelian, so $ G' \subseteq H $. Therefore, $ |G'| $ is~odd.
\end{proof}

\begin{proposition}[{\cite[Theorem 9.4.3 on page 146]{Hall}}, cf.~{\cite[Lemma 2.11]{Sixth}}] \label{Hall Theorem}
Assume $|G|$ is square-free. Then: 
\begin{enumerate}
    \item $ G' $ and $ G/G' $ are cyclic, \label{Hall Theorem1}
    \item $  Z(G) \cap G' = \{e\} $,\label{Hall Theorem2}
    \item $ G \cong C_n \ltimes G' $, for some $ n \in \mathbb{Z}^+ $,\label{Hall Theorem3}
    \item If $b$ and $\gamma$ are elements of $G$ such that $\langle bG' \rangle = G/G'$ and $ \langle \gamma \rangle = G'$, then $ \langle b,\gamma \rangle = G$, and there are integers $m$, $n$, and $\tau$, such that $ |\gamma| = m $, $ |b| = n $, $ b\gamma b^{-1} = \gamma^{\tau} $, $ mn = |G| $, $ \gcd(\tau-1,m) = 1 $, and $ \tau^n \equiv 1 \pmod{m} $.\label{Hall Theorem4}
\end{enumerate}
\end{proposition}

\begin{notation}
For $ \tau $ as defined in Proposition~\ref{Hall Theorem}(\ref{Hall Theorem4}), we use $ \tau^{-1} $ to denote the inverse of $ \tau $ modulo $ m $~(so $ \tau^{-1} \equiv \tau^{n-1} \pmod{m} $).
\end{notation}

\subsection{Cayley graphs that contain a Hamiltonian cycle}

In this subsection we show that there exists a Hamiltonian cycle in some special connected Cayley graphs. The following proposition shows that in our proof of Theorem~\ref{theorem1.1} we can assume $ |G| $ is square-free, since the cases where $ |G| $ is not square-free have already been dealt with. At the end of this subsection we prove the Proposition~\ref{pro 1.1.4}.

\begin{proposition}\label{proposition 5.9}
Assume:
\begin{itemize}
    \item $ |G| = 6pq $, where $ p $ and $ q $ are distinct prime numbers, and
    \item $ |G| $ is not square-free (i.e. $ \{p,q\} \cap \{2,3\} \neq \emptyset $).
\end{itemize}
Then every connected Cayley graph on $ G $ has a Hamiltonian cycle.
\end{proposition}

\begin{proof}
Without loss of generality we may assume $ q \in \{2,3\} $. Then $ |G| \in \{12p,18p\} $. Therefore, Theorem~\ref{theorem 1.2}(\ref{theorem 1.2.1}) applies.
\end{proof}

\begin{proposition}[{\cite[Proposition 5.5]{Thirteen}}]\label{prop 1.5.2} 
If $ n $ is divisible by at most $3$ distinct primes, then every Cayley diagram in $ D_{2n} $ has a Hamiltonian cycle. 
\end{proposition}

The following proposition demonstrates that we can assume $ |G'| $ in Theorem~\ref{theorem1.1} is a product of two distinct prime numbers.

\begin{proposition} \label{proposition 5.4}
Assume $ |G| = 2pqr $, where $ p $, $ q $ and $ r $ are distinct odd prime numbers. Now if $ |G'| \in \{1,pqr\} $ or $ |G'| $ is prime, then every connected Cayley graph on $ G $ has a Hamiltonian cycle.
\end{proposition}
\begin{proof}
If $ |G'| = 1 $, then $ G' = \{e\} $. So $ G $ is an abelian group. Therefore, Lemma~\ref{abelain group} applies. Now if $|G'|$ is prime, then Theorem~\ref{Keating-Witte} applies. Finally, if $ |G'| = pqr $, then
\begin{align*}
G = \mathcal{C}_2\ltimes(\mathcal{C}_p\times\mathcal{C}_q\times\mathcal{C}_r)\cong D_{2pqr}.
\end{align*}
So Proposition~\ref{prop 1.5.2} applies.
\qedhere
 \end{proof}
 
 The next theorem tells us that if we have a finite group that can be broken into a semidirect product of two cyclic subgroups, then there is a Hamiltonian cycle in the connected Cayley graph of this group that comes from the generators of the factors.
 
 \begin{theorem}[{}{B. Alspach \cite[Corollary 5.2]{First}}] \label{theorem 5.5}
If $ G = \langle s \rangle \ltimes \langle t \rangle $, for some elements $ s $ and $ t $ of $ G $, then $ \Cay(G ; \{s,t\}) $ has a Hamiltonian cycle.
\end{theorem}

The following lemmas show that some special Cayley graphs have a Hamiltonian cycle, and we use these facts in Section~\ref{ch:6} in order to prove our main result. 
 \begin{lemma} \label{lemma 5.16.}
Assume $ G = (\mathcal{C}_2\times\mathcal{C}_r)\ltimes G' $, and $ G' = \mathcal{C}_p\times\mathcal{C}_q $, where $ p $, $ q $ and $ r $ are distinct prime numbers and let $ S = \{a,b\} $ be a generating set of $ G $. Additionally, assume $ |\overline{a}| \in \{2,2r\} $, $ |\overline{b}| = r $ and $ \gcd(|b|,r-1) = 1 $. Then $ \Cay(G;S) $ contains a Hamiltonian cycle.
\end{lemma}
\begin{proof}
        We have $ C = (\overline{b}^{r-1},\overline{a},\overline{b}^{-(r-1)},\overline{a}^{-1}) $ as a Hamiltonian cycle in $ \Cay(\overline{G};\overline{S}) $. Now we calculate its voltage
        \begin{align*}
            \mathbb{V}(C) = b^{r-1}ab^{-(r-1)}a^{-1} = [b^{r-1},a].
        \end{align*}
        Since $ \gcd(|b|,r-1) = 1 $, then by Lemma~\ref{lemma 5.13.} we have $ [b^{r-1},a] = G' $. Therefore, Factor Group Lemma~\ref{FGL} applies.
\end{proof}

 \begin{lemma}[cf.~{\cite[Case 2 of proof of Theorem 1.1, pages 3619-3620]{Sixth}}] \label{lemma5.10}
Assume 
\begin{itemize}
    \item $ G = (\mathcal{C}_2\times\mathcal{C}_r)\ltimes (\mathcal{C}_p\times\mathcal{C}_q) $,
    \item $ |S| = 3 $,
    \item $ \widehat{S} $ is a minimal generating set of $ \widehat{G} = G/\mathcal{C}_p $,
    \item $ \mathcal{C}_r $ centralizes $ \mathcal{C}_q $,
    \item $ \mathcal{C}_2 $ inverts $ \mathcal{C}_q $.
    \end{itemize}
 Then, $ \Cay(G;S) $ contains a Hamiltonian cycle.
\end{lemma}

\begin{lemma}[{}{\cite[Lemma 2.6]{Sixth}}] \label{lemma 5.7}
Assume:
\begin{itemize}
    \item $ G = \langle a \rangle \ltimes \langle S_{0} \rangle $, where $ \langle S_{0} \rangle $ is an abelian subgroup of odd order,
    \item $ |(S_{0} \cup S_{0}^{-1})| \geq 3, $ and
    \item $ \langle S_{0} \rangle $ has a nontrivial subgroup $ H $, such that $ H \triangleleft G $ and $ H \cap Z(G) = \{e\} $.
\end{itemize}
Then $ \Cay(G; S_{0} \cup \{a\}) $ has a Hamiltonian cycle.
\end{lemma}

\begin{lemma}[{}{\cite[Lemma 2.9]{Sixth}}] \label{lemma 5.8}
If $ G = D_{2pq} \times \mathcal{C}_r $, where $ p, q $ and $ r $ are distinct odd primes, then every connected Cayley graph on $ G $ has a Hamiltonian cycle.
\end{lemma}

Now we prove the Proposition~\ref{pro 1.1.4} which is on page~\pageref{pro 1.1.4}.

\begin{proof}[Proof of Proposition~\ref{pro 1.1.4}]\label{7pq}
If $ p \neq 7 $ and $ q \neq 7 $, then Theorem~\ref{theorem 1.2}(\ref{theorem 1.2.3}) applies. So we may assume $ q = 7 $, which means $|G| = 49p $ (and $p \neq 7$). We may also assume that $G$ is not abelian, for otherwise Lemma~\ref{abelain group} applies.

If a Sylow $p$-subgroup $ P $ of $G$ is normal, then $|G/P| = 49$, so the quotient $G/P$ is abelian.~(Because if $q$ is prime, then every group of order $q^2$ is abelian). Therefore, since $P$ is normal and $G/P$ is abelian, then $G'$ is contained in $P$. So $|G'| = p$. Therefore, Theorem~\ref{Keating-Witte} applies.

Now we may assume $P$ is not normal in $G$. Then by Sylow's Theorem, $ n_p|49 $ and $ n_p \equiv 1 \pmod{p} $, where $ n_p $ is the number of Sylow $p$-subgroups in $ G $. Thus, $ p \in \{2,3\} $, so $|G|\in \{14q,21q\}$. Therefore, Theorem~\ref{theorem 1.2}(\ref{theorem 1.2.1}) applies.
\end{proof}

\subsection{Some specific sets that generate $ G $}

This Subsection presents a few results that provide conditions under which certain $2$-element subsets generate $G$. Obviously, no $3$-element minimal generating set can contain any of these subsets.

\begin{lemma} \label{lemma 5.14}
Assume $ G = (\mathcal{C}_2\times\mathcal{C}_3)\ltimes G' $, and $ G' = \mathcal{C}_p\times\mathcal{C}_q$. Also, assume $ C_{G'}(\mathcal{C}_3) = \mathcal{C}_q $ and $ \mathcal{C}_q \not\subseteq C_{G'}(\mathcal{C}_2) $. If $ (a,b) $ is one of the following ordered pairs
\begin{enumerate}
    \item $ (\gt\gq,a_{2}\gt^j\gq^k\gamma_{p}) $, \label{lemma 5.14.1}
    \item $ (a_{2}\gt,\gt^j\gq^k\gamma_{p}) $, where $ k \not\equiv 0 \pmod{q} $, \label{lemma 5.14.3}
    \item $ (a_{2}\gt\gq,\gt^j\gq^k\gamma_{p}) $, where $ k \not\equiv 0 \pmod{q} $, \label{lemma 5.14.4}
    \item $ (a_{2}\gt\gq, a_{2}\gt^j\gq^k\gamma_{p}) $, where $k \not\equiv 1 \pmod{q}$, \label{lemma 5.14.5}
\end{enumerate}
then $ \langle a,b \rangle = G $.
\end{lemma}
\begin{proof}
It is easy to see that $(\overline{a}, \overline{b}) = \overline{G}$, so it suffices to show that $\langle a,b \rangle$ contains $\mathcal{C}_p$ and $\mathcal{C}_q$. Thus, it suffices to show that $\breve{G}$ and $\widecheck{G}$ are nonabelian, where $\breve{G} = G/(\mathcal{C}_3\ltimes \mathcal{C}_p)\cong D_{2q}$ and $\widecheck{G} = G/C_q$.  

Since $\gt$ does not centralize~$\mathcal{C}_p$, it is clear in each of $(1)-(4)$ that $\widecheck{a}$ does not centralize~$\gamma_p$ (and $\gamma_p$ is one of the factors in~$\widecheck{b}$), so $\widecheck{G}$ is not abelian. 

The pair $(\breve{a},\Breve{b})$ is either $(a_q, a_2 a_q^k)$, $(a_2, a_q^k)$ where $k \not\equiv 0 \pmod{q}$, $(a_2 a_q, a_q^k)$ where $k \not\equiv 0 \pmod{q} $, or $(a_2 a_q, a_2 a_q^k)$ where $k \not\equiv  1 \pmod{q}$. Each of these is either a reflection and a nontrivial rotation or two different reflections, and therefore generates the~(nonabelian) dihedral group $D_{2q} = \breve{G}$.
\end{proof}
\begin{lemma} \label{lemma 5.15}
Assume $ G = (\mathcal{C}_2\times\mathcal{C}_3)\ltimes G' $, and $ G' = \mathcal{C}_p\times\mathcal{C}_q $. Also, assume $ C_{G'}(\mathcal{C}_3) = \{e\} $. If $ (a,b) $ is one of the following ordered pairs
\begin{enumerate}
    \item $ (a_{2}\gt,a_{2}^i\gt^j\gq^k\gamma_{p}) $, where $ k \not\equiv 0 \pmod{q} $, \label{5.15.1}
     \item $ (\gt\gq,a_{2}\gt^j\gamma_{p}) $, where $ j \not\equiv 0 \pmod{3} $, \label{5.15.2}
    \item $ (\gt,a_{2}\gt^j\gq^k\gamma_{p}) $, where $ k \not\equiv 0 \pmod{q} $, \label{5.15.3}
    \item $ (a_{2}\gt\gq,a_{2}^i\gt^j\gamma_{p}) $, where $ j \not\equiv 0 \pmod{3} $, \label{lemma 5.15.5}
\end{enumerate}
then $ \langle a,b \rangle = G $.
\end{lemma}
\begin{proof}
It is easy to see that $(\overline{a}, \overline{b}) = \overline{G}$, so it suffices to show that $\langle a,b \rangle$ contains $\mathcal{C}_p$ and $\mathcal{C}_q$. we need to show that $\widehat{G}$ and $\widecheck{G}$ are nonabelian, where $\widehat{G} = G/\mathcal{C}_p$ and $\widecheck{G} = G/\mathcal{C}_q$, as usual.

As in the proof of Lemma~\ref{lemma 5.14}, since $\gt$ does not centralize~$\mathcal{C}_p$, it is clear in each of $(1)-(4)$ that $\widecheck{a}$ does not centralize~$\gamma_p$ (and $\gamma_p$ is one of the factors in~$\widecheck{b}$), so $\widecheck{G}$ is not abelian. 

In $(1)-(4)$, $\gq$ appears in one of the generators in~$(\widehat{a},\widehat{b})$, but not the other, and the other generator does have an occurrence of $\gt$. Since $\gt$ does not centralize~$\gq$, this implies that $\widehat{G}$ is not abelian.
\end{proof}
\begin{lemma} \label{lemma 5.16}
Assume $ G = (\mathcal{C}_2\times\mathcal{C}_q)\ltimes G' $, and $ G' = \mathcal{C}_3\times\mathcal{C}_p $. Also, assume $ C_{G'}(\mathcal{C}_q) = \mathcal{C}_3 $ and $ \mathcal{C}_3 \not\subseteq C_{G'}(\mathcal{C}_2) $. If $ (a,b) $ is one of the following ordered pairs
\begin{enumerate}
    \item $ (a_{2}\gq,a_{2}^i\gq^j\gt^k\gamma_{p}) $, where $ k \not\equiv 0 \pmod{q} $, \label{lemma 5.16.1}
    \item $ (\gq\gt,a_{2}\gq^j\gt^k\gamma_{p}) $, \label{lemma 5.16.2}
    \item $ (a_2^i\gq^m\gt,a_2\gq^j\gamma_p) $, where $ m \not\equiv 0 \pmod{q} $, \label{lemma 5.16.3}
\end{enumerate}
then $ G = \langle a,b \rangle $.
\end{lemma}
\begin{proof}
It is easy to see that $(\overline{a}, \overline{b}) = \overline{G}$, so it suffices to show that $\langle a,b \rangle$ contains $\mathcal{C}_p$ and $\mathcal{C}_3$. We need to show that $\breve{G}$ and $\overleftrightarrow{G}$ are nonabelian, where $\breve{G} = G/(\mathcal{C}_q \ltimes\mathcal{C}_p) \cong D_6$ and $\overleftrightarrow{G} = G/\mathcal{C}_3$.  

In each of $(1)-(4)$, $\gq$ appears in~$\overleftrightarrow{a}$, and $\gamma_p$ appears in~$\overleftrightarrow{b}$ (but not in $\overleftrightarrow{a}$). Since $\gq$ does not centralize $\gamma_p$, this implies that $\overleftrightarrow{G}$ is not abelian.

In each of $(1)-(4)$, $(\overleftrightarrow{a},\overleftrightarrow{b})$ consists of either a reflection and a nontrivial rotation or two different reflections, so it generates the~(nonabelian) dihedral group $D_{6} = \overleftrightarrow{G}$.
\end{proof}

\section{Proof of the Main Result}\label{ch:6}

In this section we prove Theorem~\ref{theorem1.1}, which is the main result. We are given a generating set $S$ of a finite group $G$ of order $6pq$, where $p$ and $q$ are distinct prime numbers, and we wish to show $\Cay(G;S)$ contains a Hamiltonian cycle. The proof is a long case-by-case analysis.~(See Figures~\ref{outline2},~\ref{outline1} and~\ref{outline4} for outlines of the many cases that are considered.) Here are our main assumptions through the whole section.
\begin{assumption} We assume: \label{assumption 3.1}
\begin{enumerate}
    \item $ p,q > 7 $, otherwise Theorem~\ref{theorem 1.2}(\ref{theorem 1.2.1}) applies. \label{assumption 3.1.1}
    \item $ |G| $ is square-free, otherwise Proposition~\ref{proposition 5.9} applies.
    \item $ G' \cap Z(G) = \{e\} $, by Proposition~\ref{Hall Theorem}(\ref{Hall Theorem2}). \label{assumption 3.1.2}
    \item $ G \cong \mathcal{C}_n \ltimes G' $, by Proposition~\ref{Hall Theorem}(\ref{Hall Theorem3}).
    \item $ |G'| \in \{pq,3p\} $, by Corollary~\ref{odd order commutator}.
    \item For every element $ \overline{s} \in \overline{S} $, $ |\overline{s}| \neq 1 $.~Otherwise, if $ |\overline{s}| = 1 $, then $ s \in G' $, so $ G' = \langle s \rangle $ or $ |s| $ is prime. In each case $ \Cay(G/\langle s \rangle;\overline{S}) $ has a Hamiltonian cycle by part~\ref{theorem 1.2.2} or~\ref{theorem 1.2.3} of Theorem~\ref{theorem 1.2}. By Assumption~\ref{assumption 3.1}(\ref{assumption 3.1.2}), $ \langle s \rangle \cap Z(G) = \{e\} $, therefore, Lemma~\ref{lemma 5.6}(\ref{lemma 5.6 n2}) applies. \label{assumption 3.1.6}
    \item $ S $ is a minimal generating set of $ G $.~(Note that $ S $ must generate $ G $, for otherwise $\Cay(G;S)$ is not connected. Also, in order to show that every connected Cayley graph on $G$ contains a Hamiltonian cycle, it suffices to consider $\Cay(G;S)$, where $S$ is a generating set that is minimal.)
\end{enumerate}
\end{assumption}

\subsection{Assume \texorpdfstring{$ |S| = 2 $ and $ G' = \mathcal{C}_p \times \mathcal{C}_q $}{Lg}}\hfill\label{3.1}

In this subsection we prove the part of Theorem~\ref{theorem1.1} where, $ |S| = 2 $ and $ G' = \mathcal{C}_p \times \mathcal{C}_q $. Recall $ \overline{G} = G/G' $ and $ \widehat{G} = G/\mathcal{C}_p $.
\begin{figure}
\setlist[myEnumerate, 1]{itemindent = -4mm} \setlist[myEnumerate, 2]{itemindent = -8mm}
\setlist[myEnumerate, 3]{itemindent = -12mm}
\setlist[myEnumerate, 4]{itemindent = -16mm}
\setlist[myEnumerate, 5]{itemindent = -20mm}
    \begin{myEnumerate}[I.]
        \item $ |S| = 2 $
        \begin{myEnumerate}[A.]
            \item $ G' = \mathcal{C}_p\times\mathcal{C}_q $ (Section~\ref{3.1}).
            \begin{myEnumerate}[1.]
                \item $ \overline{S} $ is a minimal generating set.
                \item $ \overline{S} $ is not a minimal generating set.
            \end{myEnumerate}
            \item $ G' = \mathcal{C}_3\times\mathcal{C}_p $ (Section~\ref{3.2}).
            \begin{myEnumerate}[1.]
                \item $ |\overline{a}| = |\overline{b}| = 2q $.
                \item $ |\overline{a}| = q $.
                \item $ |\overline{a}| = 2q $ and $ |\overline{b}| = 2 $.
                \item None of the previous cases apply.
            \end{myEnumerate}
        \end{myEnumerate}
    \end{myEnumerate}
   \caption{Outline of the cases in the proof of Theorem~\ref{theorem1.1} where $ |S| = 2 $}
    \label{outline2}
\end{figure}

\begin{figure}
\setlist[myEnumerate, 1]{itemindent = -5mm} \setlist[myEnumerate, 2]{itemindent = -10mm}
\setlist[myEnumerate, 3]{itemindent = -15mm}
\setlist[myEnumerate, 4]{itemindent = -14mm}
\begin{multicols}{2}
\begin{enumerate}[I.]
\setcounter{enumi}{1}
    \item $ |S| = 3 $.
    \begin{myEnumerate}[A.]
        \item $ G' = \mathcal{C}_p\times\mathcal{C}_q $.
        \begin{myEnumerate}[a.]
            \item $ C_{G'}(\mathcal{C}_3) \neq \{e\} $ or $ \widehat{S} $ is minimal.
            \begin{myEnumerate}[i.]
                \item $ C_{G'}(\mathcal{C}_3) \neq \{e\} $~(Section~\ref{3.4}).
                \begin{myEnumerate}[1.]
                    \item $ a = a_2 $ and $ b = \gq\gt $.
                    \item $ a = a_2 $ and $ b = a_2\gq\gt $.
                    \item $ a = a_2\gt $ and $ b = a_2\gq $. 
                    \item $ a = a_2\gt $ and $ b = \gq\gt $.
                    \item $ a = a_2\gt $ and $ b = a_2\gt\gq $.
                \end{myEnumerate}
                \item $ \widehat{S} $ is minimal (Section~\ref{3.3}).
                \begin{myEnumerate}[1.]
                    \item $ C_{G'}(\mathcal{C}_2) = \mathcal{C}_p\times\mathcal{C}_q $.
                    \item $ C_{G'}(\mathcal{C}_2) = \mathcal{C}_q $.
                    \item $ C_{G'}(\mathcal{C}_2) = \mathcal{C}_p $.
                    \item $ C_{G'}(\mathcal{C}_2) = \{e\} $.
                \end{myEnumerate}
            \end{myEnumerate}
            \item $ C_{G'}(\mathcal{C}_3) = \{e\} $ and $ \widehat{S} $ is not minimal.
            \begin{myEnumerate}[i.]
                \item $ C_{G'}(\mathcal{C}_2) = \mathcal{C}_p\times\mathcal{C}_q $ (Section~\ref{3.5}).
                \begin{myEnumerate}[1.]
                    \item $ a = \gt $ and $ b = a_2\gq $.
                    \item $ a = \gt $ and $ b = a_2\gt\gq $.
                    \item $ a = a_2\gt $ and $ b = \gt\gq $.
                    \item $ a = a_2\gt $ and $ b = a_2\gq $.
                    \item $ a = a_2\gt $ and $ b = a_2\gt\gq $.
                \end{myEnumerate}
                \item $ C_{G'}(\mathcal{C}_2) \neq \{e\} $ (Section~\ref{3.6}).
                \begin{myEnumerate}[1.]
                    \item $ a = a_2\gt $ and $ b = a_2\gt\gq $.
                    \item $ a = a_2\gt $ and $ b = a_2\gq $.
                    \item $ a = a_2\gt $ and $ b = \gt\gq $.
                    \item $ a = \gt $ and $ b = a_2\gq $.
                \end{myEnumerate}
                \item $ C_{G'}(\mathcal{C}_2) = \{e\} $ (Section~\ref{3.7}).
                \begin{myEnumerate}[1.]
                    \item $ a = a_2\gt $ and $ b = a_2\gt\gq $.
                    \item $ a = a_2\gt $ and $ b = a_2\gq $.
                    \item $ a = a_2\gt $ and $ b = \gt\gq $.
                    \item $ a = \gt $ and $ b = a_2\gq $.
                \end{myEnumerate}
            \end{myEnumerate}
        \end{myEnumerate}
        \item $ G' = \mathcal{C}_3\times\mathcal{C}_p $. (Section~\ref{3.8}).
        \begin{myEnumerate}[1.]
            \item $ a = a_2\gq $ and $ b = a_2\gq^m\gt $.
            \item $ a = a_2\gq $ and $ b = a_2\gt $.
            \item $ a = a_2\gq $ and $ b = \gq^m\gt $.
            \item $ a = a_2 $ and $ b = \gq\gt $.
        \end{myEnumerate}
    \end{myEnumerate}
\end{enumerate}
\end{multicols}
 \caption{Outline of the cases in the proof of Theorem~\ref{theorem1.1} where $ |S| = 3 $}\label{outline1} 
\end{figure}

 \begin{figure}
\setlist[enumerate, 1]{itemindent = -4mm}
\setlist[enumerate, 2]{itemindent = -8mm}
\setlist[enumerate, 3]{itemindent = -12mm}
\setlist[enumerate, 4]{itemindent = -16mm}
\setlist[enumerate, 5]{itemindent = -20mm}
 \begin{enumerate}[I.]
 \setcounter{enumi}{2}
     \item $ |S| \geq 4 $~(Section~\ref{3.9}). This part of the proof applies whenever $ |G| = pqrt $ with $ p $, $ q $, $ r $, and $ t $ distinct primes.
     \begin{enumerate}[1.]
         \item $ |G'| $ has only two prime factors.
         \item $ |G'| $ has three prime factors.
     \end{enumerate}
 \end{enumerate}
 \caption{Outline of the cases in the proof of Theorem~\ref{theorem1.1} where $ |S| \geq 4 $}
    \label{outline4}
\end{figure}

\begin{caseprop}
Assume
\begin{itemize}
    \item $ G = (\mathcal{C}_2\times \mathcal{C}_3) \ltimes (\mathcal{C}_p\times \mathcal{C}_q) $,
    \item $ |S| = 2 $.
\end{itemize}
Then $ \Cay(G;S) $ contains a Hamiltonian cycle.
\end{caseprop}

\begin{figure}[hbt!]
\begin{center}
\includegraphics[width=0.7\textwidth]{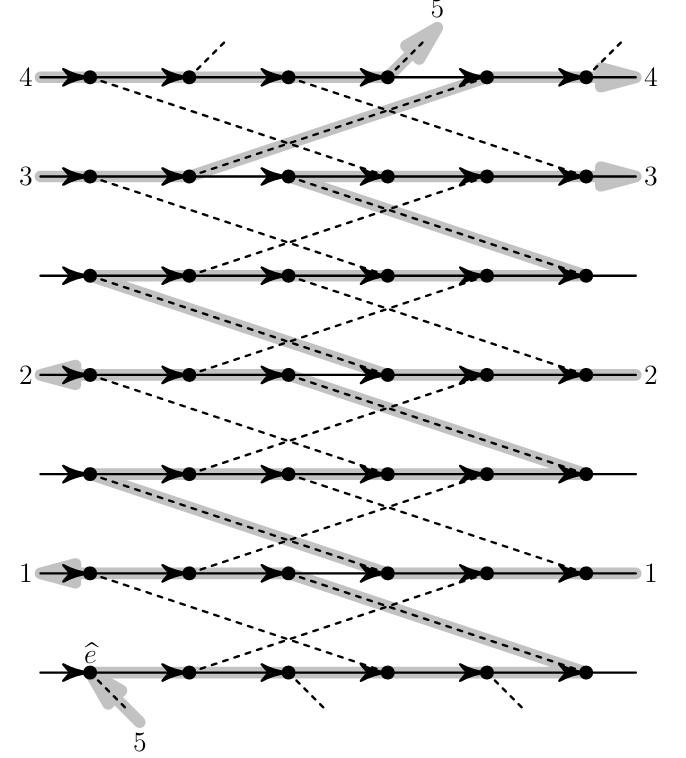}
\caption{The Hamiltonian cycle $ C_1 $: $ \widehat{a} $ edges are solid and $ \widehat{b} $ edges are dashed.}\label{Fig1}
\end{center}
    \end{figure}

\begin{proof}
   Let $ S = \{a,b\} $. For every $ s \in S $, $ |\overline{s}| \neq 1 $, by Assumption~\ref{assumption 3.1}(\ref{assumption 3.1.6}).
   \begin{case}\label{sec:3.1 case 1}
    Assume $ \overline{S} $ is minimal. Then $ |\overline{a}|,|\overline{b}| \in \{2,3\} $. When $ |\overline{a}| = |\overline{b}| = 2 $ or $  |\overline{a}| = |\overline{b}| = 3 $, then $ \overline{G} \neq \langle \overline{a},\overline{b} \rangle $. Therefore, $ G \neq \langle a, b \rangle  $ which contradicts the fact that $ G = \langle a,b \rangle $. So we may assume $ |\overline{a}| = 2 $ and $ |\overline{b}| = 3 $. Since $ |b| \in \{3,3p,3q,3pq\} $, then $ \gcd(|b|,2) = 1 $. Thus, Lemma~\ref{lemma 5.16.} applies.
\end{case}

\begin{case}
    Assume  $ \overline{S} $ is not minimal. Then $\{|\overline{a}|,|\overline{b}|\}$ is either $\{6,2\}$, $\{6,3\}$, or $\{6\}$. We may assume $ |\overline{a}| = 6 $.
\end{case}

\begin{subcase} \label{Subcase6.2.1}
Assume $ |\overline{b}| = 2 $. So we have $ \overline{b} = \overline{a}^{3}$, then $ b = a^{3}\gamma $, where $ G' = \langle \gamma \rangle $ (otherwise $ \langle a, b \rangle = \langle a, a^{3}\gamma \rangle = \langle a, \gamma \rangle \neq G $ which contradicts the fact that $ G = \langle a,b \rangle $). Now by Proposition~\ref{Hall Theorem}(\ref{Hall Theorem4}), we have $ \tau \in  \mathbb{Z}^{+} $ such that $ a\gamma a^{-1} = \gamma^{\tau} $ and $ \tau^{6}\equiv 1 \pmod{pq} $, also $ \gcd(\tau -1, pq)=1 $. This implies that $ \tau \not \equiv 1 \pmod{p} $ and $ \tau \not \equiv 1 \pmod{q} $. We have $ C_{1} = (\overline{a}^2,\overline{b},\overline{a}^{-2},\overline{b}^{-1}) $ as a Hamiltonian cycle in $ \Cay(\overline{G};\overline{S}) $. Now we calculate its voltage.
\begin{align*}
    \mathbb{V}(C_1) = a^2ba^{-2}b^{-1} = a^2a^3\gamma a^{-2}\gamma^{-1}a^{-3} = \gamma^{\tau^5-\tau^3} = \gamma^{\tau^3(\tau^2-1)}.
\end{align*}
We may assume $ \gcd(\tau^{2}-1, pq)\neq 1 $ (otherwise Factor Group Lemma~\ref{FGL} applies). Without loss of generality let $ \tau^2 \equiv 1 \pmod q $, then $ \tau \equiv -1 \pmod{q} $. We may assume $\tau\not\equiv -1 \pmod{p}$, for otherwise $G\cong D_{2pq}\times \mathcal{C}_3$, so Lemma~\ref{lemma 5.8} applies.

Consider $ \widehat{G} = G/\mathcal{C}_p = \mathcal{C}_{6} \ltimes \mathcal{C}_q $. Since $ |\overline{a}| = 6 $, then by Lemma~\ref{lemma 5.13} $ |a| = 6 $, so $ |\widehat{a}| = 6 $. We may assume $ |\widehat{b}| = 2 $, for otherwise Corollary~\ref{corollary 5.2} applies with $ s = b $ and $ t = b^{-1} $ since $ \langle \widehat{a} \rangle \neq \widehat{G} $, so any Hamiltonian cycle must use an edge labeled $ \widehat{b} $. Thus, $ \widehat{b} = \widehat{a}^3\gq $, where $ \langle a_{q} \rangle = \mathcal{C}_q $. Since $ \tau \equiv -1 \pmod{q} $, then $ \mathcal{C}_3 $ centralizes $ \mathcal{C}_q $ and $ \mathcal{C}_2 $ inverts $ \mathcal{C}_q $. Therefore, $ \widehat{G} \cong D_{2q} \times \mathcal{C}_3 $. Now we have
\begin{align*}
    C_2 = ((\widehat{a}^5,\widehat{b},\widehat{a}^{-5},\widehat{b})^{(q-3)/2},(\widehat{a}^{5},\widehat{b})^{3})
\end{align*}
 as a Hamiltonian cycle in $ \Cay(\widehat{G};\widehat{S}) $. The picture in Figure~\ref{Fig1} on page~\pageref{Fig1} shows the Hamiltonian cycle when $ q = 7 $.
If in $ C_{2} $ we change one occurrence of $ (\widehat{a}^5,\widehat{b},\widehat{a}^{-5},\widehat{b}) $ to $ (\widehat{a}^{-5},\widehat{b},\widehat{a}^{5},\widehat{b})$ we have another Hamiltonian cycle. Note that,
\begin{align*}
   a^5ba^{-5}b = a^5\cdot a^3\gamma\cdot a^{-5}\cdot a^3\gamma = a^2\gamma a^{-2}\gamma = \gamma^{\tau^2+1},
\end{align*}
and
\begin{align*}
    a^{-5}ba^5b = a^{-5}\cdot a^3\gamma\cdot a^5\cdot a^3\gamma = a^{-2}\gamma a^2\gamma = \gamma^{\tau^{-2}+1}.
\end{align*}
Since $ \tau^4 \not\equiv 0 \pmod{p} $ we see that $ \tau^2 +1 \not\equiv \tau^{-2}+1 \pmod{p} $. Therefore, the voltages of these two Hamiltonian cycles are different, so one of these Hamiltonian cycles has a nontrivial voltage. Thus, Factor Group Lemma~\ref{FGL} applies.
\end{subcase}

\begin{subcase} \label{subsubcase1.2.2}
    Assume $ |\overline{b}| = 3 $. Since $ |\overline{b}| = 3 $, then $ |b| \in \{3, 3p, 3q, 3pq\} $. Since $ |\overline{a}| = 6 $, then by \ref{lemma 5.13} $ |a| = 6 $. Since $ \gcd(|b|,2) = 1 $, then Lemma~\ref{lemma 5.16.} applies.
\end{subcase}

\begin{subcase} 
    Assume $ |\overline{b}| = 6 $. Then we have $ \overline{a} = \overline{b} $ or $ \overline{a} = \overline{b}^{-1} $. Additionally, by Lemma~\ref{lemma 5.13} we have $ |a| = |b| = 6 $. We may assume $ \overline{a} = \overline{b} $ by replacing $b$ with its inverse if necessary. Then $ b = a\gamma $, where $ G' = \langle \gamma \rangle $, because  $ G = \langle a,b \rangle  $. We have $ C = (\overline{a}^5,\overline{b}) $ as a Hamiltonian cycle in $ \Cay (\overline{G},\overline{S}) $. Now we calculate its voltage
    \begin{align*}
        \mathbb{V}(C) = a^5b = a^5a\gamma = a^6\gamma = \gamma
    \end{align*}
    which generates $ G' $. Therefore, Factor Group Lemma~\ref{FGL} applies. 
    \qedhere
\end{subcase} 
\end{proof}

\subsection{Assume \texorpdfstring{$ |S| = 2 $ and $ G' = \mathcal{C}_3 \times \mathcal{C}_p $}{Lg}}\hfill\label{3.2}

In this subsection we prove the part of Theorem~\ref{theorem1.1} where, $ |S| = 2 $ and $ G' = \mathcal{C}_3 \times \mathcal{C}_p $. Recall $ \overline{G} = G/G' $ and $ \widehat{G} = G/\mathcal{C}_p $.
\begin{caseprop}
Assume
\begin{itemize}
    \item $ G = (\mathcal{C}_2\times\mathcal{C}_q) \ltimes(\mathcal{C}_3\times\mathcal{C}_p) $,
    \item $ |S| = 2 $.
\end{itemize}
Then $ \Cay(G;S) $ contains a Hamiltonian cycle.
\end{caseprop}
\begin{proof}
        Let $ S = \{a,b\} $. Since the only non-trivial automorphism of $ \mathcal{C}_3 $ is inversion, $ \mathcal{C}_q $ centralizes $ \mathcal{C}_3 $. Since $ G' \cap Z(G) = \{e\} $ (see Proposition~\ref{Hall Theorem}(\ref{Hall Theorem4})), $ \mathcal{C}_2 $ does not centralize $ \mathcal{C}_3 $.
        \setcounter{case}{0}
        
          \begin{case}\label{Case 2.1}
            Assume $ |\overline{a}| = |\overline{b}| = 2q $. Then $ \overline{b} = \overline{a}^m $, where $ 1\leq m \leq q-1 $ by replacing $ b $ with its inverse if needed. Therefore, $ b = a^m\gamma $, where $ G' = \langle \gamma \rangle $. Also, $ \gcd(m,2q) = 1 $. So, by Proposition~\ref{Hall Theorem}(\ref{Hall Theorem4}) we have $ a\gamma a^{-1} = \gamma^{\tau} $ where $ \tau^{2q} \equiv 1 \pmod{3p} $ and $ \gcd(\tau-1,3p) = 1 $. Consider $ \overline{G} = \mathcal{C}_{2q} $.
            \begin{subcase}
                Assume $ m >3 $. Then we have 
                \begin{align*}
         C = (\overline{b}^{-2},\overline{a}^{-2},\overline{b},\overline{a},\overline{b},\overline{a}^{-(m-2)},\overline{b}^{-1},\overline{a}^{m-4},\overline{b}^{-1},\overline{a}^{-(2q-2m-3)})
     \end{align*}
     as a Hamiltonian cycle in $ \Cay(\overline{G};\overline{S}) $. Now we calculate its voltage.
    \begin{align*}
         \mathbb{V}(C) &= b^{-2}a^{-2}baba^{-(m-2)}b^{-1}a^{m-4}b^{-1}a^{-(2q-2m-3)} \\&= \gamma^{-1}a^{-m}\gamma^{-1}a^{-m}a^{-2}a^m\gamma a a^m\gamma a^{-m+2}\gamma^{-1}a^{-m}a^{m-4}\gamma^{-1}a^{-m}a^{-2q+2m+3}\\&= \gamma^{-1}a^{-m}\gamma^{-1}a^{-2}\gamma a^{m+1}\gamma a^{-m+2}\gamma^{-1}a^{-4}\gamma^{-1}a^{m+3}\\&= \gamma^{-1-\tau^{-m}+\tau^{-m-2}+\tau^{-1}-\tau^{-m+1}-\tau^{-m-3}} \\&= \gamma^{-1+\tau^{-1}-\tau^{-m+1}-\tau^{-m}+\tau^{-m-2}-\tau^{-m-3}}.
     \end{align*}
     We may assume $ \mathbb{V}(C) $ does not generate $ G' = \mathcal{C}_3\times\mathcal{C}_p $. Therefore, the subgroup generated by $ \mathbb{V}(C) $ either does not contain $ \mathcal{C}_3 $, or does not contain $ \mathcal{C}_p $. We already know $ \tau \equiv -1 \pmod{3} $, then we have
     \begin{align*}
         -1+\tau^{-1}-\tau^{-m+1}-\tau^{-m}+\tau^{-m-2}-\tau^{-m-3} &\equiv -1-1-1+1-1-1 \pmod{3} \\&= -4 = -1.
     \end{align*}
     This implies that the subgroup generated by $ \mathbb{V}(C) $ contains $ \mathcal{C}_3 $. So we may assume the subgroup generated by $ \mathbb{V}(C) $ does not contain $ \mathcal{C}_p $, then
     \begin{align*}
        0 \equiv -1+\tau^{-1}-\tau^{-m+1}-\tau^{-m}+\tau^{-m-2}-\tau^{-m-3}\pmod{p}. \tag{1.1A}\label{7.1B}
     \end{align*}
    Multiplying by $ -\tau^{m+3} $ we have
     \begin{align*}
         0 \equiv \tau^{m+3}-\tau^{m+2}+\tau^4+\tau^3-\tau+1 \pmod{p}.\tag{1.1B}\label{7.1C}
     \end{align*}
     Replacing $ \{\overline{a},\overline{b}\} $  with $ \{\overline{a}^{-1},\overline{b}^{-1}\} $  replaces $ \tau $ with $ \tau^{-1} $. Therefore, applying the above argument to $ \{\overline{a}^{-1},\overline{b}^{-1}\} $ establishes that \ref{7.1B} holds with $ \tau^{-1} $ in the place of $ \tau $, which means we have
     \begin{align*}
        0 \equiv -\tau^{m+3}+\tau^{m+2}-\tau^m-\tau^{m-1}+\tau-1  \pmod{p}. \tag{1.1C}\label{7.1D}
     \end{align*}
      By adding \ref{7.1C} and \ref{7.1D} we have
     \begin{align*}
         0 \equiv -\tau^m-\tau^{m-1}+\tau^4+\tau^3 = \tau^3(\tau+1)(1-\tau^{m-4})  \pmod{p}.
     \end{align*}
     If $ \tau \equiv -1 \pmod{p} $, then $ \mathcal{C}_{2q} $ inverts $ \mathcal{C}_{3p} $, so $ \mathcal{C}_q $ centralizes $ \mathcal{C}_p $. This implies that $ G \cong D_{6p}\times\mathcal{C}_q $, so Lemma~\ref{lemma 5.8} applies. The only other possibility is $ \tau^{m-4} \equiv 1 \pmod{p} $. Multiplying by $ \tau^4 $, we have $ \tau^m \equiv \tau^4 \pmod{p} $. We also know that $ \tau^{2q} \equiv 1 \pmod{p} $. So $ \tau^d \equiv 1 \pmod{p} $, where $ d = \gcd(m-4,2q) $. Since $ m $ is odd and $ m < q $, then $ d = 1 $.~This~contradicts~the~fact~that~$\gcd(\tau-1,3p) = 1$.
          \end{subcase}
          
          \begin{subcase}
              Assume $ m \leq 3 $. Therefore, either $ m = 1 $ or $ m = 3 $. If $ m = 1 $, then $ \overline{a} = \overline{b} $ and $ b = a\gamma $. So we have $ C_1 = (\overline{a}^{2q-1},\overline{b}) $
               as a Hamiltonian cycle in $ \Cay(\overline{G};\overline{S}) $. Now we calculate its voltage.
              \begin{align*}
                  \mathbb{V}(C_1) &= a^{2q-1}b = a^{2q-1}a\gamma = \gamma
              \end{align*}
              which generates $ G' $. Therefore, Factor Group Lemma~\ref{FGL} applies. Now if $ m = 3 $, then $ b = a^3\gamma $ and we have 
              \begin{align*}
        C_2 = (\overline{b}^2,\overline{a}^{-1},\overline{b}^{-1},\overline{a}^{-1},\overline{b}^3,\overline{a}^{-2},\overline{b},\overline{a}^{2q-11})
    \end{align*}
    as a Hamiltonian cycle in $ \Cay(\overline{G};\overline{S}) $. We calculate its voltage.
     \begin{align*}
        \mathbb{V}(C_2) &= b^2a^{-1}b^{-1}a^{-1}b^3a^{-2}ba^{2q-11} \\&= a^3\gamma a^3\gamma a^{-1}\gamma^{-1}a^{-3}a^{-1}a^3\gamma a^3\gamma a^3\gamma a^{-2}a^3\gamma a^{-11} \\&= a^3\gamma a^3\gamma a^{-1}\gamma^{-1} a^{-1}\gamma a^3\gamma a^3\gamma a\gamma a^{-11} \\&= \gamma^{\tau^3+\tau^6-\tau^5+\tau^4+\tau^7+\tau^{10}+\tau^{11}} \\&= \gamma^{\tau^{11}+\tau^{10}+\tau^7+\tau^6-\tau^5+\tau^4+\tau^3}
    \end{align*}
    We may assume $ \mathbb{V}(C_2) $ does not generate $ G' = \mathcal{C}_3\times\mathcal{C}_p $. Therefore, the subgroup generated by $ \mathbb{V}(C) $ does not contain either $ \mathcal{C}_3 $, or $ \mathcal{C}_p $. We already know $ \tau \equiv -1 \pmod{3} $, then
    \begin{align*}
        \tau^{11}+\tau^{10}+\tau^7+\tau^6-\tau^5+\tau^4+\tau^3 \equiv -1+1-1+1+1+1-1 = 1 \pmod{3}. 
    \end{align*}
    This implies that the subgroup generated by $ \mathbb{V}(C_2) $ contains $ \mathcal{C}_3 $. So we may assume the subgroup generated by $ \mathbb{V}(C_2) $ does not contain $ \mathcal{C}_p $, for otherwise Factor Group Lemma~\ref{FGL} applies. Then we have
    \begin{align*}
        0 &\equiv \tau^{11}+\tau^{10}+\tau^7+\tau^6-\tau^5+\tau^4+\tau^3 \pmod{p} \\&= \tau^{3}(\tau^8+\tau^7+\tau^4+\tau^3-\tau^2+\tau+1).
    \end{align*}
    This implies that 
    \begin{align*}
        0 \equiv \tau^8+\tau^7+\tau^4+\tau^3-\tau^2+\tau+1 \pmod{p}. \tag{1.2A}\label{Subcase7.1.2A}
    \end{align*}
   We can replace $ \tau $ with $ \tau^{-1} $ in the above equation, by replacing $ \{\overline{a},\overline{b}\} $ with $ \{\overline{a}^{-1},\overline{b}^{-1}\} $ if necessary. Then we have
   \begin{align*}
       0 \equiv \tau^{-8}+\tau^{-7}+\tau^{-4}+\tau^{-3}-\tau^{-2}+\tau^{-1}+1 \pmod{p}.
   \end{align*}
   Multiplying $ \tau^8 $, then we have
   \begin{align*}
       0 &\equiv 1+\tau+\tau^4+\tau^5-\tau^6+\tau^7+\tau^8 \pmod{p}\\&= \tau^8+\tau^7-\tau^6+\tau^5+\tau^4+\tau+1.
   \end{align*}
   Now by subtracting the above equation from \ref{Subcase7.1.2A} we have
   \begin{align*}
       0 &\equiv \tau^6-\tau^5+\tau^3-\tau^2 \pmod{p}\\&=\tau^2(\tau-1)(\tau^3+1).
   \end{align*}
   This implies that $ \tau \equiv 1 \pmod{p} $ or $ \tau^3 \equiv -1 \pmod{p} $. If $ \tau \equiv 1 \pmod{p} $, then it contradicts the fact that $ \gcd(\tau-1,3p) = 1 $. Now if $ \tau^3 \equiv -1 \pmod{p} $, then $ \tau^6 \equiv 1 \pmod{p} $. We already know $ \tau^{2q} \equiv 1 \pmod{p} $. Then $ \tau^d \equiv 1 \pmod{p} $, where $ d = \gcd(2q,6) $. Since $ \gcd(2,6) = 2 $ and $ \gcd(q,6) = 1 $, then $ d = 2 $. This implies that $ \tau^2 \equiv 1 \pmod{p} $, which means $ \mathcal{C}_q $ centralizes $ \mathcal{C}_p $. Then we have
   \begin{align*}
       G = \mathcal{C}_q\times(\mathcal{C}_2\ltimes\mathcal{C}_{3p}) \cong \mathcal{C}_q\times D_{6p}.
   \end{align*}
   So Lemma~\ref{lemma 5.8} applies.
          \end{subcase}
        \end{case}
        
        \begin{case} \label{case2.2}
        
            Assume $ |\overline{a}| = q $. Then $ |\overline{b}| \in \{2,2q\} $. Thus $ |b| \in \{2,2q,2p,2pq\}$. If $ |b| = 2pq $, then $ \mathcal{C}_q $ centralizes $ \mathcal{C}_p $. This implies that
              \begin{align*}
                  G = \mathcal{C}_q\times(\mathcal{C}_2\ltimes\mathcal{C}_{3p}) \cong \mathcal{C}_q\times D_{6p}
              \end{align*}
              so, Lemma~\ref{lemma 5.8} applies. Therefore, we may assume $ \mathcal{C}_q $ does not centralize $ \mathcal{C}_p $, so $ |a| $ is not divisible by $ p $. If $ |b| = 2p $, then Corollary~\ref{corollary 5.2} applies with $ s = b $ and $ t = b^{-1} $, because we have a Hamiltonian cycle in $ \Cay(\widehat{G};\widehat{S}) $ by Theorem~\ref{theorem 1.2}(\ref{theorem 1.2.3}).~(Since $ b $ is the only generator whose order is even, then any Hamiltonian cycle in $ \Cay(\widehat{G};\widehat{S}) $ must use some edge labeled $ \widehat{b} $.)
              
              We may now assume $ |b| \in \{2,2q\} $. We have $ C = (\overline{a}^{q-1},\overline{b},\overline{a}^{-(q-1)},\overline{b}^{-1}) $ as a Hamiltonian cycle in $ \Cay(\overline{G};\overline{S}) $. Now if $ |a| = q $, then by Lemma~\ref{lemma 5.13.} we have $ G' = \langle [a^{q-1},b]\rangle $. Therefore, Factor Group Lemma~\ref{FGL} applies. So, we may assume $ |a| = 3q $. Since $ \mathcal{C}_q $ does not centralize $ \mathcal{C}_p $, then after conjugation we can assume $ a = \gt\gq $ and $ b = a_2\gq^j\gamma_p $, where $ 0\leq j \leq q-1 $. We already know that $ C $ is a Hamiltonian cycle in $ \Cay(\overline{G};\overline{S}) $. So we can assume $ \gcd(3q,q-1) \neq 1 $ (otherwise Lemma~\ref{lemma 5.13.} applies, which implies that Factor Group Lemma~\ref{FGL} applies). This implies that $ \gcd(3,q-1) \neq 1 $ which means $ q \equiv 1 \pmod{3} $.
              
              Consider $ \widehat{G} = G/\mathcal{C}_p $. Then $ \widehat{a} = \gt\gq $ and $ \widehat{b} = a_2\gq^j $. Therefore, there exists $ 0\leq k \leq 3q-1 $ such that $ \widehat{b}^{-1}\widehat{a}\widehat{b} = \widehat{a}^k $. Since $ \widehat{b} $ inverts $ \gt $ and centralizes $ \gq $, then we must have $ \widehat{a} = \widehat{b}\widehat{a}^k\widehat{b}^{-1} = \gt^{-k}\gq^k $, so $ k \equiv -1 \pmod{3} $ and $ k \equiv 1 \pmod{q} $. Since $ q \equiv 1 \pmod{3} $, then $ k = q+1 $. Additionally, we have $ a\gamma_{p}a^{-1} = \gamma_{p}^{\widehat{\tau}} $, where $ \widehat{\tau}^q \equiv 1 \pmod{p} $. We also have $ \widehat{\tau} \not \equiv 1 \pmod{p} $, because $ \mathcal{C}_q $ does not centralize $ \mathcal{C}_p $.
              Now we have
              \begin{align*}
                  b^{-1}ab = \gamma_{p}^{-1}\gq^{-j}a_{2}a a_{2}\gq^j\gamma_{p} =  \gamma_{p}^{-1}a^{q+1}\gamma_{p}.
              \end{align*}
             This implies that
              \begin{align*}
             b^{-1}a^ib = (b^{-1}ab)^i = (\gamma_{p}^{-1}a^{q+1}\gamma_{p})^i = \gamma_{p}^{-1}a^{i(q+1)}\gamma_{p}.
              \end{align*}
              Therefore, 
              \begin{align*}
                  b^{-1}a^ib = \gamma_{p}^{-1}a^{i(q+1)}\gamma_{p} \equiv \gamma_{p}^{-1}a^i\gamma_{p} \pmod{\mathcal{C}_3}.
              \end{align*}
             We have 
              \begin{align*}
                  C_{1} &= (\widehat{a}^{q-3},\widehat{b}^{-1},\widehat{a}^{-(q-2)},\widehat{b},\widehat{a}^{-1},\widehat{b}^{-1},\widehat{a},\widehat{b},\widehat{a}^{q-2},\widehat{b}^{-1},\\&\indent\indent\widehat{a}^{-(q-3)},\widehat{b},\widehat{a}^{q-2},\widehat{b}^{-1},\widehat{a},\widehat{b},\widehat{a}^{-1},\widehat{b}^{-1},\widehat{a}^{-(q-2)},\widehat{b})
              \end{align*}
              as our first Hamiltonian cycle in $ \Cay(\widehat{G};\widehat{S}) $. The picture in Figure~\ref{1Hamcyc7} on page~\pageref{1Hamcyc7} shows the Hamiltonian cycle. In addition,
        \begin{figure}
       \includegraphics[width=0.9\textwidth]{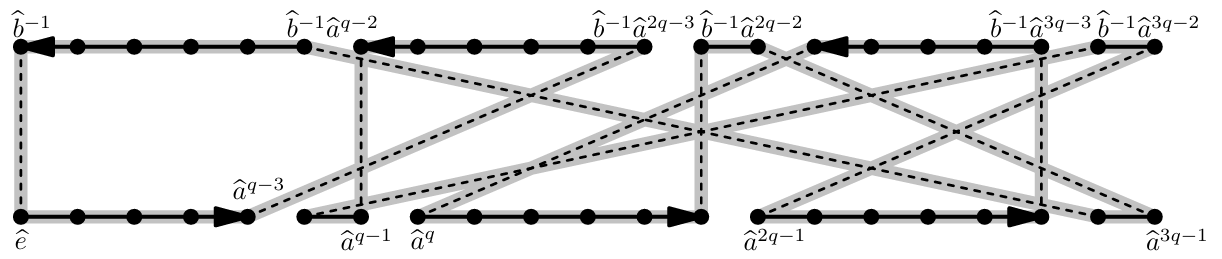}
\caption{The Hamiltonian cycle $ C_1 $: $ \widehat{a} $ edges are solid and $ \widehat{b} $ edges are dashed.}\label{1Hamcyc7}
    \end{figure}
     \begin{figure}
    \includegraphics[width=0.9\textwidth]{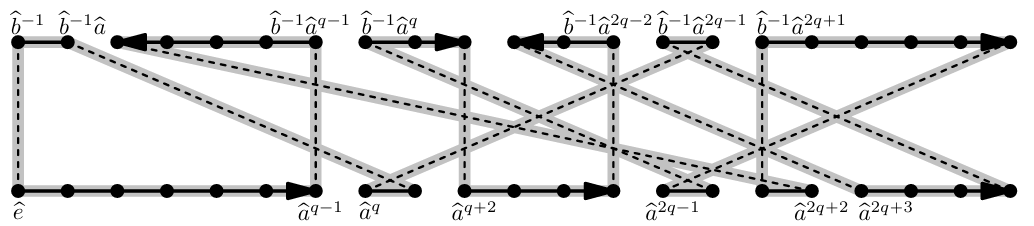}
        \caption{The Hamiltonian cycle $ C_2 $: $ \widehat{a} $ edges are solid and $ \widehat{b} $ edges are dashed.}\label{2Hamcyc7}
    \end{figure}
              \begin{align*}
                  C_{2} &= (\widehat{a}^{q-1},\widehat{b}^{-1},\widehat{a}^{-(q-3)},\widehat{b},\widehat{a}^{-1},\widehat{b}^{-1},\widehat{a}^{q-2},\widehat{b},\widehat{a},\widehat{b}^{-1},\widehat{a}^2,\widehat{b},\\&\indent\indent\widehat{a}^{q-4},\widehat{b}^{-1},\widehat{a}^{-(q-5)},\widehat{b},\widehat{a}^{q-4},\widehat{b}^{-1},\widehat{a},\widehat{b},\widehat{a},\widehat{b}^{-1},\widehat{a}^{-1},\widehat{b})
              \end{align*}
              is the second Hamiltonian cycle in $ \Cay(\widehat{G};\widehat{S}) $. The picture in Figure~\ref{2Hamcyc7} on page~\pageref{2Hamcyc7} shows the Hamiltonian cycle.
    We calculate the voltage of $ C_{1} $ in $ \overleftrightarrow{G} = G/\mathcal{C}_3 $. Since $ a^q \equiv e \pmod{\mathcal{C}_3} $, we have
     \begin{align*}
        \mathbb{V}(C_1) &\equiv a^{-3}(b^{-1}a^2b)a^{-1}(b^{-1} ab)a^{-2}(b^{-1}a^3b)a^{-2}(b^{-1}ab)a^{-1}(b^{-1}a^2b)  \pmod{\mathcal{C}_3} \\ &= a^{-3}(\gamma_{p}^{-1}a^2\gamma_{p})a^{-1}(\gamma_{p}^{-1}a\gamma_{p})a^{-2}(\gamma_{p}^{-1}a^3\gamma_{p})a^{-2}(\gamma_{p}^{-1}a\gamma_{p})a^{-1}(\gamma_{p}^{-1}a^2\gamma_{p}) \\ &= a^{-3}(\gamma_{p}^{\widehat{\tau}^2-1}a^2)a^{-1}(\gamma_{p}^{\widehat{\tau}-1}a)a^{-2}(\gamma_{p}^{\widehat{\tau}^3-1}a^3)a^{-2}(\gamma_{p}^{\widehat{\tau}-1}a)a^{-1}(\gamma_{p}^{\widehat{\tau}^2-1}a^2)  \\ &= a^{-3}\gamma_{p}^{\widehat{\tau}^2-1}a\gamma_{p}^{\widehat{\tau}-1}a^{-1}\gamma_{p}^{\widehat{\tau}^3-1}a\gamma_{p}^{\widehat{\tau}^2+\widehat{\tau}-2}a^2 \\&= \gamma_{p}^{\widehat{\tau}^{-3}(\widehat{\tau}^2-1)+\widehat{\tau}^{-2}(\widehat{\tau}-1)+\widehat{\tau}^{-3}(\widehat{\tau}^3-1)+\widehat{\tau}^{-2}(\widehat{\tau}^2+\widehat{\tau}-2)}\\&= \gamma_{p}^{-2\widehat{\tau}^{-3}-3\widehat{\tau}^{-2}+3\widehat{\tau}^{-1}+2}.
    \end{align*}
     We may assume this does not generate $ \mathcal{C}_p $, so
    \begin{align*}
      0 \equiv -2\widehat{\tau}^{-3}-3\widehat{\tau}^{-2}+3\widehat{\tau}^{-1}+2 \pmod{p}.
    \end{align*}
    Multiplying by $ \widehat{\tau}^3 $, we have
    \begin{align*}
       0 \equiv 2\widehat{\tau}^3+3\widehat{\tau}^2-3\widehat{\tau}-2 = (\widehat{\tau}-1)(\widehat{\tau}+2)(2\widehat{\tau}+1) \pmod{p}.
    \end{align*}
    Since $ \widehat{\tau} \not \equiv 1 \pmod{p} $, then we may assume $ \widehat{\tau} \equiv -2 \pmod{p} $, by replacing $ \widehat{a} $ with $ \widehat{a}^{-1} $ if needed.
    
    Now we calculate the voltage of $ C_{2} $ in $ \overleftrightarrow{G} = G/\mathcal{C}_3 $.
    \begin{align*}
        \mathbb{V}(C_2) &\equiv a^{-1}(b^{-1}a^3b)a^{-1}(b^{-1}a^{-2}b)a(b^{-1}a^2b)a^{-4}(b^{-1} a^5b)a^{-4}(b^{-1}ab)a(b^{-1}a^{-1}b) \pmod{\mathcal{C}_3} \\ &= a^{-1}(\gamma_{p}^{-1}a^3\gamma_{p})a^{-1}(\gamma_{p}^{-1}a^{-2}\gamma_{p})a(\gamma_{p}^{-1}a^2\gamma_{p}) \\& \indent\indent\cdot a^{-4}(\gamma_{p}^{-1}a^5\gamma_{p}) a^{-4}(\gamma_{p}^{-1}a\gamma_{p})a(\gamma_{p}^{-1}a^{-1}\gamma_{p}) \\&= a^{-1}(\gamma_{p}^{\widehat{\tau}^3-1}a^3)a^{-1}(\gamma_{p}^{\widehat{\tau}^{-2}-1}a^{-2})a(\gamma_{p}^{\widehat{\tau}^2-1}a^2) \\&\indent\indent\cdot a^{-4}(\gamma_{p}^{\widehat{\tau}^5-1}a^5) a^{-4}(\gamma_{p}^{\widehat{\tau}-1}a)a(\gamma_{p}^{\widehat{\tau}^{-1}-1}a^{-1}) \\&= a^{-1}\gamma_{p}^{\widehat{\tau}^3-1}a^2\gamma_{p}^{\widehat{\tau}^{-2}-1}a^{-1}\gamma_{p}^{\widehat{\tau}^2-1}a^{-2}\gamma_{p}^{\widehat{\tau}^5-1}a\gamma_{p}^{\widehat{\tau}-1}a^2\gamma_{p}^{\widehat{\tau}^{-1}-1}a^{-1} \\ &= \gamma_{p}^{\widehat{\tau}^{-1}(\widehat{\tau}^3-1)+\widehat{\tau}(\widehat{\tau}^{-2}-1)+\widehat{\tau}^2-1+\widehat{\tau}^{-2}(\widehat{\tau}^5-1)+\widehat{\tau}^{-1}(\widehat{\tau}-1)+\widehat{\tau}(\widehat{\tau}^{-1}-1)} \\&= \gamma_{p}^{\widehat{\tau}^3+2\widehat{\tau}^2-2\widehat{\tau}+1-\widehat{\tau}^{-1}-\widehat{\tau}^{-2}}.
    \end{align*}
    We may assume this does not generate $ \mathcal{C}_p $, so
    \begin{align*}
        0 \equiv \widehat{\tau}^3+2\widehat{\tau}^2-2\widehat{\tau}+1-\widehat{\tau}^{-1}-\widehat{\tau}^{-2} \pmod{p}.
    \end{align*}
    Multiplying by $ \widehat{\tau}^2 $, we have
    \begin{align*}
       0 \equiv \widehat{\tau}^5+2\widehat{\tau}^4-2\widehat{\tau}^3+\widehat{\tau}^2-\widehat{\tau}-1 \pmod{p}.
    \end{align*}
    We already know $ \widehat{\tau} \equiv -2 \pmod{p} $. By substituting this in the equation above, we have 
    \begin{align*}
      0 \equiv (-2)^5+2(-2)^4-2(-2)^3+(-2)^2-(-2)-1 = 21 = 3\cdot7 \pmod{p}.
    \end{align*}
    Since $ p>7 $, then $ 21  \not \equiv 0 \pmod{p} $. This is a contradiction.
    
        \end{case}
        
      \begin{case} \label{case 2.3}
             Assume $ |\overline{a}| = 2q $ and $ |\overline{b}| = 2 $. Since $ |\overline{a}| = 2q $, then by Lemma~\ref{lemma 5.13} $ |a| = 2q $. We have $ b = a^q\gamma $ where $ G' = \langle \gamma \rangle $. 
             
             By Proposition~\ref{Hall Theorem}(\ref{Hall Theorem4}) we have $ a\gamma a^{-1} = \gamma ^{\tau} $, where $ \tau^{2q} \equiv 1 \pmod{3p} $ and $ \gcd(\tau -1, 3p) = 1 $. This implies that $ \tau \not \equiv 0, 1 \pmod{p} $ and $ \tau \equiv -1 \pmod{3} $. 
             
             Suppose, for the moment, that $ \tau \equiv -1 \pmod{p}$. Then $ G \cong D_{6p}\times\mathcal{C}_q $, so $ \Cay(G;S) $ has a Hamiltonian cycle by Lemma~\ref{lemma 5.8}.
             
     We may now assume that $ \tau \not\equiv -1 \pmod{p} $. Recall that $ \widehat{G} = G/\mathcal{C}_p = \mathcal{C}_{2q} \ltimes \mathcal{C}_3 $. We may assume $ \widehat{a} = a_{2}\gq $ and $ \widehat{b} = a_{2}\gt $. We
     have 
     \begin{align*}
         C_{1} &= ((\widehat{a},\widehat{b},\widehat{a},\widehat{b},\widehat{a}^{-1},\widehat{b},\widehat{a},\widehat{b},\widehat{a}^{-1},\widehat{b},\widehat{a},\widehat{b})^{(q-5)/2},\widehat{a},\widehat{b},\widehat{a}^4,\\&\indent\indent\widehat{b},\widehat{a}^{-3},\widehat{b},\widehat{a}^{-1},\widehat{b},\widehat{a}^2,\widehat{b},\widehat{a}^2,\widehat{b},\widehat{a}^{-1},\widehat{b},\widehat{a}^{-3},\widehat{b},\widehat{a}^4,\widehat{b})
     \end{align*}
     as the first Hamiltonian cycle in $ \Cay(\widehat{G};\widehat{S}) $.
       The picture in Figure~\ref {1Hamcyc13} on page~\pageref{1Hamcyc13} shows the Hamiltonian cycle. We also have
       \begin{align*}
           C_2 = ((\widehat{a},\widehat{b},\widehat{a}^{-1},\widehat{b},\widehat{a},\widehat{b})^{q-5},\widehat{a}^3,\widehat{b},\widehat{a}^2,\widehat{b},\widehat{a}^{-1},\widehat{b},\widehat{a}^{-3},\widehat{b},\widehat{a}^3,\widehat{b},\widehat{a}^{-3},\widehat{b},\widehat{a}^{-1},\widehat{b},\widehat{a}^2,\widehat{b},\widehat{a}^3,\widehat{b})
       \end{align*}
       as the second Hamiltonian cycle in $ \Cay(\widehat{G};\widehat{S}) $. The picture in Figure~\ref{2Hamcyc13} on page~\pageref{2Hamcyc13} shows the Hamiltonian cycle.
    Now we calculate the voltage of $ C_{1} $.
   \begin{align*}
        \mathbb{V}(C_1) &= ((ababa^{-1}b)(aba^{-1}bab))^{(q-5)/2}(aba^4ba^{-3}ba^{-1}ba^2ba^2ba^{-1}ba^{-3}ba^4b) \\ &= ((aa^q\gamma aa^q\gamma a^{-1}a^q\gamma)(aa^q\gamma a^{-1}a^q\gamma aa^q\gamma))^{(q-5)/2} \\&\indent\indent \cdot(aa^q\gamma a^4a^q\gamma a^{-3}a^q\gamma a^{-1}a^q\gamma a^2a^q\gamma a^2a^q\gamma a^{-1}a^q\gamma a^{-3}a^q\gamma a^4a^q\gamma)\\&= ((a^{q+1}\gamma a^{q+1}\gamma a^{q-1}\gamma)(a^{q+1}\gamma a^{q-1}\gamma a^{q+1}\gamma))^{(q-5)/2}\\&\indent\indent\cdot(a^{q+1}\gamma a^{q+4}\gamma a^{q-3}\gamma a^{q-1}\gamma a^{q+2}\gamma a^{q+2}\gamma a^{q-1}\gamma a^{q-3}\gamma a^{q+4}\gamma)\\&= ((\gamma^{\tau^{q+1}+\tau^2+\tau^{q+1}}a^{q+1})(\gamma^{\tau^{q+1}+1+\tau^{q+1}}a^{q+1}))^{(q-5)/2}\\&\indent\indent\cdot(\gamma^{\tau^{q+1}+\tau^5+\tau^{q+2}+\tau+\tau^{q+3}+\tau^5+\tau^{q+4}+\tau+\tau^{q+5}}a^{q+5}) \\&= ((\gamma^{2\tau^{q+1}+\tau^2}a^{q+1})(\gamma^{2\tau^{q+1}+1}a^{q+1}))^{(q-5)/2}\\&\indent\indent\cdot(\gamma^{\tau^{q+5}+\tau^{q+4}+\tau^{q+3}+\tau^{q+2}+\tau^{q+1}+2\tau^5+2\tau}a^{q+5}) \\&= ((\gamma^{2\tau^{q+1}+\tau^2+\tau^{q+1}(2\tau^{q+1}+1)}a^2))^{(q-5)/2}\\&\indent\indent\cdot(\gamma^{\tau^{q+5}+\tau^{q+4}+\tau^{q+3}+\tau^{q+2}+\tau^{q+1}+2\tau^5+2\tau}a^{q+5}) \\&= (\gamma^{3\tau^{q+1}+3\tau^2}a^2)^{(q-5)/2}(\gamma^{\tau^{q+5}+\tau^{q+4}+\tau^{q+3}+\tau^{q+2}+\tau^{q+1}+2\tau^5+2\tau}a^{q+5}) \\&= (\gamma^{(3\tau^{q+1}+3\tau^2)(\tau^{q-5}-1)/(\tau^2-1)}a^{q-5})(\gamma^{\tau^{q+5}+\tau^{q+4}+\tau^{q+3}+\tau^{q+2}+\tau^{q+1}+2\tau^5+2\tau}a^{q+5}) \\&= \gamma^{(3\tau^{q+1}+3\tau^2)(\tau^{q-5}-1)/(\tau^2-1)+\tau^{q-5}(\tau^{q+5}+\tau^{q+4}+\tau^{q+3}+\tau^{q+2}+\tau^{q+1}+2\tau^5+2\tau)}.
    \end{align*}
         \begin{figure}
       \includegraphics[width=0.9\textwidth]{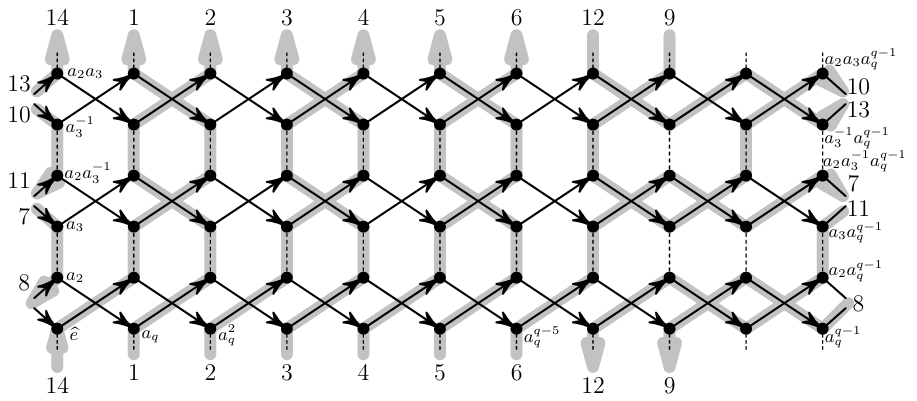}
 \caption{The Hamiltonian cycle $ C_1 $: $ \widehat{a} $ edges are solid and $ \widehat{b} $ edges are dashed.}\label{1Hamcyc13}
\end{figure}Since $ \tau^{2q} \equiv 1 \pmod{p} $, we have $ \tau^q \equiv \pm 1 \pmod{p} $. 

Let us now consider the case where $ \tau^{q}\equiv 1 \pmod{p} $, then by substituting this in the formula for the voltage of $ C_{1} $ we have
    \begin{align*}
        \mathbb{V}(C_1) &= \gamma^{(3\tau+3\tau^2)(\tau^{-5}-1)/(\tau^2-1)+\tau^{-5}(\tau^{5}+\tau^{4}+\tau^{3}+\tau^{2}+\tau+2\tau^5+2\tau)} \\&= \gamma^{3\tau(1+\tau)(\tau^{-5}-1)/(\tau+1)(\tau-1)+(1+\tau^{-1}+\tau^{-2}+\tau^{-3}+\tau^{-4}+2+2\tau^{-4})} \\&= \gamma^{3\tau(\tau^{-5}-1)/(\tau-1)+(3+\tau^{-1}+\tau^{-2}+\tau^{-3}+3\tau^{-4})} \\&= \gamma^{(-2+2\tau^{-3})/(\tau-1)}.
    \end{align*}
    We may assume this does not generate $ \mathcal{C}_p $, then
\begin{align*}
   0 \equiv -2+2\tau^{-3} \pmod{p}.
\end{align*}
 Multiplying by $ \tau^3 $, we have
 \begin{align*}
    0 \equiv -2\tau^3+2 \pmod{p}.
 \end{align*}
   This implies that $ \tau^3 \equiv 1 \pmod{p} $, which contradicts the fact that $ \tau^q \equiv 1 \pmod{p} $ but $ \tau \not\equiv 1 \pmod{p} $.
   
   Now we may assume $ \tau^q \equiv -1 \pmod{p} $, then substituting this in the formula for the voltage of $ C_{1} $ we have 
\begin{align*}
    \mathbb{V}(C_1) &= \gamma^{(-3\tau+3\tau^2)(-\tau^{-5}-1)/(\tau^2-1)-\tau^{-5}(-\tau^{5}-\tau^{4}-\tau^{3}-\tau^{2}-\tau+2\tau^5+2\tau)} \\&= \gamma^{3\tau(\tau-1)(-\tau^{-5}-1)/(\tau+1)(\tau-1)+(1+\tau^{-1}+\tau^{-2}+\tau^{-3}+\tau^{-4}-2-2\tau^{-4})} \\&= \gamma^{3\tau(-\tau^{-5}-1)/(\tau+1)+(-1+\tau^{-1}+\tau^{-2}+\tau^{-3}-\tau^{-4})} \\&= \gamma^{(-4\tau+2\tau^{-1}+2\tau^{-2}-4\tau^{-4})/(\tau+1)}.
\end{align*}
We may assume this does not generate $ \mathcal{C}_p $, then
\begin{align*}
     0 \equiv -4\tau+2\tau^{-1}+2\tau^{-2}-4\tau^{-4}  \pmod{p}.
\end{align*}
Multiplying by $ (-\tau^{4})/2 $, we have
\begin{align*}
     0 &\equiv 2\tau^5-\tau^3-\tau^2+2 \\&= (\tau+1)(2\tau^4-2\tau^3+\tau^2-2\tau+2) \pmod{p}.  
\end{align*}
Since we assumed $ \tau \not \equiv -1 \pmod{p} $, then the above equation implies that 
\begin{align*}
   0 \equiv 2\tau^4-2\tau^3+\tau^2-2\tau+2 \pmod{p} \tag{3A}\label{7.3A}.
\end{align*}
\begin{figure}
      \includegraphics[width=0.9\textwidth]{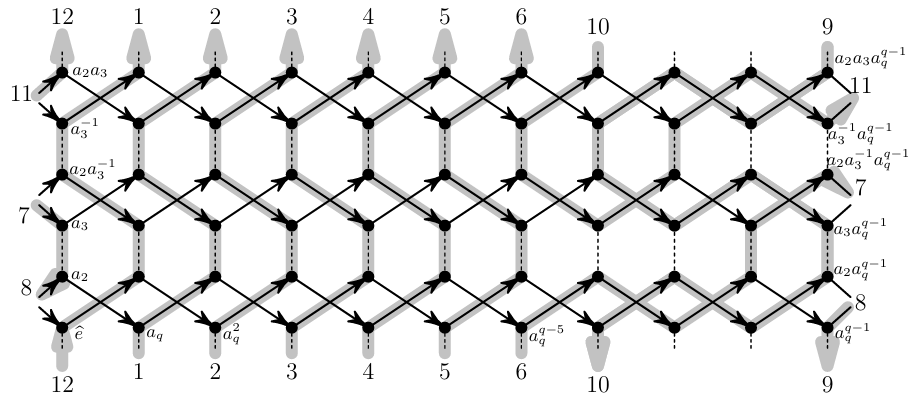}
        \caption{The Hamiltonian cycle $ C_2 $: $ \widehat{a} $ edges are solid and $ \widehat{b} $ edges are dashed.}\label{2Hamcyc13}
    \end{figure}
Now we calculate the voltage of $ C_{2} $.
\begin{align*}
    \mathbb{V}(C_2) &= (aba^{-1}bab)^{(q-5)}(a^3ba^2ba^{-1}ba^{-3}ba^3ba^{-3}ba^{-1}ba^2ba^3b)\\& = (aa^q\gamma a^{-1}a^q\gamma aa^q\gamma)^{(q-5)}(a^3a^q\gamma a^2a^q\gamma a^{-1}a^q\gamma a^{-3}a^q\gamma a^3a^q\gamma a^{-3}a^q\gamma a^{-1}a^q\gamma a^2a^q\gamma a^3a^q\gamma) \\&= (a^{q+1}\gamma a^{q-1}\gamma a^{q+1}\gamma)^{(q-5)}(a^{q+3}\gamma a^{q+2}\gamma a^{q-1}\gamma a^{q-3}\gamma a^{q+3}\gamma a^{q-3}\gamma a^{q-1}\gamma a^{q+2}\gamma a^{q+3}\gamma) \\&= (\gamma^{\tau^{q+1}+1+\tau^{q+1}}a^{q+1})^{(q-5)}(\gamma^{\tau^{q+3}+\tau^5+\tau^{q+4}+\tau+\tau^{q+4}+\tau+\tau^q+\tau^2+\tau^{q+5}}a^{q+5}) \\&= (\gamma^{2\tau^{q+1}+1}a^{q+1})^{(q-5)}(\gamma^{\tau^{q+5}+2\tau^{q+4}+\tau^{q+3}+\tau^q+\tau^5+\tau^2+2\tau}a^{q+5})\\&= (\gamma^{(2\tau^{q+1}+1)((\tau^{q+1})^{(q-5)}-1)/(\tau^{q+1}-1)}a^{(q+1)(q-5)})(\gamma^{\tau^{q+5}+2\tau^{q+4}+\tau^{q+3}+\tau^q+\tau^5+\tau^2+2\tau}a^{q+5}) \\&= \gamma^{(2\tau^{q+1}+1)((\tau^{q+1})^{(q-5)}-1)/(\tau^{q+1}-1)+\tau^{(q+1)(q-5)}(\tau^{q+5}+2\tau^{q+4}+\tau^{q+3}+\tau^q+\tau^5+\tau^2+2\tau)}.
\end{align*}
Since we are assuming $ \tau^q \equiv -1 \pmod{p} $, then by substituting this in the above formula we have 
\begin{align*}
    \mathbb{V}(C_2) &= \gamma^{(-2\tau+1)((-\tau)^{-5}-1)/(-\tau-1)-\tau^{-5}(-\tau^{5}-2\tau^{4}-\tau^{3}-1+\tau^5+\tau^2+2\tau)} \\&= \gamma^{(2\tau^{-4}+2\tau-\tau^{-5}-1)/(-\tau-1)+1+2\tau^{-1}+\tau^{-2}+\tau^{-5}-1-\tau^{-3}-2\tau^{-4}} \\&= \gamma^{(2\tau-3-3\tau^{-1}+3\tau^{-3}+3\tau^{-4}-2\tau^{-5})/(-\tau-1)}.
\end{align*}
 We may assume this does not generate $ \mathcal{C}_p $, then
 \begin{align*}
     2\tau-3-3\tau^{-1}+3\tau^{-3}+3\tau^{-4}-2\tau^{-5} \equiv 0 \pmod{p}.
 \end{align*}
 Multiplying by $ \tau^{5} $, we have
 \begin{align*}
     0 \equiv 2\tau^6-3\tau^5-3\tau^4+3\tau^2+3\tau-2 = (\tau^2-1)(2\tau^4-3\tau^3-\tau^2-3\tau+2) \pmod{p}.
 \end{align*}
 Since $ \tau^2 \not \equiv 1 \pmod{p} $, then the above equation implies that
 \begin{align*}
     0 \equiv 2\tau^4-3\tau^3-\tau^2-3\tau+2 \pmod{p}.
 \end{align*}
  Therefore, by subtracting the above equation from \ref{7.3A}, we have
  \begin{align*}
      0 \equiv (\tau^3+2\tau^2+\tau) = \tau(\tau+1)^2 \pmod{p}.
  \end{align*}
  This is a contradiction.
        \end{case}

      \begin{case}
          Assume none of the previous cases apply. Since $ \langle \overline{a},\overline{b} \rangle = \overline{G} $, we may assume $ |\overline{a}| $ is divisible by $ q $, which means $ |\overline{a}| $ is either $ q $ or $ 2q $. Since Case~\ref{case2.2} applies when $ |\overline{a}| = q $, we must have $ |\overline{a}| = 2q $. Then $ |\overline{b}| = q $, since Cases~\ref{Case 2.1} and~\ref{case 2.3} do not apply. So Case~\ref{case2.2} applies after interchanging $ a $ and $ b $.
  \qedhere
 \end{case}
    \end{proof}

\subsection{Assume \texorpdfstring{$ |S| = 3 $, $ G' = \mathcal{C}_p\times\mathcal{C}_q $ and $ C_{G'}(\mathcal{C}_3) \neq \{e\} $}{Lg}}\hfill\label{3.4}

In this subsection we prove the part of Theorem~\ref{theorem1.1} where, $ |S| = 3 $, $ G' = \mathcal{C}_p\times\mathcal{C}_q $ and $ C_{G'}(\mathcal{C}_3) \neq \{e\} $. Recall $ \overline{G} = G/G' $, $ \widecheck{G} = G/\mathcal{C}_q $ and $ \widehat{G} = G/\mathcal{C}_p $. \label{Section 9}
    \begin{caseprop}\label{prop 3.3}
    Assume
    \begin{itemize}
        \item $ G = (\mathcal{C}_2\times \mathcal{C}_3) \ltimes (\mathcal{C}_p\times\mathcal{C}_q) $,
        \item $ |S| = 3 $,
        \item $ C_{G'}(\mathcal{C}_3) \neq \{e\} $.
    \end{itemize}
    Then $ \Cay(G;S) $ contains a Hamiltonian cycle.
    \end{caseprop}
\begin{proof}
        Let $ S = \{a,b,c\} $. If $  C_{G'}(\mathcal{C}_3) = \mathcal{C}_p\times\mathcal{C}_q  $, then since $ G' \cap Z(G) = \{e\} $ (see Proposition~\ref{Hall Theorem}(\ref{Hall Theorem2})), we conclude that $ C_{G'}(\mathcal{C}_2) = \{e\} $. So we have
           \begin{align*}
               G = \mathcal{C}_3\times(\mathcal{C}_2\ltimes\mathcal{C}_{pq})\cong\mathcal{C}_3\times D_{2pq}.
           \end{align*}
           Therefore, Lemma~\ref{lemma 5.8} applies. 
           
           Since $ C_{G'}(\mathcal{C}_3) \neq \{e\} $, then
        we may assume $ C_{G'}(\mathcal{C}_3) = \mathcal{C}_q $ by interchanging $ q $ and $ p $ if necessary. Since $ G' \cap Z(G) = \{e\} $, then $ \mathcal{C}_2 $ inverts $ \mathcal{C}_q $. Since $ \mathcal{C}_3 $ centralizes $ \mathcal{C}_q $ and $ Z(G) \cap G' = \{e\} $ (by Proposition~\ref{Hall Theorem}(\ref{Hall Theorem2})), then $ \mathcal{C}_2 $ inverts $ \mathcal{C}_q $. Thus,
        \begin{align*}
            \widehat{G} = (\mathcal{C}_2\times\mathcal{C}_3)\ltimes \mathcal{C}_q \cong (\mathcal{C}_2 \ltimes \mathcal{C}_q) \times \mathcal{C}_3 = D_{2q} \times \mathcal{C}_3.
        \end{align*}
          Now if $ \widehat{S} $ is minimal, then Lemma~\ref{lemma5.10} applies. Therefore, we may assume $ \widehat{S} $ is not minimal. Choose a 2-element subset $ \{a,b\} $ of $ S $ that generates $ \widehat{G} $. From the minimality of $ S $, we see that $ \langle a, b \rangle = D_{2q} \times \mathcal{C}_3 $ after replacing $ a $ and $ b $ by conjugates. The projection of $ (a,b) $ to $ D_{2q} $ must be of the form $ (a_{2}, \gq) $ or $ (a_{2}, a_{2}\gq) $, where $ a_{2} $ is reflection and $ \gq $ is a rotation.~(Also note that $ \widehat{b} \neq \gq  $ because $ S \cap G' = \emptyset $ by Assumption~\ref{assumption 3.1}(\ref{assumption 3.1.6}).) Therefore, $ (a,b) $ must have one of the following forms:
        \begin{enumerate}
            \item $ (a_2,\gt\gq) $,
            \item $ (a_2,a_2\gt\gq) $,
            \item $ (a_2\gt,a_2\gq) $,
            \item $ (a_{2}\gt,\gt\gq) $,
            \item $ (a_{2}\gt,a_2\gt\gq) $.
        \end{enumerate}
        
        Let $ c $ be the third element of $ S $. We may write $ c = a_{2}^{i}\gt^{j}\gq^{k}\gamma_{p} $ with $ 0\leq i \leq 1 $, $ 0\leq j \leq 2 $ and $ 0\leq k \leq q-1 $. Note since $ S \cap G' = \emptyset $, we know that $ i $ and $ j $ cannot both be equal to $ 0 $. Additionally, we have $ \gt\gamma_{p} \gt^{-1} = \gamma_{p}^{\widehat{\tau}} $ where $ \widehat{\tau}^{3} \equiv 1  \pmod{\mathcal{C}_p}$. Also, $ \widehat{\tau} \not \equiv 1 \pmod{p}$ since $ C_{G'}(\mathcal{C}_3) = \mathcal{C}_q $. Therefore, we conclude that $ \widehat{\tau}^2+\widehat{\tau}+1 \equiv 0 \pmod{p} $. Note that this implies $ \widehat{\tau} \not \equiv -1 \pmod{p} $.
        \setcounter{case}{0}
        \begin{case}
        Assume $ a = a_{2} $ and $ b = \gt\gq $.
         \begin{subcase}
            Assume $ i \neq 0 $. Then, $ c = a_{2}\gt^j\gq^k\gamma_{p} $. Thus, by Lemma~\ref{lemma 5.14}(\ref{lemma 5.14.1}) $ \langle b,c \rangle = G $ which contradicts the minimality of $ S $.
        \end{subcase}
        \begin{subcase}
            Assume $ i = 0 $. Then $ j \neq 0 $. We may assume $ j = 1 $, by replacing $ c $ with $ c^{-1} $ if necessary. Thus $ c = \gt\gq^{k}\gamma_{p} $. Consider $ \overline{G} = \mathcal{C}_2\times\mathcal{C}_3 $. We have $ \overline{a} = a_{2} $, $ \overline{b} = \gt $ and $ \overline{c} = \gt $. Therefore, $ \overline{b} = \overline{c} = \gt $. We have $ (\overline{a},\overline{b}^{2},\overline{a},\overline{b}^{-2}) $ as a Hamiltonian cycle in $ \Cay(\overline{G};\overline{S}) $. Since we can replace each $ \overline{b}$ by $ \overline{c} $, then we consider $ C_{1} = (\overline{a},\overline{b}^{2},\overline{a},\overline{b}^{-1},\overline{c}^{-1}) $ and $ C_{2} =  (\overline{a},\overline{b}^{2},\overline{a},\overline{c}^{-2}) $ as Hamiltonian cycles in $ \Cay(\overline{G};\overline{S}) $. Now since there is one occurrence of $ c $ in $ C_1 $, then by Lemma~\ref{lemma 2.5.2} the subgroup generated by $ \mathbb{V}(C_1) $ contains $ \mathcal{C}_p $. Also,
        \begin{align*}
          \mathbb{V}(C_1) &= ab^2ab^{-1}c^{-1} \\&\equiv a_{2}\cdot\gt^{2}\gq^{2}\cdot a_{2}\cdot\gq^{-1}\gt^{-1}\cdot\gq^{-k}\gt^{-1} \pmod{\mathcal{C}_p} \\ &= \gq^{-2}\gt\gq^{-1-k}\gt^{-1} \\ &= \gq^{-3-k}.
        \end{align*}
         We can assume this does not generate $ \mathcal{C}_q $, for otherwise Factor Group Lemma~\ref{FGL} applies. Therefore,
         \begin{align*}
             -3-k \equiv 0 \pmod{q}.
         \end{align*}
         Thus, $ k \equiv -3 \pmod{q} $.
         
         Now we calculate the voltage of $ C_{2} $.
         \begin{align*}
             \mathbb{V}(C_2) &= ab^2ac^{-2} \\&\equiv a_{2}\cdot\gt^2\cdot a_{2}\cdot\gamma_{p}^{-1}\gt^{-1}\gamma_{p}^{-1}\gt^{-1} \pmod{\mathcal{C}_q} \\ &= \gt^2\gamma_{p}^{-1}\gt^{-1}\gamma_{p}^{-1}\gt^{-1} \\ &= \gamma_{p}^{-\widehat{\tau}^2-\widehat{\tau}}.
        \end{align*}
         Since $ \widehat{\tau}^2+\widehat{\tau}+1 \equiv 0 \pmod{p} $, then $ -\widehat{\tau}^2-\widehat{\tau} \equiv 1 \pmod{p} $. Thus, $ \gamma_{p}^{-\widehat{\tau}^2-\widehat{\tau}} = \gamma_{p} $ generates~$ \mathcal{C}_p $.
        \begin{align*}
            \mathbb{V}(C_2) &= ab^2ac^{-2} \\&\equiv a_{2}\cdot\gt^2\gq^2\cdot a_{2}\cdot\gq^{-k}\gt^{-1}\gq^{-k}\gt^{-1} \pmod{\mathcal{C}_p} \\ &= \gq^{-2}\gt^2\gq^{-k}\gt^{-1}\gq^{-k}\gt^{-1}\\ &= \gq^{-2(k+1)}.
        \end{align*}
        We know $ k \equiv -3 \pmod{q} $, therefore, $ -2(k+1) \equiv 4 \pmod{q} $, so Factor Group Lemma~\ref{FGL} applies.
        \end{subcase}
        \end{case}
        \begin{case}\label{Case 9.2}
        Assume $ a = a_{2} $ and $ b = a_{2}\gt\gq $.
        \begin{subcase}
            Assume $ i = 0 $, then $ j \neq 0 $. If $ k \neq 0 $, then $ c = \gt^{j}\gq^{k}\gamma_{p} $. Thus, by Lemma~\ref{lemma 5.14}(\ref{lemma 5.14.4}) $ \langle b,c \rangle = G $ which contradicts the minimality of $ S $. Therefore, we may assume $ k = 0 $. We may also assume $ j = 1 $, by replacing $ c $ with $ c^{-1} $ if necessary. Then $ c = \gt\gamma_{p} $.
            
            Consider $ \overline{G} = \mathcal{C}_2\times\mathcal{C}_3 $, thus $ \overline{a} = a_{2} $, $ \overline{b} = a_{2}\gt $ and $ \overline{c} = \gt $. Therefore, $ |\overline{a}| = 2 $, $ |\overline{b}| = 6 $ and $ |\overline{c}| = 3 $. Consider $ C =  (\overline{b}^{2},\overline{c},\overline{b},\overline{c}^{-1},\overline{a}) $ as a Hamiltonian cycle in $ \Cay(\overline{G};\overline{S}) $. Now we calculate its voltage.
       \begin{align*}
            \mathbb{V}(C) &= b^2cbc^{-1}a \\&\equiv a_{2}\gt\gq a_{2}\gt\gq\cdot\gt\cdot a_{2}\gt\gq\cdot\gt^{-1}\cdot a_{2} \pmod{\mathcal{C}_p} \\ &=  \gq^{-1} 
       \end{align*}
       which generates $ \mathcal{C}_q $. By considering the fact that $ \mathcal{C}_2 $ might centralize $ \mathcal{C}_p $ or not, we have
       \begin{align*}
           \mathbb{V}(C) &= b^2cbc^{-1}a \\&\equiv a_{2}\gt a_{2}\gt\cdot\gt\gamma_{p}\cdot a_{2}\gt\cdot\gamma_{p}^{-1}\gt^{-1}\cdot a_{2} \pmod{\mathcal{C}_q} \\ &= \gamma_{p}\gt\gamma_{p}^{\mp1}\gt^{-1} \\ &= \gamma_{p}^{1\mp\widehat{\tau}}.
       \end{align*}
       which generates $ \mathcal{C}_p $. Therefore, the subgroup generated by $ \mathbb{V}(C) $ is $ G' $. So, Factor Group Lemma~\ref{FGL} applies.
        \end{subcase}
        \begin{subcase}
                Assume $ j = 0 $. Then $ i \neq 0 $. If $ k \neq 1 $, then $ c = a_{2}\gq^{k}\gamma_{p} $. Thus, by Lemma~\ref{lemma 5.14}(\ref{lemma 5.14.5}) $ \langle b,c \rangle = G $ which contradicts the minimality of $ S $. We may therefore assume $ k = 1 $. Then $ c = a_{2}\gq\gamma_{p} $. 
                
                Consider $ \overline{G} = \mathcal{C}_2\times\mathcal{C}_3 $, then $ \overline{a} $  = $ \overline{c} = a_{2} $ and $ \overline{b} = a_{2}\gt $. Thus, $ |\overline{a}| = |\overline{c}| = 2 $ and $ |\overline{b}| = 6 $. We have $ C =  (\overline{b}^{2},\overline{c},\overline{b}^{-2},\overline{a}) $ as a Hamiltonian cycle in $ \Cay(\overline{G};\overline{S}) $. Since there is one occurrence of $ c $ in $ C $, and it is the only generator of $ G $ that contains $ \gamma_p $, then by Lemma~\ref{lemma 2.5.2} we conclude that the subgroup generated by $ \mathbb{V}(C) $ contains $ \mathcal{C}_p $. Also,
        \begin{align*}
            \mathbb{V}(C) &= b^2cb^{-2}a \\&\equiv a_{2}\gt\gq a_{2}\gt\gq\cdot a_{2}\gq\cdot\gq^{-1}\gt^{-1}a_{2}\gq^{-1}\gt^{-1}a_{2}\cdot a_{2} \pmod{\mathcal{C}_p} \\ &= \gq^{-1}\gt\gq\gt\gq^{-1}\gt^{-1}\gq\gt^{-1}\gq^{-1} \\ &= \gq^{-1}.
        \end{align*}
        which generates $ \mathcal{C}_q $.  Therefore, the subgroup generated by $ \mathbb{V}(C) $ is $ G' $. So, Factor Group Lemma~\ref{FGL} applies.
        \end{subcase}
        \begin{subcase}
                    Assume $ i \neq 0 $ and $ j \neq 0 $. We may assume $ j = 1 $, by replacing $ c $ with $ c^{-1} $ if necessary. So $ c = a_{2}\gt\gq^{k}\gamma_{p} $. If $ k \neq 1 $, then by Lemma~\ref{lemma 5.14}(\ref{lemma 5.14.5}) $ \langle b,c \rangle = G $ which contradicts the minimality of $ S $. We may now assume $ k = 1 $. Then $ c = a_{2}\gt\gq\gamma_{p} $.
                    
                    Consider $ \overline{G} = \mathcal{C}_2\times\mathcal{C}_3 $. Then $ \overline{a} = a_{2} $ and $ \overline{b} = \overline{c} = a_{2}\gt $. Therefore, $ |\overline{b}| = |\overline{c}| = 6 $ and $ |\overline{a}| = 2 $. We have $ C = (\overline{c},\overline{a},(\overline{b},\overline{a})^2) $ as a Hamiltonian cycle in $ \Cay(\overline{G};\overline{S}) $. Since there is one occurrence of $ c $ in $ C $, and it is the only generator of $ G $ that contains $ \gamma_p $, then by Lemma~\ref{lemma 2.5.2} we conclude that the subgroup generated by $ \mathbb{V}(C) $ is $ \mathcal{C}_p $. Also,
        \begin{align*}
           \mathbb{V}(C) &= ca(ba)^2 \\&\equiv a_{2}\gt\gq\cdot a_{2}\cdot a_{2}\gt\gq\cdot a_{2}\cdot a_{2}\gt\gq\cdot a_{2} \pmod{\mathcal{C}_p} \\ &= \gt\gq^{-2}\gt\gq^{-1}\gt\\ &= \gq^{-3} 
        \end{align*}
        which generates $ \mathcal{C}_q $. Therefore, the subgroup generated by $ \mathbb{V}(C) $ is $ G' $. So, Factor Group Lemma~\ref{FGL} applies.
        \end{subcase}
        \end{case}
        \begin{case}
            Assume $ a = a_2\gt $ and $ b = a_2\gq $. Since $ b = a_2\gq $ is conjugate to $ a_2 $ via an element of $ \mathcal{C}_q $ (which centralizes $ \mathcal{C}_3 $), then $ \{a,b\} $ is conjugate to $ \{a_2\gt\gq^m,a_2\} $ for some nonzero $ m $. So Case~\ref{Case 9.2} applies (after replacing $ \gq $ with $ \gq^m $).
        \end{case}
        \begin{case}
        Assume $ a = a_{2}\gt $ and $ b = \gt\gq $.
        \begin{subcase}
               Assume $ i \neq 0 $. Then $ c = a_{2}\gt^j\gq^{k}\gamma_{p} $. Thus, by Lemma~\ref{lemma 5.14}(\ref{lemma 5.14.1}) $ \langle b,c \rangle = G $ which contradicts the minimality of $ S $.
        \end{subcase}
        \begin{subcase}
                Assume $ i = 0 $. Then $ j \neq 0 $ and $ c = \gt^{j}\gq^{k}\gamma_{p} $. If $ k \neq 0 $, then by Lemma~\ref{lemma 5.14}(\ref{lemma 5.14.3}) $ \langle a,c \rangle = G $ which contradicts the minimality of $ S $. So we may assume $ k = 0 $. We may also assume $ j = 1 $, by replacing $ c $ with $ c^{-1} $ if necessary. Then $ c = \gt\gamma_{p} $.
                
                Consider $ \overline{G} = \mathcal{C}_2\times\mathcal{C}_3 $. Therefore, $ \overline{a} = a_{2}\gt $ and $ \overline{b} = \overline{c} = \gt $. In addition, $ |\overline{a}| = 6 $ and $ |\overline{b}| = |\overline{c}| = 3 $. We have $ C =  (\overline{c},\overline{b},\overline{a},\overline{b}^{-2},\overline{a}^{-1}) $ as a Hamiltonian cycle in $ \Cay(\overline{G};\overline{S}) $. Since there is one occurrence of $ c $ in $ C $, and it is the only generator of $ G $ that contains $ \gamma_p $, then by Lemma~\ref{lemma 2.5.2} we conclude that the subgroup generated by $ \mathbb{V}(C) $ contains $ \mathcal{C}_p $. Also,
                \begin{align*}
                    \mathbb{V}(C) &= cbab^{-2}a^{-1} \\&\equiv \gt\cdot\gt\gq\cdot a_{2}\gt\cdot\gq^{-2}\gt^{-2}\cdot\gt^{-1}a_{2} \pmod{\mathcal{C}_p} \\ &= \gt\gq\gt^2\gq^2 \\ &= \gq^3 
                \end{align*}
                which generates $ \mathcal{C}_q $. Therefore, the subgroup generated by $ \mathbb{V}(C) $ is $ G' $. Thus, Factor Group Lemma~\ref{FGL} applies.
        \end{subcase}
        \end{case}
        \begin{case}
        Assume $ a = a_{2}\gt $, $ b = a_{2}\gt\gq $.
        \begin{subcase}
            Assume $ i = 0 $. Then $ j \neq 0 $ and $ c = \gt^j \gq^k \gamma_{p} $. If $ k \neq 0 $, then by Lemma~\ref{lemma 5.14}(\ref{lemma 5.14.4}) $ \langle b,c \rangle = G $ which contradicts the minimality of $ S $. So we may assume $ k = 0 $. We may also assume $ j = 1 $, by replacing $ c $ with $ c^{-1} $ if necessary. Then $ c = \gt\gamma_{p} $. 
            
            Consider $ \overline{G} = \mathcal{C}_2\times\mathcal{C}_3 $. Therefore, $ \overline{a} = \overline{b} = a_{2}\gt $ and $ \overline{c} = \gt $. Thus, $ |\overline{a}| = |\overline{b}| = 6 $ and $ |\overline{c}| = 3 $. We have $ C = (\overline{a},\overline{c}^2,\overline{b}^{-1},\overline{c}^{-2}) $ as a Hamiltonian cycle in $ \Cay(\overline{G};\overline{S}) $. Now we calculate its voltage.
    \begin{align*}
       \mathbb{V}(C) &= ac^2b^{-1}c^{-2} \\&\equiv a_{2}\gt\cdot\gt^2\cdot\gq^{-1}\gt^{-1}a_{2}\cdot\gt^{-2} \pmod{\mathcal{C}_p} \\ &= \gt^{-1}\gq\gt^{-2}  \\ &= \gq
    \end{align*}
    which generates $ \mathcal{C}_q $. Also
    \begin{align*}
        \mathbb{V}(C) &= ac^2b^{-1}c^{-2} \\&\equiv ac^2 a^{-1}c^{-2} \pmod{\mathcal{C}_q } \text{ } (\text{because }  a \equiv b \pmod{\mathcal{C}_q}) \\&= ac^{-1}a^{-1}c \text{ }(\text{because } |c| = 3)\\& = [a,c^{-1}].
    \end{align*}
    This generates $ \mathcal{C}_p $, because $ \{a,c\} $ generates $ G/\mathcal{C}_q $. Therefore, the subgroup generated by $ \mathbb{V}(C) $ is $ G' $. So, Factor Group Lemma~\ref{FGL} applies.
        \end{subcase}
        \begin{subcase}
        Assume $ i \neq 0 $. Then $ c = a_{2}\gt^j\gq^k\gamma_{p} $. If  $ k \neq 1 $, then by Lemma~\ref{lemma 5.14}(\ref{lemma 5.14.5}) $ \langle b,c \rangle = G $ which contradicts the minimality of $ S $. So we may assume $ k = 1 $. Then $ c = a_{2}\gt^j\gq\gamma_{p} $. We show that $ \langle a,c \rangle = G $. Now, we have
        \begin{align*}
            \langle a,c \rangle &= \langle a_2,\gt,c \rangle \text{ } (\text{because } \langle a \rangle = \langle a_2\gt \rangle = \langle a_2,\gt \rangle ) \\ &= \langle a_2,\gt,a_{2}\gt^j\gq\gamma_{p} \rangle \\&= \langle a_2,\gt,\gq\gamma_p \rangle \\&= \langle a_2,\gt,\gq,\gamma_p \rangle \\&= G,
        \end{align*}
        which contradicts the minimality of $ S $.
        \qedhere
        \end{subcase}
        \end{case}
\end{proof}

 \subsection{Assume \texorpdfstring{$ |S| = 3 $, $ G' = \mathcal{C}_p\times \mathcal{C}_q $ and $ \widehat{S} $ is minimal}{Lg}}\hfill\label{3.3}
    
 In this subsection we prove the part of Theorem~\ref{theorem1.1} where, $ |S| = 3 $, $ G' = \mathcal{C}_{p} \times \mathcal{C}_q $ and $ C_{G'}(\mathcal{C}_3) = \{e\} $. Recall $ \overline{G} = G/G' $ and $ \widehat{G} = G/\mathcal{C}_p $.
           \begin{caseprop}\label{caseprop 3.4}
           Assume
           \begin{itemize}
               \item $ G = (\mathcal{C}_2\times\mathcal{C}_3)\ltimes(\mathcal{C}_p\times\mathcal{C}_q) $,
               \item $ |S| = 3 $,
               \item $ \widehat{S} $ is minimal.
           \end{itemize}
           Then $ \Cay(G;S) $ contains a Hamiltonian cycle.
           \end{caseprop}
    \begin{proof}
           
          Let $ S = \{a,b,c\} $. If $ C_{G'}(\mathcal{C}_3) \neq \{e\} $, then Proposition~\ref{prop 3.3} applies. Hence we may assume $ C_{G'}(\mathcal{C}_3) = \{e\} $. Then we have four different cases.
          \setcounter{case}{0}
             \begin{case} \label{Subcase 8.2.1}
             Assume $  C_{G'}(\mathcal{C}_2) = \mathcal{C}_p\times\mathcal{C}_q  $, thus $ G = \mathcal{C}_2 \times (\mathcal{C}_3 \ltimes \mathcal{C}_{pq}) $. Since $ \widehat{S} $ is minimal, then all three elements belonging to $ \widehat{S} $ must have prime order. There is an element $ \widehat{a} \in \widehat{S} $, such that $ |\widehat{a}| = 2 $, otherwise all elements of $ S $ belong to a subgroup of index $ 2 $ of $ G $, so $ \langle a,b,c \rangle \neq G $ which is a contradiction. If $ |a| = 2p $, then Corollary~\ref{corollary 5.2} applies with $ s = a $ and $ t = a^{-1} $, because there is a Hamiltonian cycle in $ \Cay(\widehat{G};\widehat{S}) $~(see Theorem~\ref{theorem 1.2}(\ref{theorem 1.2.3})) which uses at least one labeled edge $ \widehat{a} $ because $ \widehat{S} $ is minimal.
             
             Now we may assume $ |a| = 2 $. Replacing $ a $ by a conjugate we may assume $ \langle a \rangle = \mathcal{C}_2 $. Thus, $ \langle b,c\rangle = \mathcal{C}_3 \ltimes \mathcal{C}_{pq} $. By Theorem~\ref{theorem 1.2}(\ref{theorem 1.2.3}), there is a Hamiltonian path $ L $ in $ \Cay(\mathcal{C}_3 \ltimes \mathcal{C}_{pq},\{b,c\}) $. Therefore, $ LaL^{-1}a^{-1} $ is a Hamiltonian cycle in $ \Cay(G;S) $.
             \end{case}
             \begin{case} \label{subcase 8.2.2}
             Assume $ C_{G'}(\mathcal{C}_2) = \mathcal{C}_q $. Therefore,
             \begin{align*}
                 \widehat{G} = G/\mathcal{C}_p = \mathcal{C}_{6} \ltimes \mathcal{C}_q \cong \mathcal{C}_2 \times (\mathcal{C}_3 \ltimes \mathcal{C}_q).
             \end{align*}
              There is some $ a \in S $ such that $ |\widehat{a}| = 2 $. Thus, we can assume $ |a| = 2 $, for otherwise Corollary~\ref{corollary 5.2} applies with $ s = a $ and $ t = a^{-1} $.~(Note since $ \widehat{S} $ is minimal, then $ \widehat{a} $ must be used in any Hamiltonian cycle in $ \Cay(\widehat{G};\widehat{S}) $.) We may assume $ a = a_{2} $. Since $ \widehat{S} $ is minimal, $ S \cap G' = \emptyset $~(see Assumption~\ref{assumption 3.1}(\ref{assumption 3.1.6})) and each element belonging to $ \widehat{S} $ has prime order, then $ |\widehat{b}| = |\widehat{c}| = 3 $. We may assume $ \widehat{a} = a_{2} $, $ \widehat{b} = \gt $ and $ \widehat{c} = \gt\gq $. We have the following two Hamiltonian paths in $ \Cay(\mathcal{C}_3\ltimes \mathcal{C}_q;\{\widehat{b},\widehat{c}\}) $:
             \begin{align*}
                 L_1 = ((\widehat{c},\widehat{b}^2)^{q-1},\widehat{c}, \widehat{b})
             \end{align*}
             and
             \begin{align*}
                 L_2 = ((\widehat{b},\widehat{c},\widehat{b})^{q-1},\widehat{b}, \widehat{c}).
             \end{align*}
             These lead to the following two Hamiltonian cycles in $ \Cay(\widehat{G};\widehat{S}) $:
             \begin{align*}
                 C_{1} = (L_1,\widehat{a},L_1^{-1},\widehat{a})
             \end{align*}
             and
             \begin{align*}
                 C_{2} = (L_2,\widehat{a},L_2^{-1},\widehat{a}).
             \end{align*}
             Then if we let 
             \begin{align*}
                 \prod L_1 = (cb^2)^{q-1}cb = (cb^2)^qb^{-1} \in \gt^{-1}\mathcal{C}_p
             \end{align*}
              and 
             \begin{align*}
                 \prod L_2 = (bcb)^{q-1}bc = (bcb)^qb^{-1} = b(cb^2)^qb^{-2} = b(\prod L_1)b^{-1}
             \end{align*}
             then it is clear that $ V(C_i) = [\prod L_i,a] $ for $ i = 1,2 $. Therefore, we may assume $ a $ centralizes $ \prod L_1 $ and $ \prod L_2 $, for otherwise Factor Group Lemma~\ref{FGL} applies. Now, since $ a $ centralizes $ \prod L_1 $, and $ \prod L_1 \in \gt^{-1}\mathcal{C}_p $, we must have $ \prod L_1 = \gt^{-1} $. So $ \prod L_2 = b\gt^{-1}b^{-1} $. If $ b $ does not centralize $ \gt $, then $ \mathbb{V}(C_1) \neq \mathbb{V}(C_2) $, so the voltage of $ C_1 $ or $ C_2 $ cannot both be equal to identity. Therefore, Factor Group Lemma~\ref{FGL} applies. Now if $ b $ centralizes $ \gt $, then we can assume $ b = \gt $. Therefore, $ c = \gt\gq\gamma_p $. We calculate the voltage of $ C_1 $. We have
             \begin{align*}
                 \mathbb{V}(C_1) &= (cb^2)^qb^{-1}a((cb^2)^qb^{-1})^{-1}a \\&= (\gt\gq\gamma_p\cdot\gt^2)^q\cdot\gt^{-1}\cdot a_{2}\cdot((\gt\gq\gamma_p\cdot\gt^2)^q\cdot\gt^{-1})^{-1}\cdot a_{2} \\&= (\gt\gq\gamma_p\gt^{-1})^q\gt^{-1}a_2((\gt\gq\gamma_p\gt^{-1})\gt^{-1})^{-1}a_2\\&= \gt\gq^q\gamma_p^q\gt^{-1}\gt^{-1}a_2(\gt\gq^q\gamma_p^q\gt^{-1}\gt^{-1})^{-1}a_2\\&= \gt\gamma_p^q\gt^{-2}a_2(\gt\gamma_p^q\gt^{-2})^{-1}a_2\\&= \gt\gamma_p^q\gt^{-2}a_2\gt^2\gamma_p^{-q}\gt^{-1}a_2\\&= \gt\gamma_p^{2q}\gt^{-1}
             \end{align*}
             which generates $ \mathcal{C}_p $. Thus, Factor Group Lemma~\ref{FGL} applies.
             \end{case}
            \begin{case} \label{Subcase 8.2.3}
                Assume $ C_{G'}(\mathcal{C}_2) = \mathcal{C}_p $. Therefore,
                \begin{align*}
                    \widecheck{G} = G/\mathcal{C}_q = \mathcal{C}_{6}\ltimes\mathcal{C}_p \cong \mathcal{C}_2 \times (\mathcal{C}_3\ltimes\mathcal{C}_p).
                \end{align*}
                Now since $ S \cap G' = \emptyset $~(see Assumption~\ref{assumption 3.1}(\ref{assumption 3.1.6})) and $ \mathcal{C}_3 $ does not centralize $ \mathcal{C}_p $, then for all $ a \in S $, we have $ |\widecheck{a}| \in \{2,3,6,2p\} $. If $ |\widecheck{a}| = 6 $, then $ |\widehat{a}| $ is divisible by $ 6 $ which contradicts the minimality of $ \widehat{S} $.~(Note that every element belong to $ \widehat{S} $ has prime order.)~If $ |\widecheck{a}| = 2p $, then $ |\widehat{a}| = 2 $ (because $ \widehat{S} $ is minimal). Therefore, Corollary~\ref{corollary 5.2} applies with $ s = a $ and $ t = a^{-1} $~(Note that since $ \widehat{S}$ is minimal, then there is a Hamiltonian cycle in $ \Cay(\widehat{G};\widehat{S}) $ uses at least one labeled edge $ \widehat{a} $.)~Thus, $ |\widecheck{a}| \in \{2,3\} $ for all $ a \in S $. This implies that $ \widecheck{S} $ is minimal, because we need an $ a_2 $ and an $ \gt $ to generate $ \mathcal{C}_2\times\mathcal{C}_3 $ and two elements whose order divisible by $ 2 $ or $ 3 $ to generate $ \mathcal{C}_p $. So by interchanging $ p $ and $ q $ the proof in Case~\ref{subcase 8.2.2} applies.
            \end{case}
            \begin{case} \label{Case 3.3.4}
                Assume $ C_{G'}(\mathcal{C}_2) = \{e\} $. Consider
                \begin{align*}
                    \widehat{G} = G/\mathcal{C}_p = (\mathcal{C}_2\times\mathcal{C}_3)\ltimes\mathcal{C}_q.
                \end{align*}
                Now since $ \widehat{S} $ is minimal, every element of $ \widehat{S} $ has prime order. Since $ S\cap G' = \emptyset $~(see Assumption~\ref{assumption 3.1}(\ref{assumption 3.1.6})), then for every $ \widehat{s} \in \widehat{S} $, we have $ |\widehat{s}| \in \{2,3\} $. Since $ C_{G'}(\mathcal{C}_2) = \{e\} $ and $ C_{G'}(\mathcal{C}_3) = \{e\} $, this implies that for every $ s \in S $, we have $ |s| \in \{2,3\} $. From our assumption we know that $ S = \{a,b,c\} $. Now we may assume $ |a| = 2 $ and $ |b| = 3 $. Also, we know that $ |c| \in \{2,3\} $.
                
                If $ |c| = 2 $, then $ c = a\gamma $, where $ \gamma \in G' $. Suppose, for the moment, $ \langle \gamma \rangle \neq G' $. Since $ \langle \gamma \rangle \triangleleft G $, then we have
                \begin{align*}
                    G = \langle a,b,c \rangle = \langle a,b,\gamma \rangle = \langle a,b \rangle \langle \gamma \rangle.
                \end{align*}
                Now since $ \widehat{S} $ is minimal, $ \langle a,b \rangle $ does not contain $ \mathcal{C}_q $. So this implies that $ \langle \gamma \rangle $ contains $ \mathcal{C}_q $. Since $ \langle \gamma \rangle $ does not contain $ G' $, then $ \langle \gamma \rangle = \mathcal{C}_q $. Thus, we may assume that $ a = a_2 $~(by conjugation if necessary), $ b = \gt\gamma_p $ and $ c = a_2\gq $. So $ \langle b,c \rangle = \langle \gt\gamma_p,a_2\gq \rangle = G $~(since $ \gt\gamma_p $ and $ a_2\gq $ clearly generate $ \overline{G} $ and do not commute modulo $ p $ or modulo $ q $, they must generate $ G $). This contradicts the minimality of $ S $. Therefore, $ \langle \gamma \rangle = G' $.
                
                Consider $ \overline{G} = \mathcal{C}_2\times\mathcal{C}_3 $. Then $ \overline{a} = \overline{c} $. We have $ |\overline{a}| = |\overline{c}| = 2 $ and $ |\overline{b}| = 3 $. We also have $ C_1 = (\overline{c}^{-1},\overline{b}^{-2},\overline{a},\overline{b}^2) $ as a Hamiltonian cycle in $ \Cay(\overline{G};\overline{S}) $. Now we calculate~its~voltage. 
                \begin{align*}
                    \mathbb{V}(C_1) = c^{-1}b^{-2}ab^2 = \gamma^{-1}a^{-1}b^{-2}ab^2.
                \end{align*}
                Now, $ a^{-1}b^{-2}ab^2 \in G' $. Since $ \langle a,b \rangle \neq G $, we have $ a^{-1}b^{-2}ab^2 \in \{e,\gamma_p\} $. If $ a^{-1}b^{-2}ab^2 = e $, then $ a $ and $ b^2 $ commute, so $ a $ and $ b $ commute. Hence $ b = \gt $, so $ \langle b,c \rangle = G $, a contradiction. So $ a^{-1}b^{-2}ab^2 = \gamma_p $, and $ \mathbb{V}(C_1) = \gamma^{-1}\gamma_p $ which generates $ G' $. Therefore, Factor Group Lemma~\ref{FGL} applies.
                
                Now we can assume $ |c| = 3 $. Then $ c = b\gamma $, where $ \gamma \in G' $~(after replacing $ c $ with its inverse if necessary). Suppose, for the moment, $ \langle \gamma \rangle \neq G' $. Since $ \langle \gamma \rangle \triangleleft G $, then~we~have
                \begin{align*}
                    G = \langle a,b,c \rangle = \langle a,b,\gamma \rangle = \langle a,b \rangle \langle \gamma \rangle.
                \end{align*}
                Now since $ \widehat{S} $ is minimal, then $ \langle a,b \rangle $ does not contain $ \mathcal{C}_q $. So this implies that $ \langle \gamma \rangle $ contains $ \mathcal{C}_q $. Since $ \langle \gamma \rangle $ does not contain $ G' $, then $ \langle \gamma \rangle = \mathcal{C}_q $. Therefore, we may assume that $ a = a_2\gamma_p $~(by conjugation if necessary), $ b = \gt $ and $ c = \gt\gq $. So $ \langle a,c \rangle = \langle a_2\gamma_p,\gt\gq \rangle = G $~(since $ a_2\gamma_p $ and $ \gt\gq $ clearly generate $ \overline{G} $ and do not commute modulo $ p $ or modulo $ q $, they must generate $ G $). This contradicts the minimality of $ S $. So $ \langle \gamma \rangle = G' $.
                
                 Consider $ \overline{G} = \mathcal{C}_2\times\mathcal{C}_3 $. Then $ \overline{b} = \overline{c} $. We have $ |\overline{a}| = 2 $ and $ |\overline{b}| = |\overline{c}| = 3 $. We also have $ C_2 = (\overline{c}^{-1},\overline{b}^{-1},\overline{a}^{-1},\overline{b}^2,\overline{a}) $ as a Hamiltonian cycle in $ \Cay(\overline{G};\overline{S}) $. Now we calculate its voltage.
                \begin{align*}
                    \mathbb{V}(C_2) = c^{-1}b^{-1}a^{-1}b^2a = \gamma^{-1}b^{-1}b^{-1}a^{-1}b^2a.
                \end{align*}
                Now, $ b^{-2}a^{-1}b^2a \in G' $. Since $ \langle a,b \rangle \neq G $, we have $ b^{-2}a^{-1}b^2a \in \{e,\gamma_p\} $. If $ b^{-2}a^{-1}b^2a = e $, then $ a $ and $ b^2 $ commute, so $ a $ and $ b $ commute. Hence $ a = a_2 $, so $ \langle a,c \rangle = G $, a contradiction. So $ b^{-2}a^{-1}b^2a = \gamma_p $, and $ \mathbb{V}(C_2) = \gamma^{-1}\gamma_p $ which generates $ G' $. Therefore, Factor Group Lemma~\ref{FGL} applies.
                 \qedhere
            \end{case}
\end{proof}

\subsection{Assume \texorpdfstring{$ |S| = 3 $, $ G' = \mathcal{C}_p\times\mathcal{C}_q $ and $ C_{G'}(\mathcal{C}_2) = \mathcal{C}_p\times\mathcal{C}_q $}{Lg}}\hfill\label{3.5}

 In this subsection we prove the part of Theorem~\ref{theorem1.1} where, $ |S| = 3 $, $ G' = \mathcal{C}_p\times\mathcal{C}_q $, $  C_{G'}(\mathcal{C}_2) = \mathcal{C}_p\times\mathcal{C}_q $, and neither $ C_{G'}(\mathcal{C}_3) \neq \{e\} $ nor $ \widehat{S} $ is minimal holds. Recall $ \overline{G} = G/G' $, $ \widecheck{G} = G/\mathcal{C}_q $ and $ \widehat{G} = G/\mathcal{C}_p $.
   
   \begin{caseprop}\label{prop 3.5}
 Assume
 \begin{itemize}
     \item $ G = (\mathcal{C}_2\times\mathcal{C}_3)\ltimes(\mathcal{C}_p\times\mathcal{C}_q) $,
     \item $ |S| = 3 $,
     \item $ C_{G'}(\mathcal{C}_2) = \mathcal{C}_p\times\mathcal{C}_q $.
 \end{itemize}
 Then $ \Cay(G;S) $ contains a Hamiltonian cycle.
 \end{caseprop}
 \begin{proof}
        Let $ S = \{a,b,c\} $. If $ C_{G'}(\mathcal{C}_3) \neq \{e\} $, then Proposition~\ref{prop 3.3} applies. So we may assume $ C_{G'}(\mathcal{C}_3) = \{e\} $. Now if $ \widehat{S} $ is minimal, then Proposition~\ref{caseprop 3.4} applies. So we may assume $ \widehat{S} $ is not minimal. Consider
       \begin{align*}
           \widehat{G} = G/\mathcal{C}_p = (\mathcal{C}_2 \times\mathcal{C}_3) \ltimes \mathcal{C}_q \cong (\mathcal{C}_3\ltimes\mathcal{C}_q)\times\mathcal{C}_2.
       \end{align*}
       Choose a 2-element $ \{a,b\} $ subset of $ S $ that generates $ \widehat{G} $. From the minimality of $ S $, we see that
       \begin{align*}
           \langle a,b \rangle = (\mathcal{C}_3\ltimes\mathcal{C}_q)\times\mathcal{C}_2,
       \end{align*}
       after replacing $ a $ and $ b $ by conjugates. The projection of $ (a,b) $ to $ \mathcal{C}_3\ltimes\mathcal{C}_q $ must be of the form $ (\gt,\gq) $ or $ (\gt,\gt\gq) $ (perhaps after replacing $ a $ and/or $ b $ with its inverse; also note that $ \widehat{b} \neq \gq $ because $ S \cap G' = \emptyset $). Therefore, $ (a,b) $ must have one of the following forms:
       \begin{enumerate}
           \item $ (\gt,a_2\gq) $,
           \item $ (\gt,a_2\gt\gq) $,
           \item $ (a_2\gt,\gt\gq) $,
           \item $ (a_2\gt,a_2\gq) $,
           \item $ (a_2\gt,a_2\gt\gq) $.
       \end{enumerate}
       Let $ c $ be the third element of $ S $. We may write $ c = a_{2}^{i}\gt^{j}\gq^{k}\gamma_{p} $ with $ 0\leq i \leq 1 $, $ 0\leq j \leq 2 $ and $ 0\leq k \leq q-1 $. Note since $ S \cap G' = \emptyset $, we know that $ i $ and $ j $ cannot both be equal to $ 0 $. Additionally, we have $ \gt\gamma_{p} \gt^{-1} = \gamma_{p}^{\widehat{\tau}} $ where $ \widehat{\tau}^{3} \equiv 1  \pmod{p}$ and $ \widehat{\tau} \not \equiv 1 \pmod{p}$. Thus $ \widehat{\tau}^2+\widehat{\tau}+1 \equiv 0 \pmod{p} $. Note that this implies $ \widehat{\tau} \not \equiv -1 \pmod{p} $. Also we have $ \gt\gq\gt^{-1} = \gq^{\widecheck{\tau}} $. By using the same argument we can conclude that $ \widecheck{\tau} \not \equiv 1 \pmod{q} $ and $ \widecheck{\tau}^2+\widecheck{\tau}+1 \equiv 0 \pmod{q} $. Note that this implies $ \widecheck{\tau} \not \equiv -1 \pmod{q} $. Combining these facts with $ \widehat{\tau}^3\equiv 1 \pmod{p} $ and $ \widecheck{\tau}^3\equiv 1 \pmod{q} $, we conclude that $ \widehat{\tau}^2 \not\equiv \pm 1 \pmod{p} $, and $ \widecheck{\tau}^2 \not\equiv \pm1 \pmod{q} $.
       \setcounter{case}{0}
       \begin{case}
           Assume $ a = \gt $ and $ b = a_2\gq $.
           
           \begin{subcase} \label{Subcase 10.1.1}
               Assume $ i = 0 $. Then $ j \neq 0 $ and $ c = \gt^j\gq^k\gamma_p $. For future reference in Subcase~\ref{subcase 3.6.4.1} of Proposition~\ref{prop6.2}, we note that the argument here does not require our current assumption that $ \mathcal{C}_2 $ centralizes $ \mathcal{C}_p $. We may assume $ j = 1 $, by replacing $ c $ with $ c^{-1} $ if necessary. Then $ c = \gt\gq^k\gamma_p $. Consider $ \overline{G} = \mathcal{C}_2\times\mathcal{C}_3 $. Then we have $ \overline{a} = \overline{c} = \gt $, $ \overline{b} = a_2 $. We have $ C_1 = (\overline{c},\overline{a},\overline{b},\overline{a}^{-2},\overline{b}) $ and $ C_2 = (\overline{c}^2,\overline{b},\overline{a}^{-2},\overline{b}) $ as Hamiltonian cycles in $ \Cay(\overline{G};\overline{S}) $. Since there is one occurrence of $ c $ in $ C_1 $, then by Lemma~\ref{lemma 2.5.2} we conclude that the subgroup generated by $ \mathbb{V}(C_1) $ contains $ \mathcal{C}_p $. Also,
               \begin{align*}
                    \mathbb{V}(C_1) &= caba^{-2}b \\&\equiv \gt\gq^k\cdot\gt\cdot a_2\gq\cdot\gt^{-2}\cdot a_2\gq \pmod{\mathcal{C}_p} \\&= \gq^{k\widecheck{\tau}+\widecheck{\tau}^2+1} \\&= \gq^{\widecheck{\tau}^2+k\widecheck{\tau}+1}.
               \end{align*}
               We may assume this does not generate $ \mathcal{C}_q $, for otherwise Factor Group Lemma~\ref{FGL} applies. Therefore,
               \begin{align*}
                   0 \equiv \widecheck{\tau}^2+k\widecheck{\tau}+1 \pmod{q}. \tag{1.1A}\label{10.1.1A}
               \end{align*}
               We also have
               \begin{align*}
                   0 \equiv \widecheck{\tau}^2+\widecheck{\tau}+1 \pmod{q}. \tag{1.1B}\label{10.1.1B}
               \end{align*}
               By subtracting the above equation from \ref{10.1.1A}, we have $ 0 \equiv (k-1)\widecheck{\tau} \pmod{q} $. This implies that $ k = 1 $.
               
               Now we calculate the voltage of $ C_2 $.
               \begin{align*}
                   \mathbb{V}(C_2) &= c^2ba^{-2}b \\&\equiv \gt\gamma_p\gt\gamma_p\cdot a_2\cdot\gt^{-2}\cdot a_2 \pmod{\mathcal{C}_q} \\&= \gamma_p^{\widehat{\tau}+\widehat{\tau}^2} 
               \end{align*}
               which generates $ \mathcal{C}_p $. Also
               \begin{align*}
                   \mathbb{V}(C_2) &= c^2ba^{-2}b \\&\equiv \gt\gq\cdot\gt\gq\cdot a_2\gq\cdot\gt^{-2}\cdot a_2\gq \pmod{\mathcal{C}_p} \\&= \gq^{\widecheck{\tau}+\widecheck{\tau}^2+\widecheck{\tau}^2+1} \\&= \gq^{2\widecheck{\tau}^2+\widecheck{\tau}+1}.
               \end{align*}
               We may assume this does not generate $ \mathcal{C}_q $, for otherwise Factor Group Lemma~\ref{FGL} applies. Then 
               \begin{align*}
                  0 \equiv 2\widecheck{\tau}^2+\widecheck{\tau}+1 \pmod{q}. 
               \end{align*}
               By subtracting \ref{10.1.1B} from the above equation we have
               \begin{align*}
                   0 \equiv \widecheck{\tau}^2 \pmod{q}
               \end{align*}
               which is a contradiction.
           \end{subcase}
           \begin{subcase} \label{Subcase 10.1.2}
               Assume $ j = 0 $. Then $ i \neq 0 $ and $ c = a_2\gq^k\gamma_p $. For future reference in Subcase~\ref{subcase 3.6.4.2} of Proposition~\ref{prop6.2}, we note that the argument here does not require our current assumption that $ \mathcal{C}_2 $ centralizes $ \mathcal{C}_p $. If $ k \neq 0 $, then by Lemma~\ref{lemma 5.15}(\ref{5.15.3}) $ \langle a,c \rangle = G $ which contradicts the minimality of $ S $. 
               
               So we can assume $ k = 0 $. Then $ c = a_2\gamma_p $. Consider $ \overline{G} = \mathcal{C}_2\times\mathcal{C}_3 $. Then we have $ \overline{a} = \gt $ and $ \overline{b} = \overline{c} = a_2 $. This implies that $ |\overline{a}| = 3 $ and $ |\overline{b}| = |\overline{c}| = 2 $. We have $ C = (\overline{c}^{-1},\overline{a}^2,\overline{b},\overline{a}^{-2}) $ as a Hamiltonian cycle in $ \Cay(\overline{G};\overline{S}) $. Since there is one occurrence of $ c $ in $ C $, and it is the only generator of $ G $ that contains $ \gamma_p $, then by Lemma~\ref{lemma 2.5.2} we conclude that the subgroup generated by $ \mathbb{V}(C) $ contains $ \mathcal{C}_p $. Similarly, since there is one occurrence of $ b $ in $ C $, and it is the only generator of $ G $ that contains $ \gq $, then by Lemma~\ref{lemma 2.5.2} we conclude that the subgroup generated by $ \mathbb{V}(C) $ contains $ \mathcal{C}_q $. Therefore, the subgroup generated by $ \mathbb{V}(C) $ is $ G' $. So, Factor Group Lemma~\ref{FGL} applies.
           \end{subcase}
           \begin{subcase}
               Assume $ i \neq 0 $ and $ j \neq 0 $. Then $ c = a_2\gt^j\gq^k\gamma_p $. If $ k \neq 0 $, then by Lemma~\ref{lemma 5.15}(\ref{5.15.3}) $ \langle a,c \rangle = G $ which contradicts the minimality of $ S $. 
               
               So we can assume $ k = 0 $. We may also assume $ j = 1 $, by replacing $ c $ with $ c^{-1} $ if necessary. Then $ c = a_2\gt\gamma_p $. Consider $ \overline{G} = \mathcal{C}_2\times\mathcal{C}_3 $. Then we have $ \overline{a} = \gt $, $ \overline{b} = a_2 $ and $ \overline{c} = a_2\gt $. This implies that $ |\overline{a}| = 3 $, $ |\overline{b}| = 2 $ and $ |\overline{c}| = 6 $. We have $ C = (\overline{c},\overline{b},\overline{a},\overline{c},\overline{a}^{-1},\overline{c}) $ as a Hamiltonian cycle in $ \Cay(\overline{G};\overline{S}) $. Now we calculate its voltage.
               \begin{align*}
                   \mathbb{V}(C) &= cbaca^{-1}c \\&\equiv a_2\gt\cdot a_2\gq\cdot\gt\cdot a_2\gt\cdot\gt^{-1}\cdot a_2\gt \pmod{\mathcal{C}_p}\\&= \gt\gq\gt^2 \\&= \gq^{\widecheck{\tau}}
               \end{align*}
               which generates $ \mathcal{C}_q $. Also
               \begin{align*}
                   \mathbb{V}(C) &= cbaca^{-1}c \\&\equiv a_2\gt\gamma_p\cdot a_2\cdot\gt\cdot a_2\gt\gamma_p\cdot \gt^{-1}\cdot a_2\gt\gamma_p \pmod{\mathcal{C}_q}\\&= \gt\gamma_p\gt^2\gamma_p^2 \\&= \gamma_p^{\widehat{\tau}+2}.
               \end{align*}
               We may assume this does not generate $ \mathcal{C}_p $, for otherwise Factor Group Lemma~\ref{FGL} applies. Then $ \widehat{\tau} \equiv -2 \pmod{p} $. By substituting this in
               \begin{align*}
                   0 \equiv \widehat{\tau}^2+\widehat{\tau}+1 \pmod{p},
               \end{align*}
               we have 
               \begin{align*}
                   0 &\equiv 4-2+1 \pmod{p}\\&= 3.
               \end{align*}
               This contradicts the fact that $ p>3 $.
           \end{subcase}
       \end{case}
       \begin{case} \label{Case 10.2}
           Assume $ a = \gt $ and $ b = a_2\gt\gq $.
            \begin{subcase}
               Assume $ i \neq 0 $ and $ j \neq 0 $. Then $ c = a_2\gt^j\gq^k\gamma_p $. If $ k \neq 0 $, then by Lemma~\ref{lemma 5.15}(\ref{5.15.3}) $ \langle a,c \rangle = G $ which contradicts the minimality of $ S $. So we can assume $ k = 0 $. Then $ c = a_2\gt^j\gamma_p $. Thus, by Lemma~\ref{lemma 5.15}(\ref{lemma 5.15.5}) $ \langle b,c \rangle = G $ which contradicts the minimality of $ S $.
           \end{subcase}
           \begin{subcase}
               Assume $ i = 0 $. Then $ j \neq 0 $. We may assume $ j = 1 $, by replacing $ c $ with $ c^{-1} $ if necessary. Then $ c = \gt\gq^k\gamma_p $. 
               
               Suppose, for the moment, that $ k \neq 1 $. Then $ c = \gt\gq^k\gamma_p $. We have $ \langle \overline{b},\overline{c} \rangle = \langle \overline{a}_2\overline{a}_3,\overline{a}_3 \rangle = \overline{G} $. Consider $ \{\widehat{b},\widehat{c}\} = \{ a_2\gt\gq,\gt\gq^k\} $. Since $ \mathcal{C}_2 $ centralizes $ \mathcal{C}_q $, then
               \begin{align*}
                   [a_2\gt\gq,\gt\gq^k] = [\gt\gq,\gt\gq^k] = \gt\gq\gt\gq^k\gq^{-1}\gt^{-1}\gq^{-k}\gt^{-1} = \gq^{\widecheck{\tau}+k\widecheck{\tau}^2-\widecheck{\tau}^2-k\widecheck{\tau}} = \gq^{\widecheck{\tau}(k-1)(\widecheck{\tau}-1)}
               \end{align*}
               which generates $ \mathcal{C}_q $. Now consider $\{\widecheck{b},\widecheck{c}\} = \{a_2\gt,\gt\gamma_p\}$. Since $ \mathcal{C}_2 $ centralizes $ \mathcal{C}_p $, then
               \begin{align*}
                   [a_2\gt,\gt\gamma_p] = [\gt,\gt\gamma_p] = \gt\gt\gamma_p\gt^{-1}\gamma_p^{-1}\gt^{-1} = \gamma_p^{\widehat{\tau}^2-\widehat{\tau}} = \gamma_p^{\widehat{\tau}(\widehat{\tau}-1)}
               \end{align*}
               which generates $ \mathcal{C}_p $. Therefore, $ \langle b,c \rangle = G $ which contradicts the minimality of $ S $.
               
               Now we can assume $ k = 1 $. Then $ c = \gt\gq\gamma_p $. Consider $ \overline{G} = \mathcal{C}_2\times\mathcal{C}_3 $. We have $ \overline{a} = \overline{c} = \gt $ and $ \overline{b} = a_2\gt $. This implies that $ |\overline{a}| = |\overline{c}| = 3 $ and $ |\overline{b}| = 6 $. We have $ C = (\overline{c},\overline{b},\overline{a}^{2},\overline{b},\overline{a}) $ as a Hamiltonian cycle in $ \Cay(\overline{G};\overline{S}) $. Since there is one occurrence of $ c $ in $ C $, and it is the only generator of $ G $ that contains $ \gamma_p $, then by Lemma~\ref{lemma 2.5.2} we conclude that the subgroup generated by $ \mathbb{V}(C) $ is $ \mathcal{C}_p $. Also,
               \begin{align*}
                   \mathbb{V}(C) &= cba^2ba \\&\equiv \gt\gq\cdot a_2\gt\gq\cdot\gt^2\cdot a_2\gt\gq\cdot\gt \pmod{\mathcal{C}_p}\\&= \gt\gq\gt\gq^2\gt\\&= \gq^{\widecheck{\tau}+2\widecheck{\tau}^2}\\&= \gq^{\widecheck{\tau}(1+2\widecheck{\tau})}.
               \end{align*}
               We may assume this does not generate $ \mathcal{C}_q $, for otherwise Factor Group Lemma~\ref{FGL} applies. Therefore, $ 1+2\widecheck{\tau} \equiv 0 \pmod{q} $. This implies that $ \widecheck{\tau} \equiv -1/2 \pmod{q} $. By substituting $ \widecheck{\tau} \equiv -1/2 \pmod{q} $ in 
               \begin{align*}
                   \widecheck{\tau}^2+\widecheck{\tau}+1 \equiv 0 \pmod{q},
               \end{align*}
               then we have $ 3/4 \equiv 0 \pmod{q} $, which contradicts Assumption~\ref{assumption 3.1}(\ref{assumption 3.1.1}).
           \end{subcase}
           \begin{subcase}
               Assume $ j = 0 $. Then $ i \neq 0 $ and $ c = a_2\gq^k\gamma_p $. If $ k \neq 0 $, then by Lemma~\ref{lemma 5.15}(\ref{5.15.3}) $ \langle a,c \rangle = G $ which contradicts the minimality of $ S $.
               
               So we can assume $ k = 0 $. Then $ c = a_2\gamma_p $. Consider $ \overline{G} = \mathcal{C}_2\times\mathcal{C}_3 $. Then we have $ \overline{a} = \gt $, $ \overline{b} = a_2\gt $ and $ \overline{c} = a_2 $. This implies that $ |\overline{a}| = 3 $, $ |\overline{b}| = 6 $ and $ |\overline{c}| = 2 $. We have $ C = (\overline{c},\overline{a},\overline{b},\overline{a}^{-1},\overline{b}^2) $ as a Hamiltonian cycle in $ \Cay(\overline{G};\overline{S}) $. Since there is one occurrence of $ c $ in $ C $, and it is the only generator of $ G $ that contains $ \gamma_p $, then by Lemma~\ref{lemma 2.5.2} we conclude that the subgroup generated by $ \mathbb{V}(C) $ contains $ \mathcal{C}_p $. Also,
           \begin{align*}
               \mathbb{V}(C) &= caba^{-1}b^2 \\&\equiv a_2\cdot\gt\cdot a_2\gt\gq\cdot\gt^{-1}\cdot a_2\gt\gq a_2\gt\gq \pmod{\mathcal{C}_p} \\&= \gt^2\gq^2\gt\gq\\&= \gq^{2\widecheck{\tau}^2+1}.
           \end{align*}
           We may assume this does not generate $ \mathcal{C}_q $, for otherwise Factor Group Lemma~\ref{FGL} applies. Thus, $ \widecheck{\tau}^2 \equiv -1/2 \pmod{q} $. By substituting this in
           \begin{align*}
               \widecheck{\tau}^2+\widecheck{\tau}+1 \equiv 0 \pmod{q},
           \end{align*}
           we have $ \widecheck{\tau} \equiv -1/2 \pmod{q} $ which contradicts $ \widecheck{\tau}^2 \equiv -1/2 \pmod{q} $.
            \end{subcase}
       \end{case}
       \begin{case}
           Assume $ a = a_2\gt $ and $ b = \gt\gq $. Since $ b = \gt\gq $ is conjugate to $ \gt $ via an element of $ \mathcal{C}_q $, then $ \{a,b\} $ is conjugate to $ \{a_2\gt\gq^m,\gt\} $ for some nonzero $ m $. So Case~\ref{Case 10.2} applies (after replacing $ \gq $ with $ \gq^m $).
       \end{case}
       \begin{case} \label{Case 3.5.4}
           Assume $ a = a_2\gt $ and $ b = a_2\gq $.
           \begin{subcase}
               Assume $ i = 0 $. Then $ j \neq 0 $ and $ c = \gt^j\gq^k\gamma_p $. If $ k \neq 0 $, then by Lemma~\ref{lemma 5.15}(\ref{5.15.1}) $ \langle a,c \rangle = G $ which contradicts the minimality of $ S $.
               
               So we can assume $ k = 0 $. We may also assume $ j = 1 $, by replacing $ c $ with $ c^{-1} $ if necessary. Then $ c = \gt \gamma_{p} $. Consider $ \overline{G} = \mathcal{C}_2\times\mathcal{C}_3 $. Thus, $ \overline{a} = a_{2}\gt $, $ \overline{b} = a_{2} $ and $ \overline{c} = \gt $. This implies that $ |\overline{a}| = 6 $, $ |\overline{b}| = 2 $ and $ |\overline{c}| = 3 $. We have $ C = (\overline{a}^2,\overline{b},\overline{c},\overline{a},\overline{c}^{-1}) $ as a Hamiltonian cycle in $ \Cay(\overline{G};\overline{S}) $. Since there is one occurrence of $ b $ in $ C $, and it is the only generator of $ G $ that contains $ \gq $, then by Lemma~\ref{lemma 2.5.2} we conclude that the subgroup generated by $ \mathbb{V}(C) $ contains $ \mathcal{C}_q $. Also,
            \begin{align*}
                \mathbb{V}(C) &= a^2bcac^{-1} \\&\equiv  \gt^2\cdot a_{2}\cdot\gt\gamma_{p}\cdot a_{2}\gt\cdot\gamma_p^{-1}\gt^{-1} \pmod{\mathcal{C}_q} \\ &= \gamma_p\gt\gamma_p^{-1}\gt^{-1} \\ &= \gamma_p^{1-\widehat{\tau}}
            \end{align*}
            which generates $ \mathcal{C}_p $. Therefore, the subgroup generated by $ \mathbb{V}(C) $ is $ G' $. So, Factor Group Lemma~\ref{FGL} applies.
          \end{subcase}
          \begin{subcase}
              Assume $ j = 0 $. Then $ i \neq 0 $ and $ c = a_2\gq^k\gamma_p $. If $ k \neq 0 $, then by Lemma~\ref{lemma 5.15}(\ref{5.15.1}) $ \langle a,c \rangle = G $ which contradicts the minimality of $ S $.
              
              So we can assume $ k = 0 $. Then $ c = a_{2}\gamma_{p} $. Consider $ \overline{G} = \mathcal{C}_2\times\mathcal{C}_3 $, then $ \overline{a} = a_{2}\gt $ and $ \overline{b} = \overline{c} = a_{2} $. This implies that $ |\overline{a}| = 6 $ and $ |\overline{b}| = |\overline{c}| = 2 $. We have $ C = ((\overline{a},\overline{b})^2,\overline{a},\overline{c}) $ as a Hamiltonian cycle in $ \Cay(\overline{G};\overline{S}) $. Since there is one occurrence of $ c $ in $ C $, and it is the only generator of $ G $ that contains $ \gamma_p $, then by Lemma~\ref{lemma 2.5.2} we conclude that the subgroup generated by $ \mathbb{V}(C) $ contains $ \mathcal{C}_p $. Also, 
                \begin{align*}
                    \mathbb{V}(C) &= (ab)^2ac \\&\equiv (a_{2}\gt\cdot a_{2}\gq)^2\cdot a_{2}\gt\cdot a_{2} \pmod{\mathcal{C}_p} \\ &= \gt\gq\gt\gq\gt  \\ &= \gq^{\widecheck{\tau}+\widecheck{\tau}^2}.
                \end{align*}
                which generates $ \mathcal{C}_q $. Therefore, the subgroup generated by $ \mathbb{V}(C) $ is $ G' $. Thus, Factor Group Lemma~\ref{FGL} applies.
          \end{subcase}
          \begin{subcase}
              Assume $ i \neq 0 $ and $ j \neq 0 $. Then $ c = a_2\gt^j\gq^k\gamma_p $. If $ k \neq 0 $, then by Lemma~\ref{lemma 5.15}(\ref{5.15.1}) $ \langle a,c \rangle = G $ which contradicts the minimality of $ S $.
              
               So we can assume $ k = 0 $. We may also assume $ j = 1 $, by replacing $ c $ with $ c^{-1} $ if necessary. Then $ c = a_{2}\gt\gamma_{p} $. Consider $ \overline{G} = \mathcal{C}_2\times\mathcal{C}_3 $. Thus, $ \overline{a} = \overline{c} = a_{2}\gt $ and $ \overline{b} = a_{2} $. This implies that $ |\overline{a}| = |\overline{c}| = 6 $ and $ |\overline{b}| = 2 $. We have $ C = (\overline{a},\overline{c},\overline{b},\overline{a}^{-2},\overline{b}) $ as a Hamiltonian cycle in $ \Cay(\overline{G};\overline{S}) $. Since there is one occurrence of $ c $ in $ C $, and it is the only generator of $ G $ that contains $ \gamma_p $, then by Lemma~\ref{lemma 2.5.2} we conclude that the subgroup generated by $ \mathbb{V}(C) $ contains $ \mathcal{C}_p $. Also,
                \begin{align*}
                    \mathbb{V}(C) &= acba^{-2}b \\&\equiv a_{2}\gt\cdot a_{2}\gt\cdot a_{2}\gq\cdot \gt^{-2}\cdot a_{2}\gq \pmod{\mathcal{C}_p} \\ &= \gt^2\gq\gt^{-2}\gq  \\ &= \gq^{\widecheck{\tau}^2+1}
                 \end{align*}
                 which generates $ \mathcal{C}_q $, because $ \widecheck{\tau}^2 \not\equiv -1 \pmod{q} $. Therefore, the subgroup generated by $ \mathbb{V}(C) $ is $ G' $. So, Factor Group Lemma~\ref{FGL} applies.
          \end{subcase}
       \end{case}
       \begin{case}
           Assume $ a = a_2\gt $ and $ b = a_2\gt\gq $. If $ k \neq 0 $, then by Lemma~\ref{lemma 5.15}(\ref{5.15.1}) $ \langle a,c \rangle = G $ which contradicts the minimality of $ S $. So we can assume $ k = 0 $. Also, if $ j \neq 0 $, then by Lemma~\ref{lemma 5.15}(\ref{lemma 5.15.5}) $ \langle b,c \rangle = G $ which contradicts the minimality of $ S $. So we may also assume $ j = 0 $. Then $ i \neq 0 $. Therefore, $ c = a_2\gamma_p $. So Case~\ref{Case 3.5.4} applies, after interchanging $ b $ and $ c $, and interchanging $ p $ and $ q $.
            \qedhere
       \end{case}
 \end{proof}
 
 \subsection{Assume \texorpdfstring{$ |S| = 3 $, $ G' = \mathcal{C}_p\times\mathcal{C}_q $ and $ C_{G'}(\mathcal{C}_2) \neq \{e\} $ }{Lg}}\hfill\label{3.6}

 In this subsection we prove the part of Theorem~\ref{theorem1.1} where, $ |S| = 3 $, $ G' = \mathcal{C}_p\times\mathcal{C}_q $, $ C_{G'}(\mathcal{C}_2) \neq \{e\} $, and neither $ C_{G'}(\mathcal{C}_2) = \mathcal{C}_p\times\mathcal{C}_q $ nor $ C_{G'}(\mathcal{C}_3) \neq \{e\} $ nor $ \widehat{S} $ is minimal holds. Recall $ \overline{G} = G/G' $, $ \widecheck{G} = G/\mathcal{C}_q $ and $ \widehat{G} = G/\mathcal{C}_p $.
\begin{caseprop} \label{prop6.2}
 Assume
 \begin{itemize}
     \item $ G = (\mathcal{C}_2\times\mathcal{C}_3)\ltimes(\mathcal{C}_p\times\mathcal{C}_q) $,
     \item $ |S| = 3 $,
     \item $ C_{G'}(\mathcal{C}_2) \neq \{e\} $.
 \end{itemize}
 Then $ \Cay(G;S) $ contains a Hamiltonian cycle.
 \end{caseprop}
 
 \begin{proof}
Let $ S = \{a,b,c\} $. If $ C_{G'}(\mathcal{C}_3) \neq \{e\} $, then Proposition~\ref{prop 3.3} applies. Therefore, we may assume $ C_{G'}(\mathcal{C}_3) = \{e\}  $. Now if $ C_{G'}(\mathcal{C}_2) = \mathcal{C}_p\times\mathcal{C}_q $, then Proposition~\ref{prop 3.5} applies. Since $ C_{G'}(\mathcal{C}_2) \neq \{e\} $, then we may assume $ C_{G'}(\mathcal{C}_2) = \mathcal{C}_q $, by interchanging $ q $ and $ p $ if necessary. This implies that $ \mathcal{C}_2 $ inverts $ \mathcal{C}_p $. Now if $ \widehat{S} $ is minimal, then Proposition~\ref{caseprop 3.4} applies. So we may assume $ \widehat{S} $ is not minimal. Consider
           \begin{align*}
              \widehat{G} = G/\mathcal{C}_p =  (\mathcal{C}_2\times\mathcal{C}_3) \ltimes \mathcal{C}_q.
           \end{align*}
            Choose a 2-element subset $ \{a,b\} $ in $ S $ that generates $ \widehat{G} $. From the minimality of $ S $, we see that
            \begin{align*}
                \langle a,b \rangle = (\mathcal{C}_2\times\mathcal{C}_3) \ltimes \mathcal{C}_q
            \end{align*}
            after replacing $ a $ and $ b $ by conjugates.
            We may assume $ |\overline{a}| \geq |\overline{b}| $ and (by conjugating if necessary) $ a $ is an element of $ \mathcal{C}_2\times\mathcal{C}_3 $. Then the projection of $ (a,b) $ to $ \mathcal{C}_2\times\mathcal{C}_3 $ has one of the following forms after replacing $ a $ and $ b $ with their inverses if necessary.
          \begin{itemize}
              \item $ (a_{2}\gt,a_2\gt) $,
              \item $ (a_{2}\gt,a_{2}) $,
              \item $ (a_{2}\gt,\gt) $,
              \item $ (\gt,a_2) $.
          \end{itemize}
           So there are four possibilities for $ (a,b) $:
          \begin{enumerate}
              \item $ (a_{2}\gt,a_{2}\gt\gq ) $,
              \item $ (a_{2}\gt,a_{2}\gq) $,
              \item $ (a_{2}\gt,\gt\gq) $,
              \item $ (\gt,a_{2}\gq) $.
          \end{enumerate}
        Let $ c $ be the third element of $ S $. We may write $ c = a_{2}^{i}\gt^{j}\gq^{k}\gamma_{p} $ with $ 0\leq i \leq 1 $, $ 0\leq j \leq 2 $ and $ 0\leq k \leq q-1 $. Note since $ S \cap G' = \emptyset $, we know that $ i $ and $ j $ cannot both be equal to $ 0 $. Additionally, we have $ \gt\gamma_{p} \gt^{-1} = \gamma_{p}^{\widehat{\tau}} $ where $ \widehat{\tau}^{3} \equiv 1  \pmod{p}$ and $ \widehat{\tau} \not \equiv 1 \pmod{p}$. Thus $ \widehat{\tau}^2+\widehat{\tau}+1 \equiv 0 \pmod{p} $. Note that this implies $ \widehat{\tau} \not \equiv -1 \pmod{p} $. Also we have $ \gt\gq\gt^{-1} = \gq^{\widecheck{\tau}} $. By using the same argument we can conclude that $ \widecheck{\tau} \not \equiv 1 \pmod{q} $ and $ \widecheck{\tau}^2+\widecheck{\tau}+1 \equiv 0 \pmod{q} $. Note that this implies $ \widecheck{\tau} \not \equiv -1 \pmod{q} $. Therefore, we conclude that $ \widehat{\tau}^2 \not\equiv \pm 1 \pmod{p} $, and $ \widecheck{\tau}^2 \not\equiv \pm1 \pmod{q} $.
        \setcounter{case}{0}
        \begin{case}
        Assume $ a = a_{2}\gt $ and $ b = a_{2}\gt\gq $. If $ k \neq 0 $, then by Lemma~\ref{lemma 5.15}(\ref{5.15.1}), $ \langle a,c\rangle = G $ which contradicts the minimality of $ S $. So we can assume $ k = 0 $. Now if $ j \neq 0 $, then by Lemma~\ref{lemma 5.15}(\ref{lemma 5.15.5}), $ \langle b,c\rangle = G $ which contradicts the minimality of $ S $. 
        Therefore, we may assume $ j = 0 $. Then $ i \neq 0 $ and $ c = a_2\gamma_p $. Consider $ \overline{G} = \mathcal{C}_2\times\mathcal{C}_3 $. Thus $ \overline{a} = \overline{b} = a_{2}\gt $ and $ \overline{c} = a_{2} $. Therefore, 
            $ |\overline{a}| = |\overline{b}| = 6 $ and $ |\overline{c}| = 2 $. We have $ C = (\overline{a},\overline{b},\overline{c},\overline{a}^{-2},\overline{c}) $ as a Hamiltonian cycle in $ \Cay(\overline{G};\overline{S}) $. Since there is one occurrence of $ b $ in $ C $, and it is the only generator of $ G $ that contains $ \gq $, then by Lemma~\ref{lemma 2.5.2} we conclude that the subgroup generated by $ \mathbb{V}(C) $ contains $ \mathcal{C}_q $. Also,
            \begin{align*}
                \mathbb{V}(C) &= abca^{-2}c \\&\equiv a_{2}\gt\cdot a_{2}\gt\cdot a_{2}\gamma_{p}\cdot\gt^{-2}\cdot a_{2}\gamma_{p} \pmod{\mathcal{C}_q} \\ &= \gt^2\gamma_{p}^{-1}\gt^{-2}\gamma_{p} \\ &= \gamma_{p}^{-\widehat{\tau}^2+1}
            \end{align*}
            which generates $ \mathcal{C}_p $. Therefore, the subgroup generated by $ \mathbb{V}(C) $ is $ G' $. So, Factor Group Lemma~\ref{FGL} applies.
        \end{case}
        \begin{case}
         Assume $ a = a_{2}\gt $ and $ b = a_{2}\gq $. 
        \begin{subcase}
            Assume $ i = 0 $. Then $ j \neq 0 $ and $ c = \gt^j\gq^k\gamma_{p} $. If $ k \neq 0 $, then by Lemma~\ref{lemma 5.15}(\ref{5.15.1}), $ \langle a,c\rangle = G $ which contradicts the minimality of $ S $.
            
            So we can assume $ k = 0 $. We may also assume $ j = 1 $, by replacing $ c $ with $ c^{-1} $ if necessary. Then $ c = \gt \gamma_{p} $. Consider $ \overline{G} = \mathcal{C}_2\times\mathcal{C}_3 $. Thus, $ \overline{a} = a_{2}\gt $, $ \overline{b} = a_{2} $ and $ \overline{c} = \gt $. Therefore, $ |\overline{a}| = 6 $, $ |\overline{b}| = 2 $ and $ |\overline{c}| = 3 $. We have $ C = (\overline{a}^2,\overline{b},\overline{c},\overline{a},\overline{c}^{-1}) $ as a Hamiltonian cycle in $ \Cay(\overline{G};\overline{S}) $. Since there is one occurrence of $ b $ in $ C $, and it is the only generator of $ G $ that contains $ \gq $, then by Lemma~\ref{lemma 2.5.2} we conclude that the subgroup generated by $ \mathbb{V}(C) $ contains $ \mathcal{C}_q $. Also,
            \begin{align*}
                \mathbb{V}(C) &= a^2bc^{-1}ac \\&\equiv \gt^2\cdot a_{2}\cdot\gt\gamma_{p}\cdot a_{2}\gt\cdot\gamma_p^{-1}\gt^{-1} \pmod{\mathcal{C}_q} \\ &= \gamma_p^{-1}\gt\gamma_p^{-1}\gt^{-1} \\ &= \gamma_{p}^{-1-\widehat{\tau}}
            \end{align*}
           which generates $ \mathcal{C}_p $. Therefore, the subgroup generated by $ \mathbb{V}(C) $ is $ G' $. So, Factor Group Lemma~\ref{FGL} applies.
        \end{subcase}
        \begin{subcase}
                Assume $ j = 0 $. Then $ i \neq 0 $ and $ c = a_{2}\gq^k\gamma_{p} $. If $ k \neq 0 $, then by Lemma~\ref{lemma 5.15}(\ref{5.15.1}), $ \langle a,c\rangle = G $ which contradicts the minimality of $ S $. 
                
                So we can assume $ k = 0 $. Then $ c = a_{2}\gamma_{p} $. Consider $ \overline{G} = \mathcal{C}_2\times\mathcal{C}_3 $, then $ \overline{a} = a_{2}\gt $ and $ \overline{b} = \overline{c} = a_{2} $. We have $ C = ((\overline{a},\overline{b})^2,\overline{a},\overline{c}) $ as a Hamiltonian cycle in $ \Cay(\overline{G};\overline{S}) $. Since there is one occurrence of $ c $ in $ C $, and it is the only generator of $ G $ that contains $ \gamma_p $, then by Lemma~\ref{lemma 2.5.2} we conclude that the subgroup generated by $ \mathbb{V}(C) $ contains $ \mathcal{C}_p $. Now we calculate its voltage. Also,
                \begin{align*}
                    \mathbb{V}(C) &= (ab)^2ac \\&\equiv (a_{2}\gt\cdot a_{2}\gq)^2\cdot a_{2}\gt\cdot a_{2} \pmod{\mathcal{C}_p} \\ &= \gt\gq\gt\gq\gt  \\ &= \gq^{\widecheck{\tau}+\widecheck{\tau}^2}.
                \end{align*}
                which generates $ \mathcal{C}_q $. Therefore, the subgroup generated by $ \mathbb{V}(C) $ generates $ G' $. So, Factor Group Lemma~\ref{FGL} applies.
        \end{subcase}
        \begin{subcase}
                Assume $ i \neq 0 $ and $ j \neq 0 $. If $ k \neq 0 $, then $ c = a_{2}\gt^j\gq^k\gamma_{p} $. Thus, by Lemma~\ref{lemma 5.15}(\ref{5.15.1}), $ \langle a,c\rangle = G $ which contradicts the minimality of $ S $. 
                
                So we can assume $ k = 0 $. We may also assume $ j = 1 $, by replacing $ c $ with $ c^{-1} $ if necessary. Then $ c = a_{2}\gt\gamma_{p} $. Consider $ \overline{G} = \mathcal{C}_2\times\mathcal{C}_3 $. Thus, $ \overline{a} = \overline{c} = a_{2}\gt $ and $ \overline{b} = a_{2} $. Therefore, $ |\overline{a}| = |\overline{c}| = 6 $ and $ |\overline{b}| = 2 $. We have $ C = (\overline{a},\overline{c},\overline{b},\overline{a}^{-2},\overline{b}) $ as a Hamiltonian cycle in $ \Cay(\overline{G};\overline{S}) $. Since there is one occurrence of $ c $ in $ C $, and it is the only generator of $ G $ that contains $ \gamma_p $, then by Lemma~\ref{lemma 2.5.2} we conclude that the subgroup generated by $ \mathbb{V}(C) $ contains $ \mathcal{C}_p $. Also,
                \begin{align*}
                    \mathbb{V}(C) &= acba^{-2}b \\&\equiv a_{2}\gt\cdot a_{2}\gt\cdot a_{2}\gq\cdot \gt^{-2}\cdot a_{2}\gq \pmod{\mathcal{C}_p} \\ &= \gt^2\gq\gt^{-2}\gq  \\ &= \gq^{\widecheck{\tau}^2+1}.
                 \end{align*}
                 Since $ \widecheck{\tau}^2 \not\equiv -1 \pmod{q} $, Factor Group Lemma~\ref{FGL} applies.
        \end{subcase}
        \end{case}
        \begin{case}
        Assume $ a = a_{2}\gt $ and $ b = \gt\gq $.
        \begin{subcase}
             Assume $ i \neq 0 $ and $ j \neq 0 $. If $ k = 0 $, then $ c = a_{2}\gt^j\gamma_{p} $. Thus, by Lemma~\ref{lemma 5.15}(\ref{5.15.2}), $ \langle b,c\rangle = G $ which contradicts the minimality of $ S $. So we can assume $ k \neq 0 $. Then $ c = a_{2}\gt^j\gq^k\gamma_{p} $. Thus, by Lemma~\ref{lemma 5.15}(\ref{5.15.1}), $ \langle a,c\rangle = G $ which contradicts the minimality of $ S $.
        \end{subcase}
        \begin{subcase}
               Assume $ i = 0 $. Then $ j \neq 0 $ and $ c = \gt^j\gq^k\gamma_{p} $. If $ k \neq 0 $, then by Lemma~\ref{lemma 5.15}(\ref{5.15.1}), $ \langle a,c\rangle = G $ which contradicts the minimality of $ S $.
               
               So we can assume $ k = 0 $. We may also assume $ j = 1 $, by replacing $ c $ with $ c^{-1} $ if necessary. Then $ c = \gt\gamma_{p} $. Consider $ \overline{G} = \mathcal{C}_2\times\mathcal{C}_3 $, then $ \overline{a} = a_{2}\gt $, $ \overline{b} = \overline{c} = \gt $. Therefore, $ |\overline{a}| = 6 $ and $ |\overline{b}| = |\overline{c}| = 3 $. We have $ C = (\overline{c},\overline{b},\overline{a},\overline{b}^{-2},\overline{a}^{-1}) $ as a Hamiltonian cycle in $ \Cay(\overline{G};\overline{S}) $. Since there is one occurrence of $ c $ in $ C $, and it is the only generator of $ G $ that contains $ \gamma_p $, then by Lemma~\ref{lemma 2.5.2} we conclude that the subgroup generated by $ \mathbb{V}(C) $ contains $ \mathcal{C}_p $. Also,
             \begin{align*}
               \mathbb{V}(C) &= cbab^{-2}a^{-1} \\&\equiv \gt\cdot\gt\gq\cdot a_{2}\gt\cdot\gq^{-1}\gt^{-1}\gq^{-1}\gt^{-1}\cdot\gt^{-1}a_{2} \pmod{\mathcal{C}_p} \\ &= \gt^2\gq\gt\gq^{-1}\gt^{-1}\gq^{-1}\gt^{-2}  \\ &= \gq^{\widecheck{\tau}^2-1-\widecheck{\tau}^{-1}} \\&= \gq^{\widecheck{\tau}^2-1-\widecheck{\tau}^2} \\&= \gq^{-1}
             \end{align*}
             which generates $ \mathcal{C}_q $. Therefore, the subgroup generated by $ \mathbb{V}(C) $ is $ G' $. So, Factor Group Lemma~\ref{FGL} applies.
        \end{subcase}
        \begin{subcase}
             Assume $ j = 0 $. Then $ i \neq 0 $ and $ c = a_{2}\gq^k\gamma_{p} $. If $ k \neq 0 $, then by Lemma~\ref{lemma 5.15}(\ref{5.15.1}), $ \langle a,c\rangle = G $ which contradicts the minimality of $ S $. 
             
             So we can assume $ k = 0 $. Then $ c = a_{2}\gamma_{p} $. Consider $ \overline{G} = \mathcal{C}_2\times\mathcal{C}_3 $, then $ \overline{a} = a_{2}\gt $, $ \overline{b} = \gt $ and $ \overline{c} = a_{2} $. Therefore, $ |\overline{a}| = 6 $, $ |\overline{b}| = 3 $ and $ |\overline{c}| = 2 $. We have $  C = (\overline{a},\overline{c},\overline{b},\overline{a},\overline{b}^{-1},\overline{a}) $ as a Hamiltonian cycle in $ \Cay(\overline{G};\overline{S}) $. Since there is one occurrence of $ c $ in $ C $, and it is the only generator of $ G $ that contains $ \gamma_p $, then by Lemma~\ref{lemma 2.5.2} we conclude that the subgroup generated by $ \mathbb{V}(C) $ contains $ \mathcal{C}_p $.
             Also,
             \begin{align*}
                 \mathbb{V}(C) &= acbab^{-1}a \\&\equiv a_{2}\gt\cdot a_{2}\cdot\gt\gq\cdot a_{2}\gt\cdot\gq^{-1}\gt^{-1}\cdot a_{2}\gt \pmod{\mathcal{C}_p} \\ &= \gt^2\gq\gt\gq^{-1} \\ &= \gq^{\widecheck{\tau}^2-1} .
             \end{align*}
             Since $ \widecheck{\tau}^2 \not\equiv 1 \pmod{q} $, Factor Group Lemma~\ref{FGL} applies.
        \end{subcase}
        \end{case}
        \begin{case}
        Assume $ a = \gt $ and $ b = a_{2}\gq $.
         \begin{subcase}\label{subcase 3.6.4.1}
               Assume $ i = 0 $. Then $ j \neq 0 $ and $ c = \gt^j\gq^k\gamma_p $. Thus, the argument in Subcase~\ref{Subcase 10.1.1} of Proposition~\ref{prop 3.5} applies.
           \end{subcase}
           \begin{subcase}\label{subcase 3.6.4.2}
               Assume $ j = 0 $. Then $ i \neq 0 $ and $ c = a_2\gq^k\gamma_p $. Thus, the argument in Subcase~\ref{Subcase 10.1.2} of Proposition~\ref{prop 3.5} applies.
           \end{subcase}
               \begin{subcase}
               Assume $ i \neq 0 $ and $ j \neq 0 $. Then $ c = a_2\gt^j\gq^k\gamma_p $. If $ k \neq 0 $, then by Lemma~\ref{lemma 5.15}(\ref{5.15.3}) $ \langle a,c \rangle = G $ which contradicts the minimality of $ S $. 
               
               So we can assume $ k = 0 $. We may also assume $ j = 1 $, by replacing $ c $ with $ c^{-1} $ if necessary. Then $ c = a_2\gt\gamma_p $. Consider $ \overline{G} = \mathcal{C}_2\times\mathcal{C}_3 $. Then we have $ \overline{a} = \gt $, $ \overline{b} = a_2 $ and $ \overline{c} = a_2\gt $. This implies that $ |\overline{a}| = 3 $, $ |\overline{b}| = 2 $ and $ |\overline{c}| = 6 $. We have $ C = (\overline{c},\overline{b},\overline{a},\overline{c},\overline{a}^{-1},\overline{c}) $ as a Hamiltonian cycle in $ \Cay(\overline{G};\overline{S}) $. Since there is one occurrence of $ b $ in $ C $, and it is the only generator of $ G $ that contains $ \gq $, then by Lemma~\ref{lemma 2.5.2} we conclude that the subgroup generated by $ \mathbb{V}(C) $ contains $ \mathcal{C}_q $. Also, since $ a_2 $ inverts $ \mathcal{C}_p $
               \begin{align*}
                   \mathbb{V}(C) &= cbaca^{-1}c \\&\equiv a_2\gt\gamma_p\cdot a_2\cdot\gt\cdot a_2\gt\gamma_p\cdot \gt^{-1}\cdot a_2\gt\gamma_p \pmod{\mathcal{C}_q}\\&= \gt\gamma_p^{-1}\gt^2 \\&= \gamma_p^{-\widehat{\tau}}
               \end{align*}
               which generates $ \mathcal{C}_p $. Therefore, the subgroup generated by $ \mathbb{V}(C) $ is $ G' $. So, Factor Group Lemma~\ref{FGL} applies.
               \qedhere
           \end{subcase}
        \end{case}
        \end{proof}
        
 \subsection{Assume \texorpdfstring{$ |S| = 3 $,   $ G' = \mathcal{C}_p \times\mathcal{C}_q $ and $ C_{G'}(\mathcal{C}_2) = \{e\} $}{Lg}}\hfill\label{3.7}

  In this subsection we prove the part of Theorem~\ref{theorem1.1} where, $ |S| = 3 $, $ G' = \mathcal{C}_p \times\mathcal{C}_q $, $ C_{G'}(\mathcal{C}_2) = \{e\} $, and neither $ C_{G'}(\mathcal{C}_3) \neq \{e\} $ nor $ \widehat{S} $ is minimal holds. Recall $ \overline{G} = G/G' $, $ \widecheck{G} = G/\mathcal{C}_q $ and $ \widehat{G} = G/\mathcal{C}_p $.
        \begin{caseprop}
             Assume
             \begin{itemize}
                 \item $ G = (\mathcal{C}_2\times\mathcal{C}_3)\ltimes(\mathcal{C}_p\times\mathcal{C}_q) $,
                 \item $ |S| = 3 $,
                 \item $ C_{G'}(\mathcal{C}_2) = \{e\} $.
             \end{itemize}
             Then $ \Cay(G;S) $ contains a Hamiltonian cycle.
             \end{caseprop}
             \begin{proof}
              Let $ S = \{a,b,c\} $. If $ C_{G'}(\mathcal{C}_3) \neq \{e\} $, then Proposition~\ref{prop 3.3} applies. So we may assume $ C_{G'}(\mathcal{C}_3) = \{e\} $. Now if $ \widehat{S} $ is minimal, then Proposition~\ref{caseprop 3.4} applies. So we may assume $ \widehat{S} $ is not minimal. Consider
            \begin{align*}
               \widehat{G} = G/\mathcal{C}_p =  (\mathcal{C}_2\times\mathcal{C}_3) \ltimes \mathcal{C}_q.
            \end{align*}
             Choose a 2-element subset $ \{a,b\} $ in $ S $ that generates $ \widehat{G} $. From the minimality of $ S $, we see
             \begin{align*}
                 \langle a,b \rangle = (\mathcal{C}_2\times\mathcal{C}_3) \ltimes \mathcal{C}_q.
             \end{align*}
             after replacing $ a $ and $ b $ by conjugates. We may assume $ |a| \geq |b| $ and (by conjugating if necessary) $ a $ is in $ \mathcal{C}_2\times\mathcal{C}_3 $. Then the projection of $ (a,b) $ to $ \mathcal{C}_2\times\mathcal{C}_3 $ is one of the following forms after replacing $ a $ and $ b $ with their inverses if necessary.
          \begin{itemize}
              \item $ (a_{2}\gt,a_2\gt) $,
              \item $ (a_{2}\gt,a_{2}) $,
              \item $ (a_{2}\gt,\gt) $,
              \item $ (\gt,a_2) $.
          \end{itemize}
           There are four possibilities for $ (a,b) $:
            \begin{enumerate}
                \item $ (a_{2}\gt,a_{2}\gt\gq) $,
                \item $ (a_{2}\gt,a_{2}\gq) $,
                \item $ (a_{2}\gt,\gt\gq) $,
                \item $ (\gt,a_{2}\gq) $.
            \end{enumerate}
          Let $ c $ be the third element of $ S $. We may write $ c = a_{2}^{i}\gt^{j}\gq^{k}\gamma_{p} $ with $ 0\leq i \leq 1 $, $ 0\leq j \leq 2 $ and $ 0\leq k \leq q-1 $. Note since $ S \cap G' = \emptyset $, we know that $ i $ and $ j $ cannot both be equal to $ 0 $. Additionally, we have $ \gt\gamma_{p} \gt^{-1} = \gamma_{p}^{\widehat{\tau}} $ where $ \widehat{\tau}^{3} \equiv 1  \pmod{p}$ and $ \widehat{\tau} \not \equiv 1 \pmod{p}$. Thus $ \widehat{\tau}^2+\widehat{\tau}+1 \equiv 0 \pmod{p} $. Note that this implies $ \widehat{\tau} \not \equiv -1 \pmod{p} $. We have $ \gt\gq\gt^{-1} = \gq^{\widecheck{\tau}} $. By using the same argument we can conclude that $ \widecheck{\tau} \not \equiv 1 \pmod{q} $ and  $ \widecheck{\tau}^2+\widecheck{\tau}+1 \equiv 0 \pmod{q} $. Note that this implies $ \widecheck{\tau} \not \equiv -1 \pmod{q} $. Therefore, we conclude that $ \widehat{\tau}^2 \not\equiv \pm 1 \pmod{p} $, and $ \widecheck{\tau}^2 \not\equiv \pm1 \pmod{q} $.
          \setcounter{case}{0}
        \begin{case}
         Assume $ a = a_{2}\gt $ and $ b = a_{2}\gt\gq $. If $ k \neq 0 $, then by Lemma~\ref{lemma 5.15}(\ref{5.15.1})~$ \langle a,c\rangle =~G $ which contradicts the minimality of $ S $. So we can assume $ k = 0 $. Now if $ j \neq 0 $, then by Lemma~\ref{lemma 5.15}(\ref{lemma 5.15.5}), $ \langle b,c\rangle = G $ which contradicts the minimality of $ S $. Therefore, we may assume $ j = 0 $. Then $ i \neq 0 $ and $ c = a_2\gamma_p $. We have $ \langle \overline{b},\overline{c} \rangle = \langle \overline{a}_2\overline{a}_3,\overline{a}_2 \rangle = \overline{G} $. Consider $ \{\widecheck{b},\widecheck{c}\} = \{a_{2}\gt,a_{2}\gamma_{p}\} $. Therefore,
               \begin{align*}
                   [a_{2}\gt,a_{2}\gamma_{p}] = a_2\gt a_2\gamma_p\gt^{-1}a_2\gamma_p^{-1}a_2 = \gt\gamma_{p}\gt^{-1}\gamma_{p} = \gamma_{p}^{\widehat{\tau}+1}.
               \end{align*}
                which generates $ \mathcal{C}_p $. Now consider $ \{\widehat{b},\widehat{c}\} = \{a_{2}\gt\gq,a_{2}\} $, then
               \begin{align*}
                   [a_{2}\gt\gq,a_{2}] = a_{2}\gt\gq a_{2}\gq^{-1}\gt^{-1}a_{2}a_{2} = 
                   \gt\gq^{-2}\gt^{-1} = \gq^{-2\widecheck{\tau}}
               \end{align*}
               which generates $  \mathcal{C}_q $. Therefore, $ \langle b,c \rangle = G $ which contradicts the minimality of $ S $.
        \end{case}
        \begin{case}
         Assume $ a = a_{2}\gt $ and $ b = a_{2}\gq $. If $ k \neq 0 $, then by Lemma~\ref{lemma 5.15}(\ref{5.15.1}), $ \langle a,c\rangle = G $ which contradicts the minimality of $ S $. So we can assume $ k = 0 $.
         \begin{subcase}
              \label{Subsubsubsubcase2.2.3.2.2}
    Assume $ j \neq 0 $. We may also assume $ j = 1 $, by replacing $ c $ with $ c^{-1} $ if necessary. Then $ c = a_2^i\gt\gamma_{p} $. We have $ \langle \overline{b},\overline{c} \rangle = \langle \overline{a}_2,\overline{a}_2^i\overline{a}_3 \rangle = \overline{G} $. Consider $ \{\widehat{b},\widehat{c}\} = \{a_2\gq,a_2^i\gt\} $. We have
    \begin{align*}
        [a_2\gq,a_2^i\gt] &= a_2\gq a_2^i\gt\gq^{-1}a_2\gt^{-1}a_2^i = \gq^{-1}a_2^{i+1}\gt\gq^{-1}\gt^{-1}a_2^{i+1} \\&= \gq^{-1}\gt\gq^{\mp1}\gt^{-1} = \gq^{-1\mp\widecheck{\tau}}
    \end{align*}
    which generates $ \mathcal{C}_q $. Now consider $ \{\widecheck{b},\widecheck{c}\} = \{a_2,a_2^i\gt\gamma_p\} $. We have
    \begin{align*}
        [a_2,a_2^i\gt\gamma_p] = a_2 a_2^i\gt\gamma_p a_2\gamma_p^{-1}\gt^{-1}a_2^i = a_2^{i+1}\gt\gamma_p^2\gt^{-1}a_2^{i+1} = \gamma_p^{\pm2\widehat{\tau}}
    \end{align*}
    which generates $ \mathcal{C}_p $. Therefore, $ \langle b,c \rangle = G $ which contradicts the minimality of $ S $.
         \end{subcase}
         \begin{subcase}
                 Assume $ j = 0 $. Then $ i \neq 0 $ and $ c = a_{2}\gamma_{p} $. Consider $ \overline{G} = \mathcal{C}_2\times\mathcal{C}_3 $, then $ \overline{a} = a_{2}\gt $ and $ \overline{b} = \overline{c} = a_{2} $. Thus, $ |\overline{a}| = 6 $ and $ |\overline{b}| = |\overline{c}| = 2 $. We have $ C =  ((\overline{a},\overline{b})^2,\overline{a},\overline{c}) $ as a Hamiltonian cycle in $ \Cay(\overline{G};\overline{S}) $. Since there is one occurrence of $ c $ in $ C $, and it is the only generator of $ G $ that contains $ \gamma_p $, then by Lemma~\ref{lemma 2.5.2} we conclude that the subgroup generated by $ \mathbb{V}(C) $ contains $ \mathcal{C}_p $. Also,
    \begin{align*}
        \mathbb{V}(C) &= (ab)^2(ac) \\&\equiv a_{2}\gt\cdot a_{2}\gq\cdot a_{2}\gt\cdot a_{2}\gq\cdot a_{2}\gt\cdot a_{2} \pmod{\mathcal{C}_p} \\ &= \gt\gq\gt\gq\gt\\ &= \gq^{\widecheck{\tau}+\widecheck{\tau}^2}
    \end{align*}
     which generates $ \mathcal{C}_q $. Therefore, the subgroup generated by $ \mathbb{V}(C) $ is $ G' $. So, Factor Group Lemma~\ref{FGL} applies.
         \end{subcase}
         \end{case}
         \begin{case}
          Assume $ a = a_{2}\gt $ and $ b = \gt\gq $. If $ k \neq 0 $, then by Lemma~\ref{lemma 5.15}(\ref{5.15.1}), $ \langle a,c\rangle = G $ which contradicts the minimality of $ S $. So we can assume $ k = 0 $.
          \begin{subcase}
                 Assume $ i \neq 0 $ and $ j \neq 0 $. Then $ c = a_{2}\gt^j\gamma_{p} $. Thus, by Lemma~\ref{lemma 5.15}(\ref{5.15.2}), $ \langle b,c\rangle = G $ which contradicts the minimality of $ S $.
          \end{subcase}
            \begin{subcase}
                  Assume $ j = 0 $. Then $ i \neq 0 $ and $ c = a_{2}\gamma_{p} $. We have $ \langle \overline{b},\overline{c} \rangle = \langle \overline{a}_3,\overline{a}_2 \rangle = \overline{G} $. Consider $ \{\widecheck{b},\widecheck{c}\} = \{\gt,a_{2}\gamma_{p}\} $. Then we have
                  \begin{align*}
                     [\gt,a_{2}\gamma_{p}] = \gt a_2\gamma_p\gt^{-1}\gamma_p^{-1}a_2 =  \gt\gamma_{p}^{-1}\gt^{-1}\gamma_{p} = \gamma_{p}^{-\widehat{\tau}+1}
                  \end{align*}
                  which generates $ \mathcal{C}_p $. Now consider $ \{\widehat{b},\widehat{c}\} = \{\gt\gq,a_{2}\} $. Thus,
                   \begin{align*}
                       [\gt\gq,a_{2}] = \gt\gq a_{2}\gq^{-1}\gt^{-1}a_{2} = \gt\gq^2\gt^{-1} =  \gq^{2\widecheck{\tau}}
                   \end{align*}
                   which generates $ \mathcal{C}_q $. Therefore, $ \langle b,c \rangle = G $ which contradicts the minimality of $ S $.
          \end{subcase}
          \begin{subcase}
               Assume $ i = 0 $. Then $ j \neq 0 $. We may also assume $ j = 1 $, by replacing $ c $ with $ c^{-1} $ if necessary. Then $ c = \gt\gamma_{p} $. Consider $ \overline{G} = \mathcal{C}_2\times\mathcal{C}_3 $, then we have $ \overline{a} = a_{2}\gt $, $ \overline{b} = \overline{c} = \gt $. Thus, $ |\overline{a}| = 6 $ and $ |\overline{b}| = |\overline{c}| = 3 $. We have $ C = (\overline{c},\overline{b},\overline{a},\overline{b}^{-2},\overline{a}^{-1}) $ as a Hamiltonian cycle in $ \Cay(\overline{G};\overline{S}) $. Since there is one occurrence of $ c $ in $ C $, and it is the only generator of $ G $ that contains $ \gamma_p $, then by Lemma~\ref{lemma 2.5.2} we conclude that the subgroup generated by $ \mathbb{V}(C) $ contains $ \mathcal{C}_p $. Also,
               \begin{align*}
                   \mathbb{V}(C) &= cbab^{-2}a^{-1} \\&\equiv \gt\cdot\gt\gq\cdot a_{2}\gt\cdot\gq^{-1}\gt^{-1}\gq^{-1}\gt^{-1}\cdot\gt^{-1}a_{2} \pmod{\mathcal{C}_p} \\ &= \gt^2\gq\gt\gq\gt^{-1}\gq\gt^{-2}  \\ &= \gq^{\widecheck{\tau}^2+1+\widecheck{\tau}^{-1}} \\&= \gq^{\widecheck{\tau}^2+1-\widecheck{\tau}^2} \\&= \gq
               \end{align*}
               which generates $ \mathcal{C}_q $. Therefore, the subgroup generated by $ \mathbb{V}(C) $ is $ G' $. So, Factor Group Lemma~\ref{FGL} applies.
          \end{subcase}
         \end{case}
         \begin{case}
         Assume $ a = \gt $ and $ b = a_{2}\gq $.
         \begin{subcase}
              Assume $ i = 0 $. Then $ j \neq 0 $. We may also assume $ j = 1 $, by replacing $ c $ with $ c^{-1} $ if necessary. Then $ c = \gt\gq^k\gamma_{p} $. Consider $ \overline{G} = \mathcal{C}_2\times\mathcal{C}_3 $. Then we have $ \overline{a} = \overline{c} = \gt $ and $ \overline{b} = a_2 $. This implies that $ |\overline{a}| = |\overline{c}| = 3 $ and $ |\overline{b}| = 2 $. We have $ C = (\overline{c}^{-2},\overline{b},\overline{a}^2,\overline{b}) $ as a Hamiltonian cycle in $ \Cay(\overline{G};\overline{S}) $. Now we calculate its voltage.
              \begin{align*}
                  \mathbb{V}(C) &= c^{-2}ba^2b \\&\equiv \gamma_p^{-1}\gt^{-1}\gamma_p^{-1}\gt^{-1}\cdot a_2\cdot\gt^2\cdot a_2 \pmod{\mathcal{C}_q} \\&= \gamma_p^{-1}\gt^{-1}\gamma_p^{-1}\gt \\& = \gamma_p^{-1-\widehat{\tau}^{-1}}
              \end{align*}
              which generates $ \mathcal{C}_p $. Also
              \begin{align*}
                  \mathbb{V}(C) &= c^{-2}ba^2b \\&\equiv \gq^{-k}\gt^{-1}\gq^{-k}\gt^{-1}\cdot a_2\gq\cdot\gt^2\cdot a_2\gq \pmod{\mathcal{C}_p} \\&= \gq^{-k}\gt^{-1}\gq^{-k}\gt^{-1}\gq^{-1}\gt^2\gq \\&= \gq^{-k-k\widecheck{\tau}^{-1}-\widecheck{\tau}^{-2}+1}.
              \end{align*}
              If $ k = 2 $, then
              \begin{align*}
         \gq^{-k-k\widecheck{\tau}^{-1}-\widecheck{\tau}^{-2}+1} = \gq^{-2-2\widecheck{\tau}^{-1}-\widecheck{\tau}^{-2}+1} = \gq^{-(\widecheck{\tau}^{-1}+1)^2}
              \end{align*}
              which generates $ \mathcal{C}_q $. So
              we may assume $ k \neq 2 $ and the subgroup generated by $ \mathbb{V}(C) $ does not contain $ \mathcal{C}_q $, for otherwise Factor Group Lemma~\ref{FGL} applies. Therefore,
              \begin{align*}
                  0 &\equiv -k-k\widecheck{\tau}^{-1}-\widecheck{\tau}^{-2}+1 \pmod{q} \\&= (1-k)-k\widecheck{\tau}^{-1}-\widecheck{\tau}^{-2}.
              \end{align*}
              Multiplying by $ \widecheck{\tau}^2 $, we have
              \begin{align*}
                  0 \equiv (1-k)\widecheck{\tau}^2-k\widecheck{\tau}-1 \pmod{q}.\tag{4.1A}\label{3.7.4.1A}
              \end{align*}
              
             We can replace $ \widecheck{\tau} $ with $ \widecheck{\tau}^{-1} $ in the above equation, by replacing $ \gt $,$ a $ and $ c $ with their inverses.
             \begin{align*}
                 0 \equiv (1-k)\widecheck{\tau}^{-2}-k\widecheck{\tau}^{-1}-1 \pmod{q}.
             \end{align*}
             Multiplying by $ \widecheck{\tau}^2 $, then
             \begin{align*}
                 0 \equiv (1-k)-k\widecheck{\tau}-\widecheck{\tau}^2 \pmod{q}.
             \end{align*}
             By subtracting \ref{3.7.4.1A} from the above equation, we have
             \begin{align*}
                 0 \equiv (k-2)\widecheck{\tau}^2+(2-k) \pmod{q}.
             \end{align*}
              This implies that $ \widecheck{\tau}^2 \equiv 1 \pmod{q} $, a contradiction.

         \end{subcase}
         \begin{subcase}
                 Assume $ j = 0 $. Then $ i \neq 0 $. If $ k \neq 0 $, then $ c = a_{2}\gq^k\gamma_{p} $. Thus, by Lemma~\ref{lemma 5.15}(\ref{5.15.3}), $ \langle a,c\rangle = G $ which contradicts the minimality of $ S $. So we can assume $ k = 0 $. Then $ c = a_{2}\gamma_{p} $. Consider $ \overline{G} = \mathcal{C}_2\times\mathcal{C}_3 $, then $ \overline{a} = \gt $ and $ \overline{b} = \overline{c} = a_{2} $. We have $ C = (\overline{a}^2,\overline{b},\overline{a}^{-2},\overline{c}) $ as a Hamiltonian cycle in $ \Cay(\overline{G};\overline{S}) $. Since there is one occurrence of $ c $ in $ C $, and it is the only generator of $ G $ that contains $ \gamma_p $, then by Lemma~\ref{lemma 2.5.2} we conclude that the subgroup generated by $ \mathbb{V}(C) $ contains $ \mathcal{C}_p $. Similarly, since there is one occurrence of $ b $ in $ C $, and it is the only generator of $ G $ that contains $ \gq $, then by Lemma~\ref{lemma 2.5.2} we conclude that the subgroup generated by $ \mathbb{V}(C) $ contains $ \mathcal{C}_q $. Therefore, the subgroup generated by $ \mathbb{V}(C) $ is $ G' $. So, Factor Group Lemma~\ref{FGL} applies.
         \end{subcase}
         \begin{subcase}
               Assume $ i \neq 0 $ and $ j \neq 0 $. If $ k \neq 0 $, then $ c = a_{2}\gt^j\gq^k\gamma_{p} $. Thus, by Lemma~\ref{lemma 5.15}(\ref{5.15.3}), $ \langle a,c\rangle = G $ which contradicts the minimality of $ S $. So we can assume $ k = 0 $. We may also assume $ j = 1 $, by replacing $ c $ with $ c^{-1} $ if necessary. Then $ c = a_{2}\gt\gamma_{p} $. We have $ \langle \overline{b},\overline{c} \rangle = \langle \overline{a}_2,\overline{a}_2\overline{a}_3 \rangle = \overline{G} $. Consider $ \{\widehat{b},\widehat{c}\} = \{a_2\gq,a_2\gt\} $.~Then~we~have
               \begin{align*}
                   [a_2\gq,a_2\gt] = a_2\gq a_2\gt\gq^{-1}a_2\gt^{-1}a_2 = \gq^{-1}\gt\gq^{-1}\gt^{-1} = \gq^{-1-\widecheck{\tau}}
               \end{align*}
               which generates $ \mathcal{C}_q $. Now consider $ \{\widecheck{b},\widecheck{c}\} = \{a_2,a_2\gt\gamma_p\} $. Then
               \begin{align*}
                   [a_2,a_2\gt\gamma_p] = a_2 a_2\gt\gamma_p a_2\gamma_p^{-1}\gt^{-1}a_2 = \gt\gamma_p^2\gt^{-1} = \gamma_p^{2\widehat{\tau}}
               \end{align*}
               which generates $ \mathcal{C}_p $. Therefore, $ \langle b,c \rangle = G $ which contradicts the minimality of $ S $.
    \qedhere
         \end{subcase}
           \end{case}
        \end{proof}
        
 \subsection{Assume \texorpdfstring{$ |S| = 3 $ and $ G' = \mathcal{C}_3\times\mathcal{C}_p $}{Lg}}\hfill\label{3.8}
 
 In this subsection we prove the part of Theorem~\ref{theorem1.1} where, $ |S| = 3 $ and $ G' = \mathcal{C}_3\times\mathcal{C}_p $. Recall $ \overline{G} = G/G' $, $ \widehat{G} = G/\mathcal{C}_p $ and $ \overleftrightarrow{G} = G/\mathcal{C}_3 $.

\begin{caseprop}
Assume
\begin{itemize}
    \item $ G = (\mathcal{C}_2\times\mathcal{C}_q) \ltimes (\mathcal{C}_3\times\mathcal{C}_p) $,
    \item $ |S| = 3 $.
\end{itemize}
Then $ \Cay(G;S) $ contains a Hamiltonian cycle.
\end{caseprop}

\begin{proof}

        Let $ S = \{a,b,c\} $. Since $ \mathcal{C}_q $ centralizes $ \mathcal{C}_3 $ and $ Z(G) \cap G' = \{e\} $ (by Proposition~\ref{Hall Theorem}(\ref{Hall Theorem2})), then $ \mathcal{C}_2 $ inverts $ \mathcal{C}_3 $. Now if $ \widehat{S} $ is minimal, then Lemma~\ref{lemma5.10} applies. So we may assume $ \widehat{S} $ is not minimal. Consider 
        \begin{align*}
            \widehat{G} = G/\mathcal{C}_p = (\mathcal{C}_2\times\mathcal{C}_q) \ltimes\mathcal{C}_3.
        \end{align*}
        Choose a 2-element subset $ \{a,b\} $ in $ S $ that generates $ \widehat{G} $. From the minimality of~$ S $~we~see
        \begin{align*}
            \langle a,b \rangle = (\mathcal{C}_2\times\mathcal{C}_q) \ltimes\mathcal{C}_3.
        \end{align*}
        after replacing $ a $ and $ b $ with conjugates. Then the projection of $ (a,b) $ to $ \mathcal{C}_2\times\mathcal{C}_q $ has one of the following forms:
        \begin{itemize}
            \item $(a_2\gq,a_2\gq^m)$, where $ 1 \leq m \leq q-1 $,
            \item $ (a_2\gq,a_2) $,
            \item $ (a_2\gq,\gq^m) $, where $ 1 \leq m \leq q-1 $,
            \item $ (a_2,\gq) $.
        \end{itemize}
         Thus, there are four different possibilities for $ (a,b) $ after assuming, without loss of generality, that $ a \in \mathcal{C}_2\times\mathcal{C}_q $:
        \begin{enumerate}
            \item $ (a_{2}\gq,a_{2}\gq^m\gt) $,
            \item $ (a_{2}\gq,a_{2}\gt) $,
            \item $ (a_{2}\gq,\gq^m\gt) $,
            \item $ (a_{2},\gq\gt) $.
        \end{enumerate}
        Let $ c $ be the third element of $ S $. We may write $ c = a_{2}^i\gq^j\gt^k\gamma_{p} $ with $ 0\leq i \leq 1 $, $ 0\leq j \leq q-1 $ and $ 0\leq k \leq 2 $. Since $ \mathcal{C}_q $ centralizes $ \mathcal{C}_3 $, we may assume $ \mathcal{C}_q $ does not centralize $ \mathcal{C}_p $, for otherwise Lemma~\ref{lemma 5.8} applies. Now we have $ \gq\gamma_p\gq^{-1} = \gamma_p^{\widehat{\tau}} $, where $ \widehat{\tau}^q \equiv 1 \pmod{p} $. We also have $ \widehat{\tau} \not \equiv 1 \pmod{p} $. Since $ \widehat{\tau}^q \equiv 1 \pmod{p} $, this implies
        \begin{align*}
         \widehat{\tau}^{q-1}+\widehat{\tau}^{q-2}+\cdots+1 \equiv 0 \pmod{p}.
        \end{align*}
         Note that this implies $ \widehat{\tau} \not \equiv -1  \pmod{p} $.
         \setcounter{case}{0}
        \begin{case}
           Assume $ a = a_{2}\gq $ and $ b = a_{2}\gq^m\gt $. If $ k \neq 0 $, then by Lemma~\ref{lemma 5.16}(\ref{lemma 5.16.1}) $ \langle a,c \rangle = G $ which contradicts the minimality of $ S $. So we can assume $ k = 0 $. Now if $ i \neq 0 $, then by Lemma~\ref{lemma 5.16}(\ref{lemma 5.16.3}) $ \langle b,c \rangle = G $ which contradicts the minimality of $ S $. Therefore, we may assume $ i = 0 $. Then $ j \neq 0 $ and $ c = \gq^j\gamma_p $.
          
          Consider $ \overline{G} = \mathcal{C}_2\times\mathcal{C}_q $. Then we have $ \overline{a} = a_{2}\gq $, $ \overline{b} = a_{2}\gq^m $ and $ \overline{c} = \gq^j $. We may assume $ m $ is odd by replacing $ b $ with $ b^{-1} $ (and $ m $ with $ q-m $) if necessary. Note that this implies $ \overline{b} = \overline{a}^m $. Also, we have  $ |\overline{a}| = |\overline{b}| = 2q $ and $ |\overline{c}| = q $.
          
          \begin{subcase}
                Assume $ m = 1 $. Then $ \overline{a} = \overline{b} $. We have
          \begin{align*}
              C = (\overline{c}^{q-1},\overline{b},\overline{c}^{-(q-1)},\overline{a}^{-1})
          \end{align*}
          as a Hamiltonian cycle in $ \Cay(\overline{G};\overline{S}) $. Since there is one occurrence of $ b $ in $ C $, and it is the only generator of $ G $ that contains $ \gt $, then by Lemma~\ref{lemma 2.5.2} we conclude that the subgroup generated by $ \mathbb{V}(C) $ contains $ \mathcal{C}_3 $. Now by considering the fact that $ \mathcal{C}_2 $ might centralize $ \mathcal{C}_p $ or not we have
          \begin{align*}
              \mathbb{V}(C) &= c^{q-1}bc^{-(q-1)}a^{-1} \\&\equiv (\gq^j\gamma_p)^{q-1}\cdot a_2\gq\cdot(\gq^j\gamma_p)^{-(q-1)}\cdot \gq^{-1}a_2 \pmod{\mathcal{C}_3} \\&= \gamma_p^{\widehat{\tau}^j+\widehat{\tau}^{2j}+\cdots+\widehat{\tau}^{(q-1)j}}\gq^{(q-1)j} a_2\gq \gq^{-(q-1)j}\gamma_p^{-(\widehat{\tau}^j+\widehat{\tau}^{2j}+\cdots+\widehat{\tau}^{(q-1)j})} \gq^{-1}a_2 \\&= \gamma_p^{\widehat{\tau}^j(1+\widehat{\tau}^j+\cdots+\widehat{\tau}^{(q-2)j})} \gq \gamma_p^{\mp\widehat{\tau}^j(1+\widehat{\tau}^j+\cdots+\widehat{\tau}^{(q-2)j})}\gq^{-1}.
              \end{align*}
              Now if $ \widehat{\tau}^j \not \equiv 1 \pmod{p} $, then
              \begin{align*}
                   \mathbb{V}(C) &= \gamma_p^{\widehat{\tau}^j(1+\widehat{\tau}^j+\cdots+\widehat{\tau}^{(q-2)j})} \gq \gamma_p^{\mp\widehat{\tau}^j(1+\widehat{\tau}^j+\cdots+\widehat{\tau}^{(q-2)j})}\gq^{-1} \\&= \gamma_p^{\widehat{\tau}^j((\widehat{\tau}^{j})^{q-1}-1)/(\widehat{\tau}^j-1)\mp\widehat{\tau}^{j+1}((\widehat{\tau}^j)^{q-1}-1)/(\widehat{\tau}^j-1)} \\&= \gamma_p^{\widehat{\tau}^j((\widehat{\tau}^{-j})-1)/(\widehat{\tau}^j-1)\mp\widehat{\tau}^{j+1}((\widehat{\tau}^{-j})-1)/(\widehat{\tau}^j-1)} \\&= \gamma_p^{(1-\widehat{\tau}^j)(1\mp\widehat{\tau})/(\widehat{\tau}^j-1)} \\&= \gamma_p^{-(1\mp\widehat{\tau})}.
          \end{align*}
          We may assume this does not generate $ \mathcal{C}_p $, for otherwise Factor Group Lemma~\ref{FGL} applies. Therefore, $ \widehat{\tau}^j \equiv 1 \pmod{p} $ or $ \widehat{\tau} \equiv \pm1 \pmod{p} $. The second case is impossible. So we must have $ \widehat{\tau}^j \equiv 1 \pmod{p} $. We also know that $ \widehat{\tau}^q \equiv 1 \pmod{p} $. So $ \widehat{\tau}^d \equiv 1 \pmod{p} $, where $ d = \gcd(j,q) $. Since $ 1 \leq j \leq q-1 $, then $ d = 1 $, which contradicts the fact that $ \widehat{\tau} \not \equiv 1 \pmod{p} $.
          \end{subcase}
         
         \begin{subcase}
             Assume $ m \neq 1 $ and $ j = 2 $. Then $ c = \gq^2\gamma_p $. We have
          \begin{align*}
              C = (\overline{b},\overline{c}^{-(m-1)/2},\overline{a},\overline{c}^{(m-1)/2},\overline{a}^{2q-m-1})
          \end{align*}
          as a Hamiltonian cycle in $ \Cay(\overline{G};\overline{S}) $. Since there is one occurrence of $ b $ in $ C $, and it is the only generator of $ G $ that contains $ \gt $, then by Lemma~\ref{lemma 2.5.2} we conclude that the subgroup generated by $ \mathbb{V}(C) $ contains $ \mathcal{C}_3 $. Considering the fact that $ \mathcal{C}_2 $ might centralize $ \mathcal{C}_p $ or not we have
          \begin{align*}
              \mathbb{V}(C) &= bc^{-(m-1)/2}ac^{(m-1)/2}a^{2q-m-1} \\&\equiv a_2\gq^m\cdot(\gq^2\gamma_p)^{-(m-1)/2}\cdot a_2\gq\cdot(\gq^2\gamma_p)^{(m-1)/2}\cdot \gq^{2q-m-1} \pmod{\mathcal{C}_3} \\&= a_2\gq^m(\gamma_p^{\widehat{\tau}^2+(\widehat{\tau}^2)^2+\cdots+(\widehat{\tau}^2)^{(m-1)/2}}\gq^{(m-1)})^{-1}a_2\gq(\gamma_p^{\widehat{\tau}^2+(\widehat{\tau}^2)^2+\cdots+(\widehat{\tau}^2)^{(m-1)/2}}\gq^{(m-1)})\gq^{-m-1} \\&= a_2\gq^m\gq^{-m+1}\gamma_p^{-\widehat{\tau}^2(1+\widehat{\tau}^2+\cdots+(\widehat{\tau}^{2})^{(m-3)/2})}a_2\gq\gamma_p^{\widehat{\tau}^2(1+\widehat{\tau}^2+\cdots+(\widehat{\tau}^{2})^{(m-3)/2})}\gq^{-2} \\&= \gq\gamma_p^{\pm\widehat{\tau}^2(1+\widehat{\tau}^2+\cdots+\widehat{\tau}^{2})^{(m-3)/2}}\gq\gamma_p^{\widehat{\tau}^2(1+\widehat{\tau}^2+\cdots+\widehat{\tau}^{2})^{(m-3)/2}}\gq^{-2} \\&= \gamma_p^{\pm\widehat{\tau}^3(\widehat{\tau}^{m-1}-1)/(\widehat{\tau}^2-1)+\widehat{\tau}^4(\widehat{\tau}^{m-1}-1)/(\widehat{\tau}^2-1)} \\&= \gamma_p^{\widehat{\tau}^3(\widehat{\tau}^{m-1}-1)(\pm1+\widehat{\tau})/(\widehat{\tau}^2-1)}.
          \end{align*}
          We may assume this does not generate $ \mathcal{C}_p $, for otherwise Factor Group Lemma~\ref{FGL} applies. Therefore, $ \widehat{\tau}^{m-1} \equiv 1 \pmod{p} $. We also know that $ \widehat{\tau}^q \equiv 1 \pmod{p} $. So $ \widehat{\tau}^d \equiv 1 \pmod{p} $, where $ d = \gcd(m-1,q) $. Since $ 2 \leq m \leq q-1 $, then $ d = 1 $, which contradicts the fact that $ \widehat{\tau} \not \equiv 1 \pmod{p} $.
         \end{subcase}
         
         \begin{subcase}
             Assume $ m \neq 1 $ and $ j \neq 2 $. We may also assume $ j $ is an even number, by replacing $ c $ with its inverse and $ j $ with $ q-j $ if necessary. This implies that $ \overline{c} = \overline{a}^j $. We have
                \begin{align*}
                    C = (\overline{b},\overline{c},\overline{a},\overline{c}^{-1},\overline{b}^{-1},\overline{a}^{m-2},\overline{c},\overline{a}^{-(j-3)},\overline{c},\overline{a}^{2q-m-j-2})
                \end{align*}
                 as a Hamiltonian cycle in $ \Cay(\overline{G};\overline{S}) $. Now we calculate its voltage.
                 \begin{align*}
                     \mathbb{V}(C) &= bcac^{-1}b^{-1}a^{m-2}ca^{-(j-3)}ca^{2q-m-j-2} \\&\equiv a_2\gt\cdot a_2\cdot \gt^{-1}a_2\cdot a_2^{m-2}\cdot a_2^{-(j-3)}\cdot a_2^{2q-m-j-2} \pmod{\mathcal{C}_q\ltimes\mathcal{C}_p} \\&= a_2\gt a_2\gt^{-1} \\&= \gt^{-2}
                 \end{align*}
                 which generates $ \mathcal{C}_3 $. Also considering the fact that $ \mathcal{C}_2 $ might centralize $ \mathcal{C}_p $ or not~we~have
                 \begin{align*}
                     \mathbb{V}(C) &= bcac^{-1}b^{-1}a^{m-2}ca^{-(j-3)}ca^{2q-m-j-2} \\&\equiv a_2\gq^m\cdot\gq^j\gamma_p\cdot a_2\gq\cdot\gamma_p^{-1}\gq^{-j}\cdot\gq^{-m}a_2\\&\indent\indent\cdot a_2\gq^{m-2}\cdot\gq^j\gamma_p\cdot\gq^{-j+3}a_2\cdot\gq^j\gamma_p\cdot a_2\gq^{2q-m-j-2} \pmod{\mathcal{C}_3} \\&= \gq^{m+j}\gamma_p^{\pm 1}\gq\gamma_p^{-1}\gq^{-2}\gamma_p\gq^3\gamma_p^{\pm 1}\gq^{-m-j-2}\\&= \gamma_p^{\pm\widehat{\tau}^{m+j}-\widehat{\tau}^{m+j+1}+\widehat{\tau}^{m+j-1}\pm\widehat{\tau}^{m+j+2}} \\&=\gamma_p^{\widehat{\tau}^{m+j-1}(\pm\widehat{\tau}^3-\widehat{\tau}^2\pm\widehat{\tau}+1)}.
                 \end{align*} 
                 So we may assume this does not generate $ \mathcal{C}_p $, for otherwise Factor Group Lemma~\ref{FGL} applies. Then we have
                 \begin{align*}
                    0 &\equiv \pm\widehat{\tau}^3-\widehat{\tau}^2\pm\widehat{\tau}+1 \pmod{p}.
                 \end{align*}
               Let $ t = \widehat{\tau} $ if $ \mathcal{C}_2 $ centralizes $ \mathcal{C}_p $ and $ t = -\widehat{\tau} $ if $ \mathcal{C}_2 $ inverts $ \mathcal{C}_p $. Then
               \begin{align*}
                   0 &\equiv t^3-t^2+t+1 \pmod{p}. \tag{1.3A}\label{12.3.1A}
               \end{align*}
               We can replace $ t $ with $ t^{-1} $ in the above equation after replacing $ \{a,b,c\} $ with their inverses, then
               \begin{align*}
                   0 &\equiv t^{-3}-t^{-2}+t^{-1}+1 \pmod{p}.
               \end{align*}
               Multiplying by $ t^3 $, we have
               \begin{align*}
                   0 &\equiv 1-t+t^2+t^3 \pmod{p} \\&= t^3+t^2-t+1.
               \end{align*}
               By subtracting \ref{12.3.1A} from the above equation, we have
               \begin{align*}
                   0 &\equiv 2t^2-2t \pmod{p} \\&= 2t(t-1)
               \end{align*}
              This implies that $ t \equiv 1 \pmod{p} $ which contradicts the fact that $ \widehat{\tau} \not \equiv \pm1 \pmod{p} $.
         \end{subcase}
        \end{case}
        
        \begin{case}
            Assume $ a = a_{2}\gq $ and $ b = a_{2}\gt $. If $ k \neq 0 $, then by Lemma~\ref{lemma 5.16}(\ref{lemma 5.16.1}) $ \langle a,c \rangle = G $ which contradicts the minimality of $ S $. So we can assume $ k = 0 $.
            
            \begin{subcase}
                Assume $ i = 0 $. Then $ j \neq 0 $ and $ c = \gq^j\gamma_{p} $. We may assume $ j $ is an odd number, by replacing $ c $ with its inverse and $ j $ with $ q-j $ if necessary. Consider $ \overline{G} = \mathcal{C}_2\times\mathcal{C}_q $. Then we have $ \overline{a} = a_{2}\gq $, $ \overline{b} = a_{2} $ and $ \overline{c} = \gq^j $. Also, we have $ |\overline{a}| = 2q $, $ |\overline{b}| = 2 $ and $ |\overline{c}| = q $. Now if $ j \neq 1 $, then we have
                \begin{align*}
                   C = (\overline{c},\overline{a}^{-1},\overline{b},\overline{a}^2,\overline{b},\overline{c}^{-1},\overline{a}^{j-3},\overline{b},\overline{a}^{-(q-4)},\overline{b},\overline{a}^{q-j-2}) 
                \end{align*}
                 as a Hamiltonian cycle in $ \Cay(\overline{G};\overline{S}) $. Now we calculate the voltage of $ C $.
                \begin{align*}
                    \mathbb{V}(C) &= ca^{-1}ba^2bc^{-1}a^{j-3}ba^{-(q-4)}ba^{q-j-2} \\&\equiv a_{2}\cdot a_{2}\gt\cdot a_2^2\cdot a_{2}\gt\cdot a_2^{j-3}\cdot a_{2}\gt\cdot a_2^{-(q-4)}\cdot a_{2}\gt\cdot a_2^{q-j-2}   \pmod{\mathcal{C}_q\ltimes\mathcal{C}_p} \\&= \gt a_2\gt a_2\gt a_2 a_2\gt \\ &= \gt^2
                \end{align*}
                which generates $ \mathcal{C}_3 $. By considering the fact that $ \mathcal{C}_2 $ might centralize $ \mathcal{C}_p $ or not, we have
                \begin{align*}
                    \mathbb{V}(C) &= ca^{-1}ba^2bc^{-1}a^{j-3}ba^{-(q-4)}ba^{q-j-2} \\&\equiv \gq^j\gamma_{p}\cdot\gq^{-1}a_{2}\cdot a_{2}\cdot\gq^2\cdot a_{2}\cdot\gamma_{p}^{-1}\gq^{-j}\cdot\gq^{j-3}\cdot a_{2}\cdot a_{2}\gq^{-q+4}\cdot a_{2}\cdot \gq^{q-j-2} \pmod{\mathcal{C}_3} \\&= \gq^j\gamma_{p}\gq\gamma_{p}^{\mp1}\gq^{-j-1} \\ &= \gamma_p^{\widehat{\tau}^j\mp\widehat{\tau}^{j+1}}\\&= \gamma_{p}^{\widehat{\tau}^j(1\mp\widehat{\tau})}
                \end{align*}
                which generates $ \mathcal{C}_p $. Therefore, the subgroup generated by $ \mathbb{V}(C) $ is $ G' $. Thus, Factor Group Lemma~\ref{FGL} applies.
                
                So we may assume $ j = 1 $, then $ c = \gq\gamma_p $ and $ \overline{c} = \gq $. We have
                \begin{align*}
                    C_1 = ((\overline{b},\overline{c})^{q-1},\overline{b},\overline{a})                \end{align*}
                as a Hamiltonian cycle in $ \Cay(\overline{G};\overline{S}) $. Now we calculate its voltage. 
                \begin{align*}
                    \mathbb{V}(C_1) &= (bc)^{q-1}ba \\&\equiv (a_2\gt)^{q-1}\cdot a_2\gt\cdot a_2 \pmod{\mathcal{C}_q\ltimes\mathcal{C}_p} \\&= \gt^{-1}
                \end{align*}
                which generates $ \mathcal{C}_3 $. If $ \mathcal{C}_2 $ centralizes $ \mathcal{C}_p $, then
                \begin{align*}
                    \mathbb{V}(C_1) &= (bc)^{q-1}ba \\&\equiv (a_2\cdot \gq\gamma_p)^{q-1}\cdot a_2\cdot a_2\gq \pmod{\mathcal{C}_3} \\&= (\gq\gamma_p)^{q-1}\gq \\&= \gamma_p^{\widehat{\tau}+\widehat{\tau}^2+\cdots+\widehat{\tau}^{q-1}} \\&= \gamma_p^{-1}
                \end{align*}
                which generates $ \mathcal{C}_p $. So in this case, the subgroup generated by $ \mathbb{V}(C_1) $ is $ G' $. Thus, Factor Group Lemma~\ref{FGL} applies.
                
                Now if $ \mathcal{C}_2 $ inverts $ \mathcal{C}_p $, then
                \begin{align*}
                    \mathbb{V}(C_1) &= (bc)^{q-1}ba \\&\equiv (a_2\cdot \gq\gamma_p)^{q-1}\cdot a_2\cdot a_2\gq \pmod{\mathcal{C}_3} \\&=  \gamma_p^{-\widehat{\tau}+\widehat{\tau}^2-\cdots-\widehat{\tau}^{q-2}+\widehat{\tau}^{q-1}}.
                \end{align*}
                Since $ \widehat{\tau} \not \equiv -1 \pmod{p} $, then
                \begin{align*}
                    \mathbb{V}(C_1) &= \gamma_p^{-\widehat{\tau}+\widehat{\tau}^2-\cdots-\widehat{\tau}^{q-2}+\widehat{\tau}^{q-1}} \\&= \gamma_p^{(\widehat{\tau}^q+1)/(\widehat{\tau}+1)-1}.
                \end{align*}
                We may assume this does not generate $ \mathcal{C}_p $, for otherwise Factor Group Lemma~\ref{FGL} applies. Therefore, since $ \widehat{\tau}^q \equiv 1 \pmod{p} $, then
                \begin{align*}
                    0 &\equiv (\widehat{\tau}^q+1)/(\widehat{\tau}+1)-1 \pmod{p} \\&= 2/(\widehat{\tau}+1)-1.
                \end{align*}
                This implies that $ \widehat{\tau} \equiv 1 \pmod{p} $, which is impossible.
            \end{subcase}
            
            \begin{subcase} \label{subcase for ch:4}
                Assume $ j  = 0 $. Then $ i \neq 0 $ and $ c = a_{2}\gamma_{p} $. Consider $ \overline{G} = \mathcal{C}_2\times\mathcal{C}_q $. Then we have $ \overline{a} = a_{2}\gq $ and $ \overline{b} = \overline{c} = a_{2} $. This implies that $ |\overline{a}| = 2q $ and $ |\overline{b}| = |\overline{c}| = 2 $.~We~have
                \begin{align*}
                    C = (\overline{c},\overline{a}^{q-1},\overline{b},\overline{a}^{-(q-1)})
                \end{align*}
                 as a Hamiltonian cycle in $ \Cay(\overline{G};\overline{S}) $. Since there is one occurrence of $ b $ in $ C $, and it is the only generator of $ G $ that contains $ \gt $, then by Lemma~\ref{lemma 2.5.2} we conclude that the subgroup generated by $ \mathbb{V}(C) $ contains $ \mathcal{C}_3 $. Similarly, since there is one occurrence of $ c $ in $ C $, and it is the only generator of $ G $ that contains $ \gamma_p $, then by Lemma~\ref{lemma 2.5.2} we conclude that the subgroup generated by $ \mathbb{V}(C) $ contains $ \mathcal{C}_p $. Therefore, the subgroup generated by $ \mathbb{V}(C) $ is $ G' $. So, Factor Group Lemma~\ref{FGL} applies. 
            \end{subcase}
            
            \begin{subcase}
                Assume $ i \neq 0 $ and $ j \neq 0 $. Then $ c = a_2\gq^j\gamma_p $. Consider $ \overline{G} = \mathcal{C}_2\times\mathcal{C}_q $. Then we have $ \overline{a} = a_{2}\gq $, $ \overline{b} = a_{2} $ and $ \overline{c} = a_{2}\gq^j $. This implies that $ |\overline{a}| = |\overline{c}| = 2q $ and $ |\overline{b}| = 2 $. We may assume $ j $ is even by replacing $ c $ with its inverse and $ j $ with $ q-j $ if necessary. 
                
                Suppose, for the moment, that $ j = q-1 $, then $ c = a_2\gq^{-1}\gamma_p $ and $ \overline{c} = \overline{a}^{-1} $. We have 
                \begin{align*}
                    C_1 = (\overline{c},\overline{b},(\overline{a}^{-1},\overline{b})^{q-1})
                \end{align*}
                as a Hamiltonian cycle in $ \Cay(\overline{G};\overline{S}) $. Since there is one occurrence of $ c $ in $ C $, and it is the only generator of $ G $ that contains $ \gamma_p $, then by Lemma~\ref{lemma 2.5.2} we conclude that the subgroup generated by $ \mathbb{V}(C_1) $ contains $ \mathcal{C}_p $. Also,
                \begin{align*}
                    \mathbb{V}(C_1) &= cb(a^{-1}b)^{q-1} \\&\equiv a_2\cdot a_2\gt\cdot (a_2\cdot a_2\gt)^{q-1} \pmod{\mathcal{C}_q\ltimes\mathcal{C}_p} \\&= \gt^q
                \end{align*}
                which generates $ \mathcal{C}_3 $. Therefore, the subgroup generated by $ \mathbb{V}(C_1) $ contains $ G' $. Thus, Factor Group Lemma~\ref{FGL} applies.
                
                So we may assume $ j \neq q-1 $. Then we have
                \begin{align*}
                    C_{2} = (\overline{c},\overline{a}^{q-j-1},\overline{b},\overline{a}^{-q+j+1},(\overline{a}^{-1},\overline{b})^{j})
                \end{align*}
                 and
                 \begin{align*}
                   C_{3} = (\overline{c},\overline{a}^{q-j-2},\overline{b},\overline{a}^{-q+j+2},(\overline{a}^{-1},\overline{b})^{j-1},\overline{a}^{-2},\overline{b},\overline{a})  
                 \end{align*}
                  as Hamiltonian cycles in $ \Cay(\overline{G};\overline{S}) $. Since there is one occurrence of $ c $ in $ C_2 $, and it is the only generator of $ G $ that contains $ \gamma_p $, then by Lemma~\ref{lemma 2.5.2} we conclude that the subgroup generated by $ \mathbb{V}(C_2) $ contains $ \mathcal{C}_p $. Also,
                \begin{align*}
                    \mathbb{V}(C_2) &= ca^{q-j-1}ba^{-q+j+1}(a^{-1}b)^j \\&\equiv a_2\cdot a_{2}^{q-j-1}\cdot a_{2}\gt\cdot a_2^{-q+j+1}\cdot\gt^j \pmod{\mathcal{C}_q\ltimes\mathcal{C}_p} \\ &= \gt^{j+1}.
                \end{align*}
                We may assume this does not generate $ \mathcal{C}_3 $, for otherwise Factor Group Lemma~\ref{FGL} applies. Then $ j \equiv -1 \pmod{3} $. 
                
                Since there is one occurrence of $ c $ in $ C_3 $, and it is the only generator of $ G $ that contains $ \gamma_p $, then by Lemma~\ref{lemma 2.5.2} we conclude that the subgroup generated by $ \mathbb{V}(C_3) $ contains $ \mathcal{C}_p $. Also,
                 \begin{align*}
                    \mathbb{V}(C_3) &= ca^{q-j-2}ba^{-q+j+2}(a^{-1}b)^{j-1}a^{-2}ba \\&\equiv a_{2}\cdot a_{2}^{q-j-2}\cdot a_{2}\gt\cdot a_{2}^{-q+j+2}\cdot\gt^{j-1}\cdot a_2^{-2}\cdot a_{2}\gt\cdot a_{2} \pmod{\mathcal{C}_q\ltimes\mathcal{C}_p} \\&= a_2\gt a_2\gt^{j-1}a_2\gt a_2 \\ &= \gt^{j-3}\\&= \gt^j
                \end{align*}
                Since $ j \equiv -1 \pmod{3} $, this generates $ \mathcal{C}_3 $. So, Factor Group Lemma~\ref{FGL} applies.
            \end{subcase}
        \end{case}
        
        \begin{case}
            Assume $ a = a_{2}\gq $ and $ b = \gq^m\gt $. If $ k \neq 0 $, then by Lemma~\ref{lemma 5.16}(\ref{lemma 5.16.1}) $ \langle a,c \rangle = G $ which contradicts the minimality of $ S $. So we can assume $ k = 0 $. Now if $ i \neq 0 $, then by Lemma~\ref{lemma 5.16}(\ref{lemma 5.16.3}) $ \langle b,c \rangle = G $ which contradicts the minimality of $ S $. Therefore, we may assume $ i = 0 $. Then $ j \neq 0 $ and $ c = \gq^j\gamma_p $. 
            Consider $ \overline{G} = \mathcal{C}_2\times\mathcal{C}_q $. Then we have $ \overline{a} = a_2\gq $, $ \overline{b} = \gq^m $ and $ \overline{c} = \gq^j $.
            
            Suppose, for the moment, that $ m = j $. Then $ \overline{b} = \overline{c} $. We have
            \begin{align*}
               C_1 = (\overline{c}^{-1},\overline{b}^{-(q-2)},\overline{a}^{-1},\overline{b}^{q-1},\overline{a}) 
            \end{align*}
            as a Hamiltonian cycle in $ \Cay(\overline{G};\overline{S}) $. Since there is one occurrence of $ c $ in $ C_1 $, and it is the only generator of $ G $ that contains $ \gamma_p $, then by Lemma~\ref{lemma 2.5.2} we conclude that the subgroup generated by $ \mathbb{V}(C_1) $ contains $ \mathcal{C}_p $. Also,
            \begin{align*}
                \mathbb{V}(C_1) &= c^{-1}b^{-(q-2)}a^{-1}b^{q-1}a \\&\equiv \gt^{-(q-2)}\cdot a_2\cdot\gt^{q-1}\cdot a_2 \pmod{\mathcal{C}_q\ltimes\mathcal{C}_p}\\&= \gt^{-2q+3} \\&= \gt^{-2q}
            \end{align*}
            which generates $ \mathcal{C}_3 $, because $ \gcd(2q,3) = 1 $. So, the subgroup generated by $ \mathbb{V}(C_1) $ is $ G' $. Therefore, Factor Group Lemma~\ref{FGL} applies.
           
           So we may assume $ m \neq j $. We may also assume $ m $ and $ j $ are even, by replacing $ \{b,c\} $ with their inverses, $ m $ with $ q-m $, and $ j $ with $ q-j $ if necessary. Now suppose, for the moment, $ j = 2 $. Then we have $ c = \gq^2\gamma_p $. We also have
           \begin{align*}
               C_2 = (\overline{b},\overline{c}^{-(m-2)/2},\overline{a}^{-1},\overline{c}^{m/2},\overline{a}^{2q-m-1})
           \end{align*}
           as a Hamiltonian cycle in $ \Cay(\overline{G};\overline{S}) $. Since there is one occurrence of $ b $ in $ C_2 $, and it is the only generator of $ G $ that contains $ \gt $, then by Lemma~\ref{lemma 2.5.2} we conclude that the subgroup generated by $ \mathbb{V}(C_2) $ contains $ \mathcal{C}_3 $. Now by considering the fact that $ \mathcal{C}_2 $ might centralize $ \mathcal{C}_p $ or not, we have
           \begin{align*}
               \mathbb{V}(C_2) &= bc^{-(m-2)/2}a^{-1}c^{m/2}a^{2q-m-1} \\&\equiv \gq^m\cdot(\gq^2\gamma_p)^{-(m-2)/2}\cdot \gq^{-1}a_2\cdot(\gq^2\gamma_p)^{m/2}\cdot a_2^{2q-m-1}\gq^{2q-m-1} \pmod{\mathcal{C}_3} \\&= \gq^m(\gamma_p^{\widehat{\tau}^2+(\widehat{\tau}^2)^2+\cdots+(\widehat{\tau}^2)^{(m-2)/2}}\gq^{(m-2)})^{-1}\gq^{-1}a_2(\gamma_p^{\widehat{\tau}^2+(\widehat{\tau}^2)^2+\cdots+(\widehat{\tau}^2)^{m/2}}\gq^m)a_2\gq^{-m-1} \\&= \gq^m\gq^{-(m-2)}\gamma_p^{-\widehat{\tau}^2(1+\widehat{\tau}^2+\cdots+(\widehat{\tau}^2)^{(m-4)/2})}\gq^{-1}\gamma_p^{\pm\widehat{\tau}^2(1+\widehat{\tau}^2+\cdots+(\widehat{\tau}^2)^{(m-2)/2})}\gq^{m}\gq^{-m-1}.
           \end{align*}
           Since $ \widehat{\tau}^2 -1 \not \equiv 0 \pmod{p} $, then
           \begin{align*}
               \mathbb{V}(C_2) &= \gq^2\gamma_p^{-\widehat{\tau}^2(\widehat{\tau}^{m-2}-1)/(\widehat{\tau}^2-1)}\gq^{-1}\gamma_p^{\pm\widehat{\tau}^2(\widehat{\tau}^m-1)/(\widehat{\tau}^2-1)}\gq^{-1} \\&= \gamma_p^{-\widehat{\tau}^4(\widehat{\tau}^{m-2}-1)/(\widehat{\tau}^2-1)\pm\widehat{\tau}^3(\widehat{\tau}^m-1)/(\widehat{\tau}^2-1)} \\&= \gamma_p^{\widehat{\tau}^3(1\mp\widehat{\tau})(-\widehat{\tau}^{m-1}\mp1)/(\widehat{\tau}^2-1)}.
           \end{align*}
           We may assume this does not generate $ \mathcal{C}_p $, for otherwise Factor Group Lemma~\ref{FGL} applies. Therefore, $ \widehat{\tau} \equiv \pm1 \pmod{p} $ or $ \widehat{\tau}^{m-1} \equiv \pm1 \pmod{p} $. The first case is impossible. So we may assume $ \widehat{\tau}^{m-1} \equiv \pm1 \pmod{p} $. Thus, $ \widehat{\tau}^{2(m-1)} \equiv 1 \pmod{p} $. We also know that $ \widehat{\tau}^q \equiv 1 \pmod{p} $. So we have $ \widehat{\tau}^d \equiv 1 \pmod{p} $, where $ d = \gcd(2(m-1),q) $. Since $ \gcd(2,q) = 1 $ and $ 2 \leq m \leq q-1 $, then $ d = 1 $, which contradicts the fact that $ \widehat{\tau} \not \equiv 1 \pmod{p} $.
           
           So we may assume $ j \neq 2 $. We have 
           \begin{align*}
               C_3 = (\overline{b},\overline{c},\overline{a},\overline{c}^{-1},\overline{b}^{-1},\overline{a}^{m-2},\overline{c},\overline{a}^{-(j-3)},\overline{c},\overline{a}^{2q-m-j-2})
           \end{align*}
           as a Hamiltonian cycle in $ \Cay(\overline{G};\overline{S}) $. Now we calculate its voltage.
           \begin{align*}
               \mathbb{V}(C_3) &= bcac^{-1}b^{-1}a^{m-2}ca^{-(j-3)}ca^{2q-m-j-2} \\&\equiv \gt\cdot a_2\cdot\gt^{-1}\cdot a_2^{m-2}\cdot a_2^{-j+3}\cdot a_2^{2q-m-j-2} \pmod{\mathcal{C}_q\ltimes\mathcal{C}_p}\\&= \gt^2
           \end{align*}
           which generates $ \mathcal{C}_3 $. Also, by considering the fact that $ \mathcal{C}_2 $ might centralize $ \mathcal{C}_p $ or not, we have
           \begin{align*}
               \mathbb{V}(C_3) &= bcac^{-1}b^{-1}a^{m-2}ca^{-(j-3)}ca^{2q-m-j-2} \\&\equiv \gq^m\cdot\gq^j\gamma_p\cdot a_2\gq\cdot\gamma_p^{-1}\gq^{-j}\cdot\gq^{-m}\cdot a_2^{m-2}\gq^{m-2}\\&\indent\indent\cdot\gq^j\gamma_p\cdot \gq^{-j+3}a_2^{-j+3}\cdot\gq^j\gamma_p\cdot a_2^{2q-m-j-2}\gq^{2q-m-j-2} \pmod{\mathcal{C}_3} \\&= \gq^{m+j}\gamma_p a_2\gq\gamma_p^{-1}\gq^{-2}\gamma_p\gq^3 a_2\gamma_p\gq^{-m-j-2} \\&= \gq^{m+j}\gamma_p\gq\gamma_p^{\mp1}\gq^{-2}\gamma_p^{\pm1}\gq^3\gamma_p\gq^{-m-j-2}\\&= \gamma_p^{\widehat{\tau}^{m+j}\mp\widehat{\tau}^{m+j+1}\pm\widehat{\tau}^{m+j-1}+\widehat{\tau}^{m+j+2}}\\&= \gamma_p^{\widehat{\tau}^{m+j-1}(\widehat{\tau}^3\mp\widehat{\tau}^2+\widehat{\tau}\pm1)}.
           \end{align*}
           We may assume this does not generate $ \mathcal{C}_p $, for otherwise Factor Group Lemma~\ref{FGL} applies. Therefore,
           \begin{align*}
               0 \equiv \widehat{\tau}^3\mp\widehat{\tau}^2+\widehat{\tau}\pm1 \pmod{p}.
           \end{align*}
           If $ \mathcal{C}_2 $ centralizes $ \mathcal{C}_p $, then
           \begin{align*}
               0 \equiv \widehat{\tau}^3-\widehat{\tau}^2+\widehat{\tau}+1 \pmod{p}. \tag{3A}\label{3.8.3A}
           \end{align*}
           We can replace $ \widehat{\tau} $ with $ \widehat{\tau}^{-1} $ in the above equation after replacing $ \{a,b,c\} $ with their inverses in the Hamiltonian cycle, then
           \begin{align*}
               0 \equiv \widehat{\tau}^{-3}-\widehat{\tau}^{-2}+\widehat{\tau}^{-1}+1 \pmod{p}.
           \end{align*}
           Multiplying by $ \widehat{\tau}^3 $, we have
           \begin{align*}
               0 &\equiv 1-\widehat{\tau}+\widehat{\tau}^2+\widehat{\tau}^3 \pmod{p} \\&= \widehat{\tau}^3+\widehat{\tau}^2-\widehat{\tau}+1.
           \end{align*}
           Subtracting \ref{3.8.3A} from the above equation we have
           \begin{align*}
               0 &\equiv 2\widehat{\tau}^2-2\widehat{\tau} \pmod{p} \\&= 2\widehat{\tau}(\widehat{\tau}-1)
           \end{align*}
           which is impossible, because $ \widehat{\tau} \not \equiv 1 \pmod{p} $.
           
           Now if $ \mathcal{C}_2 $ inverts $ \mathcal{C}_p $, then
           \begin{align*}
               0 \equiv \widehat{\tau}^3+\widehat{\tau}^2+\widehat{\tau}-1 \pmod{p}. \tag{3B}\label{3.8.3B}
           \end{align*}
           We can replace $ \widehat{\tau} $ with $ \widehat{\tau}^{-1} $ in the above equation after replacing $ \{a,b,c\} $ with their inverses. Then
           \begin{align*}
               0 \equiv \widehat{\tau}^{-3}+\widehat{\tau}^{-2}+\widehat{\tau}^{-1}-1 \pmod{p}.
           \end{align*}
           Multiplying by $ \widehat{\tau}^3 $, then
           \begin{align*}
               0 &\equiv 1+\widehat{\tau}+\widehat{\tau}^2-\widehat{\tau}^3 \pmod{p} \\&= -\widehat{\tau}^3+\widehat{\tau}^2+\widehat{\tau}+1.
           \end{align*}
           By adding \ref{3.8.3B} and the above equation, we have
           \begin{align*}
               0 &\equiv 2(\widehat{\tau}^2+\widehat{\tau}) \pmod{p} \\&= 2\widehat{\tau}(\widehat{\tau}+1) 
           \end{align*}
           which is also impossible, because $ \widehat{\tau} \not \equiv -1 \pmod{p} $.
        \end{case}
        \begin{case}
            Assume $ a = a_{2} $ and $ b = \gq\gt $.
            
            \begin{subcase}
                Assume $ i \neq 0 $. Then $ c = a_{2}\gq^j\gt^k\gamma_{p} $. By Lemma~\ref{lemma 5.16}(\ref{lemma 5.16.2}) $ \langle b,c \rangle = G $ which contradicts the minimality of $ S $.
            \end{subcase}
            
            \begin{subcase}
                Assume $ i = 0 $. Then $ j \neq 0 $ and $ c = \gq^j\gt^k\gamma_{p} $. 
                 We may assume $ j $ is even by replacing $ c $ with its inverse and $ j $ with $ q-j $ if necessary. Consider $ \overline{G} = \mathcal{C}_2\times\mathcal{C}_q $. Then we have $ \overline{a} = a_{2} $, $ \overline{b} = \gq $ and $ \overline{c} = \gq^j $. This implies that $ |\overline{a}| = 2 $ and $ |\overline{b}| = |\overline{c}| = q $. We have 
                 \begin{align*}
                     C_1 = (\overline{c},\overline{b}^{q-j-1},\overline{c},\overline{b}^{-(j-2)},\overline{a},\overline{b}^{q-1},\overline{a})
                 \end{align*}
                 as a Hamiltonian cycle in $ \Cay(\overline{G};\overline{S}) $. Now we calculate its voltage.
                 \begin{align*}
                     \mathbb{V}(C_1) &= cb^{q-j-1}cb^{-(j-2)}ab^{q-1}a \\&\equiv \gq^j\gamma_p\cdot\gq^{q-j-1}\cdot\gq^j\gamma_p\cdot\gq^{-j+2}\cdot a_2\cdot\gq^{q-1}\cdot a_2 \pmod{\mathcal{C}_3} \\&= \gq^j\gamma_p\gq^{-1}\gamma_p\gq^{-j+1} \\&= \gamma_p^{\widehat{\tau}^{j-1}(\widehat{\tau}+1)}
                 \end{align*}
                 which generates $ \mathcal{C}_p $. Also
                 \begin{align*}
                     \mathbb{V}(C_1) &= cb^{q-j-1}cb^{-(j-2)}ab^{q-1}a \\&\equiv \gt^k\cdot\gt^{q-j-1}\cdot\gt^k\cdot\gt^{-j+2}\cdot a_2\cdot\gt^{q-1}\cdot a_2 \pmod{\mathcal{C}_q\ltimes\mathcal{C}_p} \\&= \gt^{k+q-j-1+k-j+2-q+1} \\&= \gt^{2(k-j+1)}.
                 \end{align*}
                 We may assume this does not generate $ \mathcal{C}_3 $, for otherwise Factor Group Lemma~\ref{FGL} applies. Then
                 \begin{align*}
                     0 \equiv k-j+1 \pmod{3}. \tag{4.2A} \label{12.4.2B}
                 \end{align*}
                 
                 We also have
                \begin{align*}
                   C_{2} = (\overline{c},\overline{a},(\overline{b},\overline{a})^{q-j-1},\overline{b}^j,\overline{a},\overline{b}^{-(j-1)}) 
                \end{align*}
                 as a Hamiltonian cycle in $ \Cay(\overline{G};\overline{S}) $. We calculate its voltage. Since there is one occurrence of $ c $ in $ C_2 $, and it is the only generator of $ G $ that contains $ \gamma_p $, then by Lemma~\ref{lemma 2.5.2} we conclude that the subgroup generated by $ \mathbb{V}(C_2) $ contains $ \mathcal{C}_p $. Also,
                \begin{align*}
                    \mathbb{V}(C_2) &= ca(ba)^{q-j-1}b^jab^{-(j-1)} \\&\equiv \gt^k\cdot a_{2}\cdot(\gt a_2)^{q-j-1}\cdot\gt^j\cdot a_{2}\cdot\gt^{-j+1} \pmod{\mathcal{C}_q\ltimes\mathcal{C}_p} \\
                    &= \gt^{k-2j+1}.
                \end{align*}
                We may assume this does not generate $ \mathcal{C}_3 $, for otherwise Factor Group Lemma~\ref{FGL} applies. Therefore,
                \begin{align*}
                    0 \equiv k-2j+1 \pmod{3}. 
                \end{align*}
             By subtracting the above equation from \ref{12.4.2B} we have $ j \equiv 0 \pmod{3} $.
             
              Now we have
              \begin{align*}
                 C_{3} = (\overline{c},\overline{a},\overline{b}^{q-j-1},\overline{a},\overline{b}^{-(q-j-2)},\overline{c}^{-1},\overline{b}^{j-2},\overline{a},\overline{b}^{-(j-1)},\overline{a}) 
              \end{align*}
              as a Hamiltonian cycle in $ \Cay(\overline{G};\overline{S}) $. We calculate its voltage.
              \begin{align*}
                \mathbb{V}(C_3) &=cab^{q-j-1}ab^{-(q-j-2)}c^{-1}b^{j-2}ab^{-(j-1)}a \\&\equiv \gq^j\gamma_{p}\cdot a_{2}\cdot\gq^{q-j-1}\cdot a_{2}\cdot\gq^{-q+j+2}\cdot\gamma_{p}^{-1}\gq^{-j}\cdot\gq^{j-2}\cdot a_{2}\cdot\gq^{-j+1}\cdot a_{2} \pmod{\mathcal{C}_3} \\ &= \gq^j\gamma_{p}\gq\gamma_{p}^{-1}\gq^{-j-1} \\ &= \gamma_{p}^{\widehat{\tau}^j(1-\widehat{\tau})}. 
            \end{align*}
            which generates $ \mathcal{C}_p $. Also
            \begin{align*}
            \mathbb{V}(C_3) &= cab^{q-j-1}ab^{-(q-j-2)}c^{-1}b^{j-2}ab^{-(j-1)}a \\&\equiv \gt^k\cdot a_{2}\cdot\gt^{q-j-1}\cdot a_{2}\cdot\gt^{-q+j+2}\cdot\gt^{-k}\cdot\gt^{j-2}\cdot a_{2}\cdot\gt^{-j+1}\cdot a_{2} \pmod{\mathcal{C}_q\ltimes\mathcal{C}_p} \\ &= \gt^{k-q+j+1-q+j+2-k+j-2+j-1}  \\ &= \gt^{-2q+4j}.
            \end{align*}
            We may assume this does not generate $ \mathcal{C}_3 $, for otherwise Factor Group Lemma~\ref{FGL} applies. Then
            \begin{align*}
                0 &\equiv -2q+4j \pmod{3} \\&= q+j
            \end{align*}
            We already know $ j \equiv 0 \pmod{3} $. By substituting this in the above equation, we have $ q \equiv 0 \pmod{3} $ which contradicts the fact that $ \gcd(q,3) = 1 $.
            \qedhere
            \end{subcase}
        \end{case}
\end{proof}

\subsection{Assume \texorpdfstring{$ |S| \geq 4 $}{Lg}}\hfill\label{3.9}

In this subsection we prove the following general result that includes the part of Theorem~\ref{theorem1.1}, where $ |S| \geq 4 $ (see Assumption~\ref{assumption 3.1}). Unlike in the other subsections of this section, we do not assume $ |G| = 6pq $.

\begin{caseprop} \label{prop 3.9}
Assume $ |G| $ is a product of four distinct primes and $ S $ is a minimal generating set of $ G $, where $ |S| \geq 4 $. Then $ \Cay(G;S) $ contains a Hamiltonian cycle.
\end{caseprop}

    \begin{proof}
    Suppose $ S = \{s_{1}, s_{2},..., s_{k}\} $ and let $ G_{i} = \langle s_{1}, s_{2},..., s_{i} \rangle $ for $ i = 1, 2,..., k $. Since $ S $ is minimal, we know $ \{e\}\subset G_{1}\subset G_{2}\subset...G_{k} = G $. Therefore, the number of prime factors of $ |G_{i}| $ is at least $ i $. Since $ |G| = p_1p_2p_3q $ is the product of only $ 4 $ primes, and $ k = |S|\geq4 $, we can conclude that $ |G_{i}| $ has exactly $ i $ prime factors, for all $ i $. This implies that $ |S| = 4 $. This also implies every element of $S$ has prime order.
    
   Since $ |G| $ is square-free, we know that $ G' $ is cyclic (see Proposition~\ref{Hall Theorem}(\ref{Hall Theorem1})), so $ G' \neq G $. We may assume $ |G'| \neq 1 $, for otherwise $ G $ is abelian, so Lemma~\ref{abelain group} applies. Also, if $ |G'| $ is equal to a prime number, then Theorem~\ref{Keating-Witte} applies. So we may assume $ |G'| $ has at least two prime factors. Therefore, the number of prime factors of $ |G'| $ is either~$ 2 $~or~$ 3 $.
    \setcounter{case}{0}
    \begin{case}
        Assume $ |G'| $ has only two prime factors. This implies $ |\overline{G}| = p_1p_2 $, where $ p_1 $ and $ p_2 $ are two distinct primes. Suppose $ s \in S $, then $ \overline{s} \in \overline{S} $. We know that $ |\overline{s}| \neq 1 $ (see Assumption~\ref{assumption 3.1}(\ref{assumption 3.1.6})). Now since every element of $ S $ has prime order, then $ |s| $ is either $ p_1 $ or $ p_2 $. Also, every element of order $ p_1 $ must commute with every element of order $ p_2 $, because the subgroup $ H $ generated by any element of $ S $ that has order $ p_1 $, together with any element of $ S $ that has order $ p_2 $ has exactly two prime factors, so $ |H| = p_1p_2 $, $ H' \subseteq G' $, and $ |G'| = p_3p_4 $. Thus, $ |H'| = 1 $. Let $ S_{p_1} $ be the elements of order $ p_1 $ in $ S $,  and let $ S_{p_2} $ be the elements of order $ p_2 $. Also let $ H_{p_1} $ and $ H_{p_2} $ be the subgroups generated by $ S_{p_1} $ and $ S_{p_2} $, respectively. This implies that $ \Cay(G;S) \cong \Cay(G_{p_1};S_{p_1}) \cartprod  \Cay(G_{p_2};S_{p_2}) $. Therefore, $ \Cay(G;S) $ contains a Hamiltonian cycle~(see Corollary~\ref{Cartesian}).
    \end{case}
    
    \begin{case}
        Assume $ |G'| $ has three prime factors. We may write (see Proposition~\ref{Hall Theorem}(\ref{Hall Theorem3}))
        \begin{align*}
           G = \mathcal{C}_q \ltimes G' = \mathcal{C}_q \ltimes (\mathcal{C}_{p_1}\times\mathcal{C}_{p_2}\times\mathcal{C}_{p_3}),
        \end{align*}
        where $ p_1 $, $ p_2 $, $ p_3 $ and $ q $ are distinct primes. Note that $ G' \cap Z(G) = \{e\} $ (see Proposition~\ref{Hall Theorem}(\ref{Hall Theorem2})). Now we may assume $ \langle s_4 \rangle = \mathcal{C}_q $. Since $ |\langle s_i,s_4 \rangle| $ has only two prime factors (for $1 \leq i \leq 3$), we must have $ s_i = s_{4}^{k_i}a_{p_i} $ (after permuting $ p_1 $, $ p_2 $, $ p_3 $), where $ a_{p_i} $ is a generator of $ \mathcal{C}_{p_i} $. We may also assume $ S \cap G' = \emptyset $ (see Lemma~\ref{lemma 5.6}), so $ k_i \not \equiv 0 \pmod{q} $. Now consider
        \begin{align*}
           G_2 = \langle s_1,s_2 \rangle =  \langle s_{4}^{k_1}a_{p_1},s_{4}^{k_2}a_{p_2} \rangle.
        \end{align*}
        Since $ \mathcal{C}_{p_1} $ is a normal subgroup in $ G $, we can consider $ \overline{G}_2 = G_{2}/\mathcal{C}_{p_1} $, then $ \{\overline{s}_1,\overline{s}_2\} = \{\overline{s}_4^{k_1},\overline{s}_4^{k_2}\overline{a}_{p_2}\} $. We have
        \begin{align*}
            \overline{s}_4^{k_2^{-1}} = (\overline{s}_4^{k_1})^{k_1^{-1}k_2^{-1}} = \overline{s}_1^{k_1^{-1}k_2^{-1}}.
        \end{align*}
        Multiplying by  $ \overline{s}_2 $, then
        \begin{align*}
           \overline{a}_{p_2} = \overline{s}_4^{k_{2}^{-1}}\cdot\overline{s}_4^{k_2}a_{p_2} = \overline{s}_1^{k_1^{-1}k_2^{-1}}\overline{s}_2 \in \overline{G}_2.
        \end{align*}
        Since $a_{p_2}$ generates $\mathcal{C}_{p_2}$, this implies $|G_2|$ is divisible by $p_2$. Similarly, we can show that $|G_2|$ is divisible by $p_1$. Also, $|s_1| = q$, so $|G_2|$ is divisible by $q$. Therefore, $ |G_2| $ has three prime factors, which is a contradiction. 
       \qedhere
    \end{case}
    \end{proof}
    
\textbf{Acknowledgements}. Theorem~\ref{theorem1.1} and Proposition~\ref{pro 1.1.4} are the main results of my masters thesis (University of Lethbridge, 2020). I would like to express my sincere gratitude to my supervisor, professor Joy Morris who always supported me throughout my graduate journey. I am especially grateful to my co-supervisor, professor Dave Morris, for the patient guidance and advice he has provided during my graduate study. I have been extremely lucky to have a co-supervisor who cared so much about my research, and who responded to my questions so promptly.
I am also thankful to professor Hadi Kharaghani and professor Amir Akbary and cannot forget their valuable help and motivation during my graduate years.
I am truly grateful to my family for their immeasurable love and care.

\bibliographystyle{abbrv}
\bibliography{paper}
\end{document}